\documentclass[reqno]{amsart}

\newtheorem{theorem}{Theorem}[section]
\newtheorem{remark}[theorem]{Remark}
\newtheorem{lemma}[theorem]{Lemma}
\newtheorem{corollary}[theorem]{Corollary}

\newtheorem{definition}{Definition}[section]

\RequirePackage{fix-cm}

 \usepackage{mathptmx}      
\usepackage{latexsym}
\usepackage{amssymb}

\newcommand{\nat}{\mathbb{N}}
\newcommand{\rzecz}{\mathbb{R}}

\newcommand{\eps}{\varepsilon }

\newcommand{\rd}{{\rzecz }^{d}}

\newcommand{\lin}{\mbox{\rm span}}

\newcommand{\tr}{\mbox{\rm Tr}}

\newcommand{\diver}{\mbox{\rm div}}
\newcommand{\supp}{\mbox{\rm supp}\, }

\newcommand{\acal}{\mathcal{A}}
\newcommand{\bcal}{\mathcal{B}}
\newcommand{\ccal}{\mathcal{C}}

\newcommand{\fcal}{\mathcal{F}}

\newcommand{\hcal}{\mathcal{H}}

\newcommand{\kcal}{\mathcal{K}}
\newcommand{\lcal}{\mathcal{L}}
\newcommand{\mcal}{\mathcal{M}}

\newcommand{\ocal}{\mathcal{O}}
\newcommand{\pcal}{\mathcal{P}}

\newcommand{\scal}{\mathcal{S}}
\newcommand{\tcal}{\mathcal{T}}
\newcommand{\vcal}{\mathcal{V}}

\newcommand{\xcal}{\mathcal{X}}
\newcommand{\ycal}{\mathcal{Y}}
\newcommand{\zcal}{\mathcal{Z}}

\newcommand{\cmath}{\mathbb{C}}
\newcommand{\dmath}{\mathbb{D}}
\newcommand{\fmath}{\mathbb{F}}

\newcommand{\hmath}{\mathbb{H}}

\newcommand{\smath}{\mathbb{S}}

\newcommand{\nlim}{\lim_{n \to \infty }}
\newcommand{\kinf}{k \to \infty }
\newcommand{\ninf}{n \to \infty }

\newcommand{\ball}{\mathbb{B}}
\newcommand{\ind}[1]{{1\! \!  1 }_{#1}  }

\newcommand{\norm}[3]{{\|  #1 \| }_{#2}^{#3}}
\newcommand{\Norm}[3]{{\Bigl\|  #1 \Bigr\| }_{#2}^{#3}}
\newcommand{\ilsk}[3]{{\bigl( #1 | #2 \bigr) }_{#3}}

\newcommand{\dual}[3]{{\bigl< #1 | #2 \bigr>}_{#3}}
\newcommand{\Dual}[3]{{\Bigl< #1 \bigl| #2 \Bigr>}_{#3}}
\newcommand{\dirilsk}[3]{{\bigl( \! \bigl( #1 | #2 \bigr) \! \bigr) }_{#3}}

\newcommand{\p}{\mathbb{P}}

\newcommand{\e}{\mathbb{E}}

\newcommand{\Xn}{{X}_{n}}

\newcommand{\Pn}{{P}_{n}}
\newcommand{\un}{{u}_{n}}
\newcommand{\bun}[1]{{\bar{u}}_{n}({#1})}
\newcommand{\Jn}[1]{{J}^{n}_{#1}}

\newcommand{\unk}{{u}_{{n}_{k}}}
\newcommand{\taun}{{\tau}_{n}}

\newcommand{\lhs}{{\lcal }_{HS}}

\newcommand{\Bn}{{B}_{n}}

\begin{document}

\title[Stochastic Navier-Stokes Equations]{Stochastic Navier-Stokes Equations driven by L\'{e}vy noise in unbounded 3D domains}


\author{El\.zbieta Motyl}


\maketitle

\begin{abstract}
Martingale solutions of the stochastic Navier-Stokes equations in 2D and 3D possibly unbounded domains, driven by the L\'{e}vy noise consisting of the  compensated  time homogeneous Poisson random measure 
and the Wiener process are considered.  
Using the classical Faedo-Galerkin approximation and the compactness method we prove existence of a martingale solution. We prove also the compactness and tighness criteria in a certain space contained in  some spaces of \it c\`{a}dl\`{a}g \rm functions, \it weakly  c\`{a}dl\`{a}g \rm functions and some Fr\'{e}chet spaces. Moreover, we use a version of the Skorokhod Embedding Theorem for nonmetric spaces.  



\end{abstract}

\footnote{Department of Mathematics and Computer Sciences,
         University of \L \'{o}d\'{z}, Poland,  \\
        emotyl@math.uni.lodz.pl}

\section{Introduction} \label{S:Introduction}

\noindent
Let $\ocal \subset \rd $ be an open  connected possibly unbounded subset with smooth boundary $\partial \ocal $, where $d=2,3$.
We will consider the Navier-Stokes equations
\begin{eqnarray} 
 &  du(t) 
  = & \bigl[  \Delta u  - (u \cdot \nabla ) u + \nabla p +  f(t) \bigr] \, dt +  \int_{Y}F(t,u) \, \tilde{\eta}(dt,dy)  \nonumber \\
& & + G(t,u(t))\, dW(t), \qquad t \in [0,T],   \label{E:NS_intr}
\end{eqnarray}
in $\ocal $, with the incompressibility condition
\begin{equation} \label{E:incompressibility}
  \diver u =0 ,
\end{equation}
the initial condition 
\begin{equation}  \label{E:initial}
   u(0) = {u}_{0},
\end{equation}
and with the homogeneous boundary condition ${u}_{| \partial \ocal }=0$.
In this problem $u=u(t,x)=({u}_{1}(t,x), ...{u}_{d}(t,x))$ and $p=p(t,x)$
represent the velocity and the pressure of the fluid, respectively. Furthermore, $f$ stands for the deterministic
external forces. The terms 

\noindent
$\int_{Y}F(t,u) \, \tilde{\eta}(dt,dy) $, where $\tilde{\eta}$ is a compensated time homogeneous Poisson random measure on a certain measurable space $(Y, \ycal )$,
and $G(t,u(t))\, dW(t)$, where $W$ is a cylindrical Wiener process on some separable Hilbert space ${Y}_{W}$, stand for the random forces.

\bigskip  \noindent
The problem (\ref{E:NS_intr})-(\ref{E:initial}) can be written as the following stochastic evolution equation
\begin{eqnarray*} 
 &du(t) & +\acal u(t) \, dt+B\bigl( u(t) \bigr)\, dt = f(t) \, dt 
 + \int_{Y} F(t,u({t}^{-});y) \tilde{\eta} (dt,dy)  \\
& \qquad \quad &+ G(t,u(t)) \, dW(t)  \qquad t \in [0,T] ,  \\
& u(0) &= {u}_{0} . 
\end{eqnarray*}
We will prove the existence of a martingale solution of the problem (\ref{E:NS_intr})-(\ref{E:initial}) understood as a system 
$(\Omega ,\fcal , \p , \fmath , \eta ,W, u)$, where $(\Omega ,\fcal , \p , \fmath )$ is a filtered probability space, $\eta $ is a time homogeneous Poisson random measure, $W$ is a cylindrical Wiener process 
and $u= ({u}_{t}{)}_{t \in [0,T]}$ is 
a stochastic process with trajectories in the space $\dmath \bigl( [0,T], {H}_{w} \bigr) \cap {L}^{2}(0,T;V )$ and satisfying appropriate  integral equality, see Definition \ref{D:solution} in Section \ref{S:Statement}. 
Here, $V$ and $H $ denote the closures  in ${H}^{1}(\ocal , \rd )$ and ${L}^{2}(\ocal , \rd )$, respectively of the space $\vcal $ of the divergence-free $\rd $ valued vector fields of class ${\ccal }^{\infty }$ with compact supports contained in $\ocal $.
The symbol $\dmath \bigl( [0,T], {H}_{w} \bigr)$
stands for the space of $H$ valued weakly c\`{a}dl\`{a}g  functions.

\bigskip  \noindent
To construct this solution we use the classical Faedo-Galerkin method, i.e.,
\begin{eqnarray*} 
 &  d \un (t) & =   - \bigl[ \Pn \acal \un (t)  + \Bn  \bigl(\un (t) \bigr) - \Pn f (t)  \bigr] \, dt  \\  
&   & \quad  +  \int_{Y} \Pn F(t,\un ({t}^{-}),y) \tilde{\eta} (dt,dy)  
 + \Pn G(t,u(t)) \, dW(t) ,   
\, \,  t \in [0,T]  , \\
&  \un (0) &=  \Pn {u}_{0}.
\end{eqnarray*} 
The solutions $\un $ to the Galerkin scheme generate a sequence of laws $\{ \lcal (\un)$, $n \in \nat \} $ on appropriate functional spaces. To prove that this sequence of probability measures is weakly compact we need appropriate tightness criteria.

\bigskip  \noindent
We concentrate first  on the compactness and tightness criteria. If the domain $\ocal $ is unbounded,
then the embedding $V \subset H$ is not compact. However using Lemma 2.5 in \cite{Holly_Wiciak_1995}, see Appendix C, we can find a separable Hilbert space $U$ such that $U \subset V $, the embedding being dense and compact.

\bigskip  \noindent
We consider the intersection 
$$ 
 {\zcal }_{q}: =   {L}_{w}^{q}(0,T;V)  \cap {L}^{q}(0,T;{H}_{loc }) \cap \dmath ([0,T];U') 
  \cap \dmath ([0,T],{H}_{w}),
$$
where $q\in (1,\infty )$. (The letter $w$ indicates the weak topology.)
By $\dmath ([0,T];U')$ we denote the space of $U'$-valued \it c\`{a}dl\`{a}g \rm  functions 
equipped with the Skorokhod topology and
${L}^{q}(0,T;{H}_{loc})$ stands for the Fr\'{e}chet space defined by (\ref{E:seminorms} ), see Section \ref{S:Det_comp_criterion}.

\bigskip  \noindent
Using the compactness criterion in the space of \it c\`{a}dl\`{a}g \rm functions,
we prove that a set $\kcal $ is relatively compact in ${\zcal }_{q}$
if  the following three conditions hold
\begin{itemize}
\item[(a)] for all $u \in \kcal $ and  all $t \in [0,T]$, $u(t) \in H  $ and 
$\, \, \sup_{u\in \kcal } \sup_{s \in[0,T]} {|u(s)|}_{H} < \infty  $, 
\item[(b)] $ \sup_{u\in \kcal } \int_{0}^{T} \norm{u(s)}{V}{q} \, ds < \infty  $,
  i.e. $\kcal $ is bounded in ${L}^{q}(0,T;V )$,
\item[(c)] $\lim{}_{\delta \to 0 } \sup_{u\in \kcal } {w}_{[0,T],U'}(u;\delta ) =0 $.
\end{itemize} 
Here ${w}_{[0,T],U'}(u;\delta )$ stands for the modulus of the function $u:[0,T] \to U'$.
The above result is a straightforward generalization of the compactness results of \cite{Brzezniak_Motyl_NS} and 
\cite{Motyl_NS_Poisson_pre}. In the paper \cite{Motyl_NS_Poisson_pre} the analogous result is proved in the case when the embedding $V \subset H$ is dense and compact (in the Banach space setting). In  \cite{Brzezniak_Motyl_NS} the embedding $V \subset H$ is only dense and continuous. However, instead of the spaces of c\`{a}dl\`{a}g functions, appropriate spaces of continuous functions are used.
 The present paper generalizes both \cite{Brzezniak_Motyl_NS} and \cite{Motyl_NS_Poisson_pre} in the sense that the embedding $V \subset H$ is dense and continuous and appropriate spaces of c\`{a}dl\`{a}g functions are considered, i.e. $\dmath ([0,T]; U')$ and $\dmath ([0,T],{H}_{w})$.
This approach were strongly inspired by the results due to M\'{e}tivier and Viot, especially the choice of the spaces $\dmath ([0,T]; U')$ and $\dmath ([0,T],{H}_{w})$, see \cite{Metivier_Viot_88} and \cite{Metivier_88}. It is also closely related to the result due to Mikulevicius and Rozovskii \cite{Mikulevicius_Rozovskii} and to the classical Dubinsky compactness criterion, 
\cite{Vishik_Fursikov_88}. However,  both in \cite{Vishik_Fursikov_88} and 
\cite{Mikulevicius_Rozovskii},  the spaces of continuous functions are used. 

\bigskip \noindent
Using the above deterministic compactness criterion and
the Aldous condition in the form given by Joffe and M\'{e}tivier \cite{Joffe_Metivier_86}, see also \cite{Metivier_88}, we obtain the corresponding tightness criterion for the laws on the space ${\zcal }_{q}$, 
see Corollary \ref{C:tigthness_criterion_cadlag_unbound}.

\bigskip  \noindent
We will prove that the set of probalility measures induced by the Galerkin solutions is tight on the space $\zcal$, where
$$ 
  \zcal  : = {L}_{w}^{2}(0,T;V) \cap {L}^{2}(0,T;{H}_{loc}) \cap \dmath ([0,T];U')
  \cap \dmath ([0,T];{H}_{w}) ,
$$
which is not metrizable. 
Further construction a martingale solutions is based on
the Skorokhod Embedding Theorem in nonmetric spaces. In fact, we use the result proved in 
\cite{Motyl_NS_Poisson_pre} and following easily from  the Jakubowski's version of the Skorokhod Theorem \cite{Jakubowski_1998} and the version of the Skorokhod Theorem due to Brze\'{z}niak and Hausenblas 
\cite{Brzezniak_Hausenblas_2010}, see Appendix B. 
This will allow us to construct a stochastic  process $\bar{u}$ with trajectories in the space $\zcal $,
 a time homogeneous Poisson random measure $\bar{\eta }$ and a cylindrical Wiener proces $\bar{W}$ defined on some filtered probability space $(\bar{\Omega },\bar{\fcal }, \bar{\p }, \bar{\fmath })$ such that the system $(\bar{\Omega },\bar{\fcal }, \bar{\p }, \bar{\fmath }, \bar{\eta } , \bar{W}, \bar{u})$ is a martingale solution of the problem (\ref{E:NS_intr})-(\ref{E:initial}). 
In fact, $\bar{u}$ is a process with trajectiories in the space $\zcal $. In particular, the trajectories of $\bar{u}$ are \it weakly c\`{a}dl\`{a}g \rm if  $\bar{u}$ is considered as a $H$-valued process and  \it c\`{a}dl\`{a}g \rm in the bigger space $U'$.

\bigskip  \noindent
The Navier-Stokes equations driven by the compensated Poisson random measure in the 3D bounded domains
 were studied in Dong and Zhai \cite{{Dong_Zhai_2011}}. The authors consider the martingale problem associated to the 
Navier-Stokes equations, i.e. a solution is defined to be a probability measure satisfying appropriate conditions, see Definition 3.1 in \cite{{Dong_Zhai_2011}}. 
The 2D Navier-Stokes equations were considered in \cite{Dong_Xie_2009}, \cite{Dong_Xie_2011} and \cite{Xu_Zhang_2009}. 
In the  present paper, using a different approach we generalize the existence resuls to the case of unbounded 2D and 3D domains. Moreover, we consider  more general noise term.

\bigskip  \noindent
Stochastic Navier-Stokes equations in unbounded 2D and 3D domains were usually considered with the Gaussian noise term, see e.g. \cite{Capinski_Peszat_2001}, \cite{Brzezniak_Peszat_2001}, \cite{Brzezniak_Li_2006} and
\cite{Brzezniak_Motyl_NS}.
Martingale solutions of the stochastic Navier-Stokes equations driven by white noise in  the whole space $\rd $, ($d \ge 2 $), are investigated in \cite{Mikulevicius_Rozovskii}.

\bigskip  \noindent
The present paper is organized as follows. In Section \ref{S:Functional_setting} we recall basic definitions and properties of the spaces and operators appearing in the Navier-Stokes equations.
Section \ref{S:Comp_tight} is devoted to the compactness and tightness results.
Some auxilliary results about the Aldous condition and tightness are contained in Appendix A.
Precise statement of the Navier-Stokes problem driven by L\'{e}vy noise is contained in Section \ref{S:Statement}.  
The main Theorem about existence of a martingale solution of the problem (\ref{E:NS_intr})-(\ref{E:initial}) is proved in Section \ref{S:Existence}. 
Some versions the Skorokhod Embedding Theorems  are recalled in Appendix B. 
In Appendix C we recall  Lemma 2.5 in \cite{Holly_Wiciak_1995} together with the proof.

\section{Functional setting} \label{S:Functional_setting}

\subsection{Basic definitions}

\noindent
Let $\ocal \subset \rd $ be an open connected subset with smooth boundary $\partial \ocal $, $d=2,3$.
Let
\begin{eqnarray}
 &  \vcal := \{ u \in {\ccal }^{\infty }_{c} (\ocal , \rd ) : \, \, \diver u= 0 \} , \nonumber \\
 &  H := \mbox{the closure of $\vcal $ in ${L}^{2}(\ocal , \rd )$} , \label{E:H_space}\\
 &  V := \mbox{the closure of $\vcal $ in ${H}^{1}(\ocal , \rd )$} .  \label{E:V_space}
\end{eqnarray}
In the space $H$ we consider the scalar product and the norm inherited from ${L}^{2}(\ocal , \rd \!)$ and
denote them by $\ilsk{\cdot }{\cdot }{H}$ and $|\cdot {|}_{H}$, respectively, i.e.
$$
\ilsk{u}{v}{H}:= \ilsk{u}{v}{{L}^{2}} , \qquad
 | u {|}_{H} := \norm{u}{{L}^{2}}{} , \qquad u, v \in H.
$$
In the space $V$ we consider the scalar product
inherited from the Sobolev space ${H}^{1}(\ocal , \rd )$, i.e.
\begin{equation} \label{E:V_il_sk}
  \ilsk{u}{v}{V} := \ilsk{u}{v}{{L}^{2}} + \dirilsk{u}{v}{} ,
\end{equation}
where
\begin{equation} \label{E:il_sk_Dir}
  \dirilsk{u}{v}{} :=
\ilsk{\nabla u}{\nabla v}{{L}^{2}} , \qquad u,v \in  V.
\end{equation}
and the norm
\begin{equation} \label{E:norm_V}
  \norm{u}{V}{2} := |u {|}_{H}^{2} + \norm{u}{}{2},
\end{equation}
where $\norm{u}{}{2} := \norm{\nabla u}{{L}^{2}}{2}$.

\noindent
\subsection{The form $b$}
\noindent
Let us consider the following three-linear form, see Temam \cite{Temam79}, 
$$
     b(u,w,v ) = \int_{\ocal  }\bigl( u \cdot \nabla w \bigr) v \, dx .
$$
We will recall those fundamental properties of the form $b$ that are valid both in bounded and unbounded domains.
By the Sobolev embedding Theorem, see \cite{Adams}, and the H\H{o}lder inequality, we obtain the following estimates
\begin{equation}
    |b(u,w,v )|
    \le c \norm{u }{V}{} \norm{w }{V}{} \norm{v }{V}{} , \qquad u,w,v \in V   \label{E:b_estimate_V}
\end{equation}
for some positive constant $c$. Thus the form $b$ is continuous on $V$, see also \cite{Temam79}.
Moreover, if we define a bilinear map $B$ by $B(u,w):=b(u,w, \cdot )$, then by inequality (\ref{E:b_estimate_V}) we infer that $B(u,w) \in {V}_{}'$ for all $u,w\in V$ and that the following inequality holds
\begin{equation}  \label{E:estimate_B}
 |B(u,w) {|}_{V'} 
   \le c  \norm{u }{V}{}\norm{w }{V}{},\qquad u,w \in V .
\end{equation}
Moreover, the mapping $B: V \times V \to V' $ is bilinear and continuous.
Let us also recall the following properties of the form $b$, see Temam \cite{Temam79}, Lemma II.1.3,
$$  \label{E:antisymmetry_b}
b(u,w, v ) =  - b(u,v ,w), \ \ \  u,w,v \in V .
$$
In particular,
$$  \label{E:wirowosc_b}
b(u,v,v) =0   \qquad u,v \in V.
$$
Hence
$$  \label{E:antisymmetry_B}
\dual{B(u,w)}{v}{}  =  - \dual{B(u,v)}{w}{}, \qquad u,w,v \in V 
$$
and
\begin{equation}  \label{E:wirowosc_B}
\dual{B(u,v)}{v}{}  = 0, \qquad u,v \in V .
\end{equation}
Let us, for any $m>0$ define the following standard scale of Hilbert spaces
$$
  {V}_{m} := \mbox{the closure of $\vcal $ in ${H}^{m}(\ocal , \rd )$} .
$$
If $m > \frac{d}{2} +1$ then by the Sobolev embedding Theorem, see \cite{Adams},
$$
    {H}^{m-1}(\ocal  , \rd ) \hookrightarrow  {\ccal }_{b}(\ocal  , \rd )
   \hookrightarrow {L}^{\infty } (\ocal , \rd ),
$$
where ${\ccal }_{b}(\ocal ,\rd )$ denotes the space of $\rd $-valued continuous and bounded functions defined on $\ocal $.
If $u,w \in V$ and $v \in {V}_{m}$ with $m > \frac{d}{2} +1$ then
\begin{eqnarray*}
 |b(u,w,v)| & = & |b(u,v,w)|
 = \Bigl| \sum_{i=1}^{n} \int_{\ocal } {u}_{i} w \frac{\partial v}{\partial {x}_{i}} \, dx \Bigr|  \\
& \le & \norm{u}{{L}^{2}}{} \norm{w}{{L}^{2}}{} \norm{\nabla v}{{L}^{\infty }}{}
 \le  {c}_{} \norm{u}{{L}^{2}}{} \norm{w}{{L}^{2}}{} \norm{v}{{V}_{m}}{}
\end{eqnarray*}
for some constant ${c}_{} >0 $. Thus, $b$ can be uniquely extented to the three-linear form (denoted by the same letter)
$$
   b : H \times H \times {V}_{m} \to \rzecz
$$
and $|b(u,w,v)| \le  {c}_{} \norm{u}{{L}^{2}}{} \norm{w}{{L}^{2}}{} \norm{v}{{V}_{m}}{}$
for $u,w \in H$ and $v \in {V}_{m}$. At the same time the operator $B$ can be uniquely extended
to a bounded bilinear operator
$$
    B : H \times H \to {V}_{m}' .
$$
In particular, it satisfies the following estimate
\begin{equation}  \label{E:estimate_B_ext}
 |B(u,w) {|}_{{V}_{m}'} \le c {|u|}_{H}  {|w|}_{H} ,\qquad u,w \in H.
\end{equation}
See Vishik and Fursikov\cite{Vishik_Fursikov_88}.
We will also use the following notation, $B(u):=B(u,u)$.
Let us also recall the  well known result that
the map $B:V \to V'$ is locally Lipschitz continuous, i.e. for every $r>0$ there exists a constant ${L}_{r}$ such that
\begin{equation} \label{E:estimate_B_Lipsch}
   \bigl| B(u) - B(\tilde{u}) {\bigr| }_{V'} \le {L}_{r} \norm{u - \tilde{u}}{V}{} ,
   \qquad u , \tilde{u } \in V , \quad \norm{u}{V}{}, \norm{\tilde{u}}{V}{} \le r  .
\end{equation}

\subsection{The space $U$ and some operators} \label{S:Some_operators}
We recall operators and their properties used in \cite{Brzezniak_Motyl_NS}.
Here we also recall the definition of a Hilbert space $U$  compactly embedded in appropriate space ${V}_{m}$. This is possible thanks to the result due to Holly and Wiciak, \cite{Holly_Wiciak_1995} which
we recall  with the proof in Appendix C, see Lemma \ref{L:2_5_Holly_Wiciak}.
This space will be of crucial importance in further investigations. 

\bigskip  \noindent
Consider the natural embedding $j: V \hookrightarrow H $ and its adjoint ${j}^{*} : H \to V $.
Since the range of $j$ is dense in $H$, the map ${j}^{*}$ is one-to-one.
Let us put
\begin{eqnarray} 
 & D(A) &:= {j}^{*}(H) \subset V , \nonumber \\
 & Au &:= \bigl( {j}^{*} {\bigr) }^{-1} u , \qquad u \in D(A) .   \label{E:op_A}
\end{eqnarray}
and
\begin{equation} \label{E:op_Acal}
     \acal u :=\dirilsk{u}{\cdot }{} , \qquad u \in V,
\end{equation}
where $\dirilsk{\cdot }{\cdot }{}$ is defined by (\ref{E:il_sk_Dir}). Let us notice that
if $u \in V$, then $\acal u \in V'$ and
$$ |\acal u{|}_{V'}\le \norm{u}{}{}.  $$
Indeed, this follows immediately from (\ref{E:norm_V}) and the following inequalities
$$
  |\dirilsk{u}{v}{}| \le \norm{u}{}{} \cdot \norm{v}{}{}
   \le \norm{u}{}{} (\norm{v}{}{2} + |v{|}_{H}^{2} {)}^{\frac{1}{2}}
   = \norm{u}{}{} \cdot \norm{v}{V}{}, \quad v \in V.
$$

\bigskip
\begin{lemma} \label{L:A_acal_rel}\ \rm (Lemma 2.2 in \cite{Brzezniak_Motyl_NS}) \it 
\begin{itemize}
\item[(a)] For any $u \in D(A)$ and $v \in V$:
$$
((A-I)u|v{)}_{H} =  \dirilsk{u}{v}{} = \dual{\acal  u }{v}{} ,
$$
where $I$ stands for the identity operator on $H$ and $\dual{}{}{}$ is the standard duality pairing.
In particular,
$$
         |\acal u{|}_{V'}\le |(A-I)u{|}_{H}.
$$
\item[(b)] $D(A)$ is dense in $H$.
\end{itemize}
\end{lemma}

\proof
To prove assertion (a), let $u\in D(A)$ and $v\in V$. Then
\begin{eqnarray*}
&(Au|v{)}_{H} &= (( {j}^{*} {)}^{-1} u |v{)}_{H} = (( {j}^{*} {)}^{-1} u |jv{)}_{H}
        = \ilsk{{j}^{*}( {j}^{*} {)}^{-1}u}{v}{V} = \ilsk{u}{v}{V} \\
&       &= \ilsk{u}{v}{H} + \dirilsk{u}{v}{}
        = \ilsk{Iu}{v}{H} + \dual{ \acal  u}{v}{}.
\end{eqnarray*}
Let us move to the proof of part (b). Since $V$ is dense in $H$, it is sufficient to prove that $D(A)$ is dense in $V$. Let $w\in V$ be an arbitrary element orthogonal to $D(A)$ with respect to the scalar product in $V$. Then
$$
     \ilsk{u}{w}{V} = 0 \qquad \mbox{for} \quad u \in D(A).
$$
On the other hand, by (a) and (\ref{E:V_il_sk}),
$\ilsk{u}{w}{V} = \ilsk{Au}{w}{H}$ for $u \in D(A)$.
Hence $\ilsk{Au}{w}{H}=0$ for $u\in D(A)$. Since $A:D(A) \to H$ is onto, we infer that $w=0$,
which completes the proof. \qed

\bigskip  \noindent
Let us assume that $m > 1$. It is clear that ${V}_{m}$ is dense in $V$ and the embedding ${j}_{m}: {V}_{m} \hookrightarrow V$ is continuous. Then by Lemma \ref{L:2_5_Holly_Wiciak} in Appendix C,
there exists a Hilbert space $U$ such that $U \subset  {V}_{m}$, $U$ is dense in $ {V}_{m}$ and
\begin{equation} \label{E:U_comp_V_m}
 \mbox{\it the natural embedding ${\iota }_{m}: U \hookrightarrow  {V}_{m}$ is compact \rm }.
\end{equation}
Then we have
\begin{equation} \label{E:embeddings}
  U \stackrel{{\iota }_{m}}{\hookrightarrow } {V}_{m}
  \stackrel{{j }_{m}}{\hookrightarrow }V
 \stackrel{{j}_{}}{\hookrightarrow} H \cong H'
\stackrel{j'}{\hookrightarrow } V'
\stackrel{{j}_{m}'}{\hookrightarrow } {V}_{m}'
  \stackrel{{\iota }_{m}^{'}}{\hookrightarrow }U' .
\end{equation}
Since the embedding  ${\iota }_{m}$ is compact, ${\iota }_{m}'$ is compact as well.
Consider the composition
$$
   \iota := j \circ {j}_{m} \circ {\iota }_{m} : U \hookrightarrow H
$$
and its adjoint
$$
   {\iota }^{*} := (j \circ {j}_{m} \circ {\iota }_{m} )^{*}
  ={\iota }_{m} ^{*} \circ {j}_{m}^{*} \circ {j}^{*} :  H \to U .
$$
Note that $\iota $ is compact and since the range of $\iota $ is dense in $H$, ${\iota }^{*} : H \to U $ is one-to-one. Let us put
\begin{eqnarray} 
 & D(L) &:= {\iota }^{*}(H) \subset U , \nonumber \\
 & Lu &:= \bigl( {\iota }^{*} {\bigr) }^{-1} u , \qquad u \in D(L) .  \label{E:op_L}
\end{eqnarray}
It is clear that $L:D(L) \to H $ is onto. Let us also notice that
\begin{equation} \label{E:op_L_ilsk}
    \ilsk{Lu}{w}{H} = \ilsk{u}{w}{U}, \qquad u \in D(L), \quad w \in U .
\end{equation}
By equality (\ref{E:op_L_ilsk}) and the densiness of $U$ in $H$, we infer similarly as in the proof of assertion (b) in Lemma \ref{L:A_acal_rel} that $D(L)$ is dense in $H$.
Moreover, for $u \in D(L)$,
\begin{eqnarray*}
&Lu & = \bigl( {\iota }^{*} {\bigr) }^{-1} u
   = \bigl( {\iota }_{m} ^{*} \circ {j}_{m}^{*} \circ {j}^{*}  {\bigr) }^{-1} u
  = A \circ \bigl( {j}_{m}^{*} {\bigr) }^{-1} \circ   \bigl( {\iota }_{m} ^{*} {\bigr) }^{-1} u,
\end{eqnarray*}
where $A$ is defined by (\ref{E:op_A}).

\bigskip  \noindent
Since $L$ is self-adjoint and ${L}^{-1}$ is compact, there exists an orthonormal basis $\{ {e}_{i} {\} }_{i \in \nat }$ of $H$ composed of the eigenvectors of operator $L$. Let us fix $n \in \nat $ and let $\Pn $ be the operator from $U'$ to $span \{ {e}_{1},..., {e}_{n}\} $ defined by
\begin{equation} \label{E:P_n}
  \Pn {u}^{*} := \sum_{i=1}^{n} \bigl< {u}^{*}| {e}_{i}\bigr> {e}_{i}, \qquad {u}^{*} \in U',
\end{equation}
where $\dual{\cdot }{\cdot }{}$ denotes the duality pairing between the space $U$ and its dual $U'$.
Note that the restriction of $\Pn $ to $H$, denoted still by $\Pn $, is given by 
$$
   \Pn u = \sum_{i=1}^{n} \ilsk{ u}{ {e}_{i}}{H}  {e}_{i}, \qquad  u \in H ,
$$
and thus it is the $\ilsk{\cdot }{\cdot }{H}$-orthogonal projection
onto  $span \{ {e}_{1},..., {e}_{n}\} $. Restrictions of $\Pn $ to other spaces considered in (\ref{E:embeddings}) will also be denoted by $\Pn $. Moreover, it is easy to see that
$$
   \ilsk{\Pn {u}^{*}}{v}{H} = \dual{ {u}^{*}}{\Pn v}{} , \qquad {u}^{*} \in U', \quad v \in U.
$$
It is easy to prove that the system $\bigl\{ \frac{{e}_{i}}{\norm{{e}_{i}}{U}{}} {\bigr\} }_{n \in \nat }$ is the $\ilsk{\cdot }{\cdot }{U}$-orthonormal basis in the space $U$ and that  the restriction of ${P}_{n}$ to $U$ is the $\ilsk{\cdot }{\cdot }{U}$-projection onto the subspace
$span \{ {e}_{1},...,{e}_{n}  \} $.
In particular, for every $u\in U$
\begin{description}
\item[(i)] $\lim_{n \to \infty } \norm{{P}_{n}u-u}{U}{} =0$,
\item[(ii)] $\lim_{n \to \infty } \norm{{P}_{n}u-u}{{V}_{m}}{} =0$, where $m>0$,
\item[(iii)] $\lim_{n \to \infty } \norm{{P}_{n}u-u}{V}{} =0$.
\end{description}
See Lemma 2.4 in \cite{Brzezniak_Motyl_NS} for details.

\bigskip  \noindent
We will use the basis $\{ {e}_{i} {\} }_{i \in \nat }$ and the operators $\Pn $ in the Faedo-Galerkin approximation.

\section{Compactness results }  \label{S:Comp_tight}

\subsection{The space of c\`{a}dl\`{a}g functions}

Let $(\smath ,\varrho )$ be a separable and complete metric space.
Let $\dmath ([0,T];\smath )$  the space of all $\smath $-valued \it c\`{a}dl\`{a}g \rm functions defined on $[0,T]$, i.e. the functions  which are right continuous and have left limits at every $t\in [0,T]$ . 
The space $\dmath ([0,T];\smath )$ is endowed with the Skorokhod topology.

\bigskip  \noindent
\begin{remark} \label{R:cadlag_conv}
A sequence $({u}_{n}) \subset \dmath ([0,T];\smath )$ converges to $u \in \dmath ([0,T];\smath )$ iff there exists a sequence $({\lambda }_{n})$ of homeomorphisms of  $[0,T]$ such that ${\lambda }_{n} $ tends to the identity uniformly on $[0,T]$ and ${u}_{n} \circ {\lambda }_{n} $ tends to $u$ uniformly on $[0,T]$.
\end{remark} 

\bigskip \noindent
This topology is  metrizable by the following metric ${\delta }_{T}$
$$ \label{E:II_1.1.1_[Metivier_88]}
 {\delta }_{T}(u,v) := \inf_{\lambda \in {\Lambda }_{T}} \Bigl[ \sup_{t\in [0,T]}
  \varrho  \bigl( u(t), v \circ \lambda (t)\bigr)  
  + \sup_{t\in [0,T]} |t-\lambda (t)|
 + \sup_{s\ne t} \Bigl| \log \frac{\lambda (t)-\lambda (s)}{t-s} \Bigr| \Bigr]  ,
$$
where ${\Lambda }_{T}$ is the set of increasing homeomorphisms of $[0,T]$.
Moreover, 

\noindent
$\bigl(  \dmath ([0,T];\smath ),{\delta }_{T}\bigr) $ is a complete metric space,
see \cite{Joffe_Metivier_86}.

\bigskip  \noindent
Let us recall the notion of a \it modulus \rm of the function. It plays analogous role in the space $\dmath ([0,T];\smath )$ as \it the modulus of continuity \rm in the space of continuous functions $\cmath ([0,T];\smath )$.

\bigskip  \noindent
\begin{definition}  \rm (see \cite{Metivier_88})
Let $u \in \dmath ([0,T];\smath )$ and let $ \delta >0 $ be given. A \bf modulus \rm of $u$ is defined by 
\begin{equation}  \label{E:modulus_cadlag}
  {w}_{[0,T],\smath }(u,\delta ) : = \inf_{{\Pi }_{\delta }} \, \max_{{t}_{i} \in \bar{\omega }} \,
  \sup_{{t}_{i} \le s <t < {t}_{i+1} \le T } \varrho  \bigl( u(t), u(s) \bigr) ,
\end{equation}
where ${\Pi }_{\delta }$ is the set of all increasing sequences 
$
   \bar{\omega } = \{ 0= {t}_{0} < {t}_{1} < ... < {t}_{n} =T  \}
$
with the following property
$$
      {t}_{i+1} - {t}_{i}   \ge \delta , \qquad i=0,1,...,n-1.
$$
If no confusion seems likely, we will denote the modulus by ${w}_{[0,T]}(u,\delta )$.

\end{definition}

\bigskip  \noindent
We have the following criterion for relative compactness of a subset of the space $\dmath ([0,T];\smath )$,
see \cite{Joffe_Metivier_86},\cite{Metivier_88}, Ch.II, and \cite{Billingsley}, Ch.3,
analogous to the Arzel\`{a}-Ascoli Theorem for the space of continuous functions.

\bigskip
\begin{theorem} \label{T:cadlag_compactness} \it 
A set $A \subset \dmath ([0,T];\smath ) $ has compact closure iff it satisfies the following two conditions:
\begin{itemize}
\item[(a) ] there exists a dense subset $J \subset [0,T]$ such that 
for every $ t \in J  $  the set $\{ u(t), \, \, u \in A  \} $ has compact closure in $\smath $.
\item[(b) ] $ \lim_{\delta \to 0} \, \sup_{u \in A}  \, {w}_{[0,T]}(u,\delta ) =0 $.
\end{itemize}
\end{theorem}

\subsection{Deterministic compactness criterion} \label{S:Det_comp_criterion}

\noindent
Let us recall that $V$ and $H$ are Hilbert spaces defined by (\ref{E:H_space})-(\ref{E:norm_V}).
Since $\ocal $ is an arbitrary domain of $\rd $, ($d=2,3$), the embedding $V\hookrightarrow H$ is dense and continuous. We have defined a Hilbert space $U\subset V$ such that the embedding $U\hookrightarrow V$ is dense and compact, see (\ref{E:U_comp_V_m}). In particular, we have
$$
   U \hookrightarrow V \hookrightarrow H\cong H'  \hookrightarrow U',
$$
the embedding $H\hookrightarrow U'$ being compact as well. 
Let $\bigl( {\ocal }_{R} {\bigr) }_{R \in \nat } $ be a sequence of open and bounded subsets of $\ocal $ with
regular boundaries $\partial {\ocal }_{R}$  such that
${\ocal }_{R} \subset {\ocal }_{R+1}$ and $\bigcup_{R=1}^{\infty } {\ocal }_{R} = \ocal $.
We will consider the following spaces of restrictions of functions defined on $\ocal $ to subsets ${\ocal }_{R}$, i.e.
\begin{equation}  \label{E:HR_VR}
 {H}_{{\ocal }_{R}} := \{ {u}_{|{\ocal }_{R}} ; \, \, \,  u \in H  \}
 \qquad {V}_{{\ocal }_{R}} := \{ {v}_{|{\ocal }_{R}} ; \, \, \,  v \in V  \}
\end{equation}
with appropriate scalar products and norms, i.e.
\begin{eqnarray*}
 & \ilsk{u}{v}{{H}_{{\ocal }_{R}}} &:= \int_{{\ocal }_{R}} u v \, dx , \qquad  u,v \in {H}_{{\ocal }_{R}}, \\
& \ilsk{u}{v}{{V}_{{\ocal }_{R}}} &:=\int_{{\ocal }_{R}}uv\, dx+\int_{{\ocal }_{R}}\nabla u\nabla v \, dx,
 \qquad  u,v \in {V}_{{\ocal }_{R}}
\end{eqnarray*}
and $|u{|}_{{H}_{{\ocal }_{R}}}^{2}:=\ilsk{u}{u}{{H}_{{\ocal }_{R}}}$ for
$u \in {H}_{{\ocal }_{R}}$ and $\norm{u}{{V}_{{\ocal }_{R}}}{2}:=\ilsk{u}{u}{{V}_{{\ocal }_{R}}}$ for
$u \in {V}_{{\ocal }_{R}}$.
The symbols  ${H}_{{\ocal }_{R}}'$ and ${V}_{{\ocal }_{R}}'$ will stand for the corresponding dual spaces.

\bigskip \noindent
Since the sets ${\ocal }_{R}$ are bounded,
\begin{equation} \label{E:comp_VR_HR}
\mbox{\it the embeddings ${V}_{{\ocal }_{R}} \hookrightarrow  {H}_{{\ocal }_{R}}$ are compact. \rm }
\end{equation}

\bigskip  \noindent
Let $q\in (1,\infty )$.
Let us consider the following three functional spaces, analogous to those considered in 
\cite{Motyl_NS_Poisson_pre} and
\cite{Brzezniak_Motyl_NS}, see also \cite{Metivier_88} , \cite{Metivier_Viot_88}:
\begin{eqnarray}
\dmath ([0,T],U')&:=& \mbox{the space of  c\`{a}dl\`{a}g functions } u:[0,T] \to U'
                          \mbox{ with the }  \nonumber \\
          &               & \mbox{ topology }  {\tcal }_{1}\mbox{ induced by the Skorokhod metric ${\delta }_{T}$}, \nonumber \\
{L}_{w}^{q}(0,T;V)&:=& \mbox{the space } {L}^{q} (0,T;V) \mbox{ with the weak topology }
                      {\tcal}_{2},  \nonumber  \\
{L}^{q}(0,T;{H}_{loc})&:=& \mbox{the space of measurable functions } u:[0,T] \to H 
                          \nonumber \\
                    &    &  \mbox{such that for all } R \in \nat
                       \nonumber \\
      {p}_{T,R}^{}&(u):=&\norm{u}{{L}^{q}(0,T;{H}_{{\ocal }_{R}})}{} :=
     \Bigl( \! \int_{0}^{T}\! \! \int_{{\ocal }_{R}} \! |u(t,x){|}^{q} dxdt \! {\Bigr) }^{\frac{1}{q}} <\infty  ,
      \label{E:seminorms} \\
   &  & \! \! \! \mbox{with the topology }  {\tcal }_{3} \mbox{ generated by the seminorms } \nonumber \\
    &  &     ({p}_{T,R}{)}_{R\in \nat } . \nonumber
\end{eqnarray}

\noindent
Let ${H}_{w}$ denote the Hilbert space $H$ endowed with the weak topology. 
Let us consider the fourth space, see \cite{Motyl_NS_Poisson_pre},
\begin{eqnarray} 
\dmath ([0,T];{H}_{w})& : = & \mbox{the space of weakly c\`{a}dl\`{a}g functions } 
                          u : [0,T] \to H \mbox{ with the} \nonumber \\
                   &    & \mbox{weakest topology ${\tcal }_{4}$ such that for all 
                          $h \in h $   the  mappings } \nonumber 
\end{eqnarray}
\begin{eqnarray} 
      \dmath ([0,T];&{H}_{w}) \ni u  \mapsto \ilsk{u(\cdot )}{h}{H} 
      \in \dmath  ([0,T];\rzecz ) 
     \mbox{ are continuous. }       \label{E:D([0,T];H_w)_cadlag}       
\end{eqnarray}
In particular,  
$\un \to u $ in $\dmath ([0,T];{H}_{w}) $ iff  for all $ h \in H $:
$$
  \ilsk{\un (\cdot )}{h}{H}  \to \ilsk{u(\cdot )}{h}{H}  \quad \mbox{in the space } \quad 
  \dmath ([0,T];\rzecz ).
$$ 
Let us consider the ball
$$
    \ball := \{ x \in H : \, \, \, {|x|}_{H} \le r \} .
$$
Let ${\ball }_{w}$ denote the ball $\ball $ endowed with the weak topology.
It is well-known that the ${\ball }_{w}$ is metrizable, see \cite{Brezis}. Let ${q}_{r}$ denote the metric compatible with the weak topology on $\ball $. Let us consider the following space
\begin{eqnarray} 
\dmath ([0,T]; {\ball }_{w})  = & \mbox{ the space of weakly c\`{a}dl\`{a}g functions } 
                          u : [0,T] \to H \nonumber \\
                         & \mbox{ and such that }  
                        \sup_{t \in [0,T]} {|u(t)|}_{H} \le r  .  \label{E:D([0,T];B_w)}           
\end{eqnarray}
Then $\dmath ([0,T]; {\ball }_{w})$ is metrizable with 
\begin{equation} \label{E:metric_D([0,T];B_w)}
     {\delta }_{T,r}(u,v) = \inf_{\lambda \in {\Lambda }_{T}} \! \biggl\{ \sup_{t\in [0,T]} \! {q}_{r}(u(t),v\circ \lambda (t))\! + \! \sup_{t \in [0,T]} \! |t-\lambda (t)|  + \sup_{s \ne t } \! \Bigl| \log \frac{\lambda (t)-\lambda (s)}{t-s} \Bigr|  \biggr\} .
\end{equation}
Since by the Banach-Alaoglu Theorem  ${\ball }_{w}$ is compact, $(\dmath ([0,T];{\ball }_{w}),{\delta }_{T,r} )$ is a complete metric space. 

\bigskip  \noindent
The following lemma says that any  sequence $(\un ) \subset {L}^{\infty } (0,T;H)$  convergent in 

\noindent
$\dmath ([0,T];U')$ is also convergent in the space $\dmath ([0,T];{\ball }_{w})$. 

\bigskip
\begin{lemma} (see Lemma 4.3 in \cite{Motyl_NS_Poisson_pre}) \label{L:D(0,T,{H}_{w})_conv}
Let ${u}_{n}:[0,T] \to H $, $n \in \nat $, be functions such that
\begin{itemize}
\item[(i)] $\sup_{n \in \nat } \sup_{s \in [0,T]} {|\un (s)|}_{H} \le r  $,
\item[(ii)] $\un \to u $ in $\dmath ([0,T];U')$.
\end{itemize}
Then  $u, \un \in \dmath ([0,T];{\ball }_{w}) $ and $\un \to u$ in $\dmath ([0,T];{\ball }_{w})$ as $n \to \infty $.
\end{lemma}

\bigskip \noindent
We recall the proof in Appendix E.

\bigskip  \noindent
The following Theorem is a generalization of the results of \cite{Brzezniak_Motyl_NS} and 
\cite{Motyl_NS_Poisson_pre}. In the paper \cite{Motyl_NS_Poisson_pre} the analogous result is proved in the case when the embedding $V \subset H$ is dense and compact. In  \cite{Brzezniak_Motyl_NS} the embedding $V \subset H$ is only dense and continuous. However, instead of the spaces of c\`{a}dl\`{a}g functions, appropriate spaces of continuous functions are used. The following result generalizes both \cite{Brzezniak_Motyl_NS} and \cite{Motyl_NS_Poisson_pre} in the sense that the embedding $V \subset H$ is dense and continuous and appropriate spaces of c\`{a}dl\`{a}g functions are considered, i.e.
$\dmath ([0,T]; U')$ and $\dmath ([0,T],{H}_{w})$.

\bigskip  \noindent
\begin{theorem} \rm 
\label{T:Dubinsky_cadlag_unbound} \it
Let $q\in (1,\infty )$ and let 
\begin{equation} \label{E:Z_cadlag}
 {\zcal }_{q}: =   {L}_{w}^{q}(0,T;V)  \cap {L}^{q}(0,T;{H}_{loc}) \cap \dmath ([0,T]; U') 
  \cap \dmath ([0,T],{H}_{w})
\end{equation}
and let $\tcal $ be  the supremum of the corresponding topologies. Then a set $\kcal \subset {\zcal }_{q}$ is $\tcal $-relatively compact if  the following three conditions hold
\begin{itemize}
\item[(a)  ] for all $u \in \kcal $ and  all $t \in [0,T]$, $u(t) \in H  $ and 
$\, \, \sup_{u\in \kcal } \sup_{s \in[0,T]} {|u(s)|}_{H} < \infty  $, 
\item[(b)] $ \sup_{u\in \kcal } \int_{0}^{T} \norm{u(s)}{V}{q} \, ds < \infty  $,
  i.e. $\kcal $ is bounded in ${L}^{q}(0,T;V)$,
\item[(c)] $\lim{}_{\delta \to 0 } \sup_{u\in \kcal } {w}_{[0,T],U'}(u;\delta ) =0 $.
\end{itemize}
\end{theorem}

\proof
We can  assume that $\kcal $ is a closed subset of ${\zcal }_{q}$. Because of the assumption (b), the weak topology in  ${L}_{w}^{q}(0,T;V)$ induced on ${\zcal }_{q}$ is metrizable. 
Since the topology in 
${L}^{q}(0,T;{H}_{loc})$ is defined by the countable family of seminorms (\ref{E:seminorms}), this space is also metrizable.
By assumption (a), it is sufficient to consider the metric  subspace $ \dmath ([0,T]; {\ball }_{w}) \subset \dmath ([0,T],{H}_{w}) $ defined by (\ref{E:D([0,T];B_w)}) and 
(\ref{E:metric_D([0,T];B_w)}) 
with $r:=\sup_{u\in \kcal } $ $\sup_{s \in[0,T]} $ $ {|u(s)|}_{H}$.
Thus compactness of a subset of ${\zcal }_{q}$ is equivalent to its sequential compactness. Let $(\un )$ be a sequence in $\kcal $. By the Banach-Alaoglu Theorem condition (b) yields that the set $\kcal $ is compact in $ {L}_{w}^{q}(0,T;V) $.

\bigskip  \noindent
Using the compactness criterion in the space of c\`{a}dl\`{a}g functions contained in Theorem 
\ref{T:cadlag_compactness}, we will prove that $(\un )$ is compact in  $ \dmath ([0,T]; U')$.
Indeed, by (a) for every $t\in [0,T]$ the set $\{ \un (t), n\in \nat  \} $ is bounded in $H$. Since the embedding $H \subset U'$ is compact, the set $\{ \un (t), n\in \nat  \} $ is compact in $U'$.
This together with condition (c) implies compactness of the sequence $(\un )$ in the space 
$\dmath ([0,T];U')$.

\bigskip  \noindent
Therefore there exists a subsequence $(\unk ) \subset (\un )$ such that 
$$
  \unk \to u \quad \mbox{in} \quad {L}_{w}^{q}(0,T;V) \cap \dmath ([0,T]; U') 
\quad \mbox{as } \quad \kinf .
$$ 
Since $\unk \to u$ in $\dmath ([0,T]; U')$, $\unk (t) \to u(t)$ in $U'$ for all continuity points of function $u$, (see \cite{Billingsley}). 
By condition (a) and the Lebesgue dominated convergence theorem, we infer that
for all $p \in [1,\infty )$
$$
  \unk \to u \quad \mbox{in} \quad {L}^{p}(0,T;U')  
\quad \mbox{as } \quad \kinf .
$$
We claim that
$$ \label{E:unk_u_L^q_H_loc}
  \unk \to u \quad \mbox{in} \quad  {L}^{q}(0,T; {H}_{loc})
\quad \mbox{as } \quad \kinf .
$$
In order to prove it let us fix $R >0$.
Since, by (\ref{E:comp_VR_HR}) the embedding ${V}_{{\ocal }_{R}} \hookrightarrow {H}_{{\ocal }_{R}}$ is compact and the embeddings 
${H}_{{\ocal }_{R}}  \hookrightarrow H' \hookrightarrow U'$ are continuous,  by the Lions Lemma, \cite{Lions_69}, for every $\eps >0 $ there exists a costant $C = {C}_{\eps , R}>0 $ such that
$$
  {|u|}_{{H}_{{\ocal }_{R}}}^{q} \le \eps \norm{u}{{V}_{{\ocal }_{R}}}{q}  + {C}_{\eps } {|u|}_{U'}^{q} ,
\qquad u \in V .
$$
Thus for almost all $s \in [0,T]$
$$
  {|\unk (s) - u(s )|}_{{H}_{{\ocal }_{R}}}^{q} \le \eps \norm{\unk (s) - u(s)}{{V}_{{\ocal }_{R}}}{q}
 + {C}_{\eps } {|\unk (s) - u(s)|}_{U'}^{q} ,
\quad k \in \nat ,
$$
and so for all $k \in \nat $
$$
  {\| \unk  - u \| }_{{L}^{q}(0,T;{H}_{{\ocal }_{R}})}^{q}
  \le \eps \norm{\unk  - u }{{L}^{q}(0,T;{V}_{{\ocal }_{R}})}{q}
 + {C}_{\eps } {\| \unk  - u \| }_{{L}^{q}(0,T;U')}^{q}  .
$$
Passing to the upper limit as $\kinf $ in the above inequality and using the estimate
$$
 \norm{\unk  - u }{{L}^{q}(0,T;{V}_{{\ocal }_{R}})}{q}
   \le q \bigl( \norm{\unk   }{{L}^{q}(0,T;{V}_{{\ocal }_{R}})}{q}
    + \norm{ u }{{L}^{q}(0,T;{V}_{{\ocal }_{R}})}{q} \bigr)
 \le 2q {c}_{q} ,
$$
where ${c}_{q}= \sup_{u \in \kcal } \norm{u}{{L}^{q}(0,T;V)}{q}$, we infer that
$$
  \limsup_{\kinf }  {\| \unk  - u \| }_{{L}^{q}(0,T;{H}_{{\ocal }_{R}})}^{q} \le 2q {c}_{q} \eps  ,
$$
By the arbitrariness of $\eps $,
$$
  \lim_{\kinf }  {\| \unk  - u \| }_{{L}^{q}(0,T;{H}_{{\ocal }_{R}})}^{q} =0 .
$$
The proof of Theorem is thus complete. \qed

\subsection{Tightness criterion}

\noindent
Let us recall that $U,V,H$ are separable Hilbert spaces such that
$$
    U \hookrightarrow V \hookrightarrow H , 
$$ 
where the embedding $U\hookrightarrow V$ is compact and $V\hookrightarrow H$ is continuous.
Using the compactness criterion  formulated in  Theorem  \ref{T:Dubinsky_cadlag_unbound}
 we obtain the corresponding tightness criterion in the space ${\zcal }_{q}$. Let us first recall that the space ${\zcal }_{q}$ is defined by
$$
 {\zcal }_{q}: =   {L}_{w}^{q}(0,T;V)  \cap {L}^{q}(0,T;{H}_{loc}) \cap \dmath ([0,T]; U') 
  \cap \dmath ([0,T],{H}_{w})
$$
and it is equipped with the topology $\tcal $, see (\ref{E:Z_cadlag}).

\begin{corollary} \bf (tightness criterion) \it \label{C:tigthness_criterion_cadlag_unbound}
Let $(\Xn {)}_{n \in \nat }$ be a sequence of c\`{a}dl\`{a}g $\mathbb{F}$-adapted $U'$-valued processes such that 
\begin{itemize}
\item[(a)] there exists a positive constant ${C}_{1}$ such that 
$$
         \sup_{n\in \nat}\e \bigl[ \sup_{s \in [0,T]} {|\Xn (s) |}_{H}  \bigr]  \le {C}_{1} ,
$$ 
\item[(b)] there exists a positive constant ${C}_{2}$ such that 
$$
    \sup_{n\in \nat}\e \Bigl[  \int_{0}^{T} \norm{\Xn (s)}{V}{q} \, ds    \Bigr]  \le {C}_{2} ,
$$ 
\item[(c)]  $(\Xn {)}_{n \in \nat }$ satisfies the Aldous condition \bf [A] \rm in $U'$. 
\end{itemize}
Let ${\tilde{\p }}_{n}$ be the law of $\Xn $ on ${\zcal }_{q}$.
Then for every $\eps >0 $ there exists a compact subset ${K}_{\eps }$ of ${\zcal }_{q}$ such that 
$$
   {\tilde{\p }}_{n} ({K}_{\eps })  \ge 1 - \eps .
$$ 
\end{corollary}

\noindent
We recall the Aldous condition \bf [A] \rm in Appendix A, see Definition \ref{D:Aldous}. The proof of Corrollary \ref{C:tigthness_criterion_cadlag_unbound} is postponed to Appendix A, as well.  

\section{Stochastic Navier-Stokes equations driven by  L\'{e}vy noise}

\subsection{Time homogeneous Poisson random measure}

\noindent
We follow the approach due to Brze\'{z}niak and Hausenblas \cite{Brzezniak_Hausenblas_2009}, 
\cite{Brzezniak_Hausenblas_2010}, see also \cite{Ikeda_Watanabe_81} and \cite{Peszat_Zabczyk_2007}. 
Let us denote
 $\nat :=\{ 0,1,2,... \} , \, \,  \overline{\nat }:= \nat \cup \{ \infty  \} , \, \, 
 {\rzecz }_{+}:=[0,\infty ) $.
Let  $(S, \scal )$ be a measurable space and let 
 ${M}_{\overline{\nat }}(S) $ be the set of all $\overline{\nat }$ valued measures on $(S, \scal )$.
On the set  ${M}_{\overline{\nat }}(S)$ we consider the $\sigma $-field
 $  {\mcal }_{\overline{\nat }}(S) $ defined as
the smallest $\sigma $-field  such that for all $ B \in \scal $: the map
$$
      {i}_{B} : {M}_{\overline{\nat }}(S) \ni \mu \mapsto \mu (B) \in \overline{\nat }
$$ 
is measurable.

\bigskip  \noindent
Let $(\Omega , \fcal ,\p  )$ be a complete probability space with filtration $\mathbb{F}:=({\fcal }_{t}{)}_{t\ge 0}$ satisfying the usual hypotheses, see \cite{Metivier_82}.
\bigskip 
\begin{definition} \rm (see Appendix C in \cite{Brzezniak_Hausenblas_2009}). Let $(Y, \ycal )$ be a measurable space. A
\bf  time homogeneous Poisson random measure $\eta $  \rm on $(Y, \ycal )$ over $(\! \Omega ,\fcal , \fmath ,\p )$ is a measurable function 
$$
\eta : (\Omega , \fcal ) \to \bigl( {M}_{\overline{\nat }} ({\rzecz }_{+}\times Y),{\mcal }_{\overline{\nat }} ({\rzecz }_{+}\times Y) \bigr) 
$$
such that 
\begin{itemize}
\item[(i) ] for all $ B \in \bcal ({\rzecz }_{+}) \otimes \ycal $, $\eta (B):= {i}_{B} \circ \eta : \Omega \to \overline{\nat }$ is a Poisson random measure with parameter $\e [\eta (B)]$;
\item[(ii) ] $\eta $ is independently scattered, i.e. if the sets ${B}_{j}\in \bcal ({\rzecz }_{+}) \otimes \ycal $, $j=1,...,n$, are disjoint then the random variables $\eta ({B}_{j})$, $j=1,...,n$, are independent;
\item[(iii) ]  for all $ U \in \ycal $  the $\overline{\nat }$-valued process 
$\bigl( N(t,U) {\bigr) }_{t \ge 0} $ defined by 
$$
   N(t,U):= \eta ((0,t]\times U) , \qquad  t \ge 0 
$$
is $\fmath $-adapted and its increments are independent of the past, i.e. if $t>s\ge 0$, then 
$N(t,U)-N(s,U)= \eta ((s,t]\times U)$ is independent of ${\fcal }_{s}$.
\end{itemize} 
\end{definition}
\noindent
If $\eta $ is a time homogeneous Poisson random measure then the formula
$$
  \nu (A) := \e [\eta ((0,1]\times A )] , \qquad A \in \ycal 
$$ 
defines a measure on $(Y, \ycal )$ called an \bf intensity measure \rm of $\eta $.
Moreover, for all $T<\infty $ and all $A\in \ycal $ such that $\e \bigl[ \eta ((0,T]\times A) \bigr] <\infty $, the 
$\rzecz $-valued process $\{ \tilde{N} (t,A) {\} \! }_{t \in (0,T]\! }$ defined by 
$$
   \tilde{N} (t,A) := \eta ((0,t]\times A)  - t \nu (A) , \qquad t \in (0,T],
$$ 
is an integrable martingale on $(\Omega ,\fcal , \fmath ,\p )$.
The random measure $l \otimes \nu $ on $\bcal ({\rzecz }_{+}) $ $\otimes $ $\ycal $, where $l$ stands for the Lebesgue measure, is called an \bf compensator \rm of $\eta $ and 
the difference between a time homogeneous Poisson random measure $\eta $ and its compensator, i.e. 
$$
    \tilde{\eta } := \eta - l \otimes \nu   ,
$$
is called a \bf compensated time homogeneous Poisson random measure.\rm 

\bigskip  \noindent
Let us also recall basic properties of the stochastic integral with respect to  $\tilde{\eta }$,
see \cite{Brzezniak_Hausenblas_2009}, \cite{Ikeda_Watanabe_81} and \cite{Peszat_Zabczyk_2007} for details.
Let $\hmath $ be a separable Hilbert space and let $\pcal $ be a predictable $\sigma $-field on $[0,T] \times \Omega $.
Let ${\mathfrak{L}}^{2}_{\nu ,T} (\pcal \otimes \ycal ,l \otimes \p \otimes \nu ;\hmath )$ be a space of all
$\hmath $-valued, $\pcal \otimes \ycal $-measurable processes such that
$$
  \e \Bigl[ \int_{0}^{T}\int_{Y} \norm{\xi (s, \cdot ,y)}{\hmath }{2} \,  ds d\nu (y) \Bigr] < \infty .
$$ 
If $\xi \in {\mathfrak{L}}^{2}_{\nu ,T} (\pcal \otimes \ycal ,l \otimes \p \otimes \nu ;\hmath )$
then the integral process $\int_{0}^{t} \int_{Y} \xi (s, \cdot ,y) \, \tilde{\eta }(ds,dy)$, 
$t\in [0,T]$, is a \it c\`{a}dl\`{a}g \rm ${L}^{2}$-integrable martingale. Moreover, the following isometry formula holds
\begin{equation} \label{E:isometry}
  \e \biggl[ \Norm{\int_{0}^{t} \int_{Y} \xi (s, \cdot ,y)  \tilde{\eta }(ds,dy) }{\hmath }{2} \biggr]
  =\e \Bigl[ \int_{0}^{t}\int_{Y} \norm{\xi (s, \cdot ,y)}{\hmath }{2}   ds d\nu (y) \Bigr] ,
   \, \,  t \in [0,T].
\end{equation}

\subsection{Statement of the problem} \label{S:Statement}  
\noindent
Problem (\ref{E:NS_intr})-(\ref{E:initial}) can be written as the following stochastic evolution equation
\begin{eqnarray} 
& du(t) & +\bigl[ \acal u(t) +B\bigl( u(t) \bigr) \bigr] \, dt = f(t) \, dt 
 + \int_{Y} F(t,u(t);y) \tilde{\eta } (dt,dy) \nonumber \\
& &\qquad \qquad + G(t,u(t))\, dW(t),  \qquad t \in [0,T] ,\nonumber  \\
& u(0) &= {u}_{0} .   \label{E:NS}
\end{eqnarray}

\bigskip \noindent
\bf Assumptions. \rm  We assume that
\begin{itemize}
\item[(A.1)] ${u}_{0} \in H$, $f \in {L}^{2}([0,T];V')$,
\item[(F.1)] $\tilde{\eta }$ is a compensated time homogeneous Poisson random measure on a measurable space $(Y, \ycal )$ over
$(\Omega , \fcal , \fmath , \p )$ with a $\sigma $-finite intensity measure $\nu $,
\item[(F.2)] $F:[0,T]\times H \times Y \to H $ is  a measurable function such that 

\noindent
$\int_{Y} \ind{\{ 0\} } (F(t,x;y))  \nu (dy) =0 $ for all $x \in H $ and $t \in [0,T]$. Moreover,
there exists a constant $L$ such that
\begin{equation}
    \int_{Y} |F(t,{u}_{1};y)- F(t,{u}_{2};y) {|}_{H}^{2}  \nu (dy) \le  L |{u}_{1}-{u}_{2}{|}_{H}^{2} 
   , \quad {u}_{1}, \! {u}_{2} \in H , \, \,  t \in [0,T] \label{E:F_Lipschitz_cond} ,
\end{equation}
and for each $p \in \{ 2,4,4+\gamma ,8+2\gamma \} $ there exists a constant ${C}_{p}$ such that
\begin{equation}   
    \int_{Y} |F(t,u;y) {|}_{H}^{p} \, \nu (dy) \le  {C}_{p} (1 + |u{|}_{H}^{p}), \qquad  u \in H , \quad t \in [0,T],  \label{E:F_linear_growth}
\end{equation}
where $\gamma >0$ is some positive constant.
\item[(F.3)] Moreover, for all $v \in \vcal $ the mapping ${\tilde{F}}_{v}$ defined by 
\begin{equation} \label{E:F**}
     \bigl( {\tilde{F}}_{v}(u)\bigr) (t,y):= \ilsk{F(t,u({t}^{-});y)}{v}{H}, 
     \quad u \in {L}^{2}(0,T;H), \quad (t,y) \in [0,T] \times Y 
\end{equation}
is a continuous from ${L}^{2}(0,T;H) $ into $ {L}^{2}([0,T]\times Y, dl\otimes \nu ; \rzecz ) $
 if in the space ${L}^{2}(0,T;H)$ we consider the Fr\'{e}chet topology inherited from the space 

\noindent
${L}^{2}(0,T;{H}_{loc})$.
 \footnote{Here $l$ denotes the Lebesgue measure on the interval $[0,T]$.} 
\item[(G.1)] $ W(t)$  is a cylindrical  Wiener process in a separable Hilbert space ${Y}_{W}$ defined on the stochastic basis $\bigl( \Omega , \fcal , \fmath  , \p  \bigr) $ ; 
\item[(G.2)]  $G: [0,T] \times V \to \lhs ({Y}_{W},H) $ and there exists a constant ${L}_{G}>0$
 such that 
\begin{equation}
   \norm{G(t,{u}_{1}) - G(t,{u}_{2})}{\lhs ({Y}_{W},H)}{2} \le {L}_{G} \norm{{u}_{1}-{u}_{2}}{V}{2} ,
   \quad {u}_{1}, {u}_{2} \in V , \, \,  t \in [0,T] .
\end{equation} 
Moreover there exist ${\lambda }_{}$, $\kappa \in \rzecz $ and $a\in \bigl( 2-\frac{2}{3+\gamma },2\bigr]$  such that
\begin{equation} \label{E:G}
     2 \dual{\acal u }{u }{} -  \norm{G(t,u )}{\lhs ({Y}_{W},H)}{2}
     \ge  a \norm{u}{}{2} -{\lambda }_{} {|u |}_{H}^{2} - \kappa  , \quad u \in V  , \, \, t \in [0,T].
\end{equation}
\item[(G.3)]
Moreover, $G $ extends to a continuous mapping $G :[0,T] \times H \to \lhs \! (\! {Y}_{W}, \! {V'}) $  such that
\begin{equation} \label{E:G*}
   \norm{G(t,u)}{\lhs ({Y}_{W}, {V'})}{2} \le C (1 + {|u|}_{H}^{2}) , \qquad u \in H . 
\end{equation}
for some $C>0$. Moreover, for every $v \in \vcal $ the mapping ${\tilde{G}}_{v}$ defined by
\begin{equation} \label{E:G**}
  \bigl( {\tilde{G}}_{v}(u)\bigr)  (t) := \ilsk{G(t,u(t))}{v}{H} ,
 \qquad u \in {L}^{2}(0,T;H) , \quad t \in [0,T]
\end{equation}
is a continuous mapping from ${L}^{2}(0,T;H) $ into $ {L}^{2}([0,T];\lhs ({Y}_{W},\rzecz ) ) $
 if in the space 
${L}^{2}(0,T;H)$ we consider the Fr\'{e}chet topology inherited from the space 

\noindent
${L}^{2}(0,T;{H}_{loc})$.
\end{itemize}
Let us recall that the space 
${L}^{2}(0,T;{H}_{loc})$ is defined by (\ref{E:seminorms}). For any Hilbert space $E$ the symbol $\lhs ({Y}_{W};E)$ denotes the space of Hilbert-Schmidt operators from ${Y}_{W}$ into  $E$. 

\bigskip
\begin{definition}  \rm  \label{D:solution}
 \bf A martingale solution \rm of  equation (\ref{E:NS})
is a system 

\noindent
$\bigl( \bar{\Omega }, \bar{\fcal },  \bar{\fmath } ,\bar{\p },\bar{u}, \bar{\eta }, \bar{W}\bigr) $,
where
\begin{itemize}
\item[$\bullet $]  $\bigl( \bar{\Omega }, \bar{\fcal },  \bar{\fmath } ,\bar{\p }  \bigr) $ is a filtered probability space with a filtration $\bar{\fmath } = \{ {\bar{\fcal }_{t}}{\} }_{t \ge 0} $, 
\item[$\bullet $] $\bar{\eta }$  is a time homogeneous Poisson random measure on $(Y, \ycal )$ over
$\bigl( \! \bar{\Omega }, \bar{\fcal },  \bar{\fmath } ,\bar{\p } \! \bigr) $ with the intensity measure $\nu $,
\item[$\bullet $] $\bar{W}$ is a cylindrical Wiener process on the space ${Y}_{W}$ over
$\bigl( \bar{\Omega }, \bar{\fcal },  \bar{\fmath } ,\bar{\p }  \bigr) $,
\item[$\bullet $] $\bar{u}: [0,T] \times \Omega \to H$ is a predictable process with $\bar{\p } $ - a.e. paths
$$
  \bar{u}(\cdot , \omega ) \in \dmath \bigl( [0,T], {H}_{w} \bigr)
   \cap {L}^{2}(0,T;V )
$$
such that for all $ t \in [0,T] $ and all $v \in V $ the following identity holds $\bar{\p }$ - a.s.
\begin{eqnarray*}
&\ilsk{\bar{u}(t)}{v}{H} & + \int_{0}^{t} \dual{\acal \bar{u}(s)}{v}{}  ds
+ \int_{0}^{t} \dual{B(\bar{u}(s))}{v}{}  ds 
 \nonumber \\
& =\ilsk{{u}_{0}}{v}{H} &+ \int_{0}^{t} \dual{f(s)}{v}{} ds 
 + \int_{0}^{t} \int_{Y} \ilsk{F(s,\bar{u}(s);y)}{v}{H} \,  \tilde{\bar{\eta }} (ds,dy) \nonumber  \\
& & + \dual{\int_{0}^{t}G(s,\bar{u}(s))\,  d\bar{W}(s)}{v}{} . 
\end{eqnarray*}
\end{itemize}
\end{definition}

\bigskip  \noindent
We will prove existence of a martingale solution of the equation (\ref{E:NS}).
To this end we use the Faedo-Galerkin method. The Galerkin approximations generate a sequence of probability measures on appropriate functional space. We will prove that this sequence is tight. Let us emphasize that to prove the tightness,  assumption (F.2)  with $p=2$ in inequality (\ref{E:F_linear_growth}) is sufficient. The stronger condition on $p$, i.e. inequality (\ref{E:F_linear_growth}) for a certain $p>4$, is connected with the construction of the process $\bar{u}$ to deal with the nonlinear term.  
Assumptions (G.2)-(G.3) allow to consider the Gaussian noise term $G$ dependent both on $u$ and $\nabla u$.
This corresponds to inequality (\ref{E:G}) with  $a<2$. The case when $a=2$ is related to the noise term $G$ dependent on $u$ but not on its gradient.
Moreover, assumptions (F.3) and (G.3) are important in the case of unbounded domain $\ocal $. In the case when $\ocal $ is bounded, they can be omitted, see \cite{Motyl_NS_Poisson_pre}.

\section{Existence of solutions} \label{S:Existence}

\bigskip 
\begin{theorem} \label{T:existence}  
There exists a martingale solution of the problem (\ref{E:NS}) provided assumptions (A.1), (F.1)-(F.3)
and (G.1)-(G.3) are satisfied.
\end{theorem}

\subsection{Faedo-Galerkin approximation}

\bigskip 
\noindent
Let $\{ {e}_{i} {\} }_{i =1}^{\infty  }$ be the orthonormal basis in $H$ composed of eigenvectors of the operator $L$ defined by (\ref{E:op_L}).
Let ${H}_{n}:= span \{ {e}_{1}, ..., {e}_{n} \} $ be the subspace with the norm inherited from $H$ and
let $\Pn : H \to {H}_{n} $ be defined by (\ref{E:P_n}).
Let us fix $m > \frac{d}{2}+1$ and let $U$ be the space defined by (\ref{E:U_comp_V_m}).
Consider the following mapping
$$
  \Bn (u):= \Pn B({\chi }_{n}(u),u) , \qquad u \in {H}_{n},
$$
where ${\chi }_{n}:H \to H $ is defined by ${\chi }_{n}(u) = {\theta }_{n}(|u {|}_{U'})u$, where  ${\theta }_{n } : \rzecz \to [0,1]$  of class ${\ccal }^{\infty }$ such that
\begin{eqnarray*}
 {\theta }_{n}(r)  = 1 \quad \mbox{if} \quad  r \le n  \quad \mbox{ and } \quad 
  {\theta }_{n}(r)  = 0 \quad \mbox{if} \quad  r \ge n+1 .
\end{eqnarray*}
Since ${H}_{n} \subset H$, ${B}_{n}$ is well defined. Moreover, ${B}_{n}:{H}_{n} \to {H}_{n}$ is globally Lipschitz continuous.

\bigskip  \noindent
Let us consider the classical Faedo-Galerkin approximation in the space $ {H}_{n}$
\begin{eqnarray} 
&    \un (t) & =  \Pn {u}_{0} - \int_{0}^{t}\bigl[ \Pn \acal \un (s)  + \Bn  \bigl(\un (s) \bigr) - \Pn f (s)  \bigr] \, ds  \nonumber   \\ 
&  &+  \int_{0}^{t} \int_{Y} \Pn F(s,\un ({s}^{-}),y) \tilde{\eta } (ds,dy) \nonumber \\
&  & + \int_{0}^{t} \Pn G(s,\un (s)) \, dW(s) ,   
\quad t \in [0,T]  .   \label{E:Galerkin}
\end{eqnarray}

\bigskip
\begin{lemma} \label{L:Galerkin_existence}
For each $n \in \nat $,  there exists a unique  $\fmath $-adapted,  c\`{a}dl\`{a}g ${H}_{n}$ valued process  ${u}_{n}$ satisfying the Galerkin equation (\ref{E:Galerkin}).
\end{lemma}

\proof
The assertion follows from Theorem 9.1 in \cite{Ikeda_Watanabe_81}. \qed

\bigskip  
\noindent
Using the It\^{o} formula, see \cite{Ikeda_Watanabe_81} or \cite{Metivier_82},
and the Burkholder-Davis-Gundy inequality, see \cite{Peszat_Zabczyk_2007}, we will prove the following lemma about \it a priori \rm estimates of the solutions ${u}_{n} $ of (\ref{E:Galerkin}).
In fact, these estimates hold provided the noise terms satisfy only condition (\ref{E:F_linear_growth}) in assumption (F.2) and condition (\ref{E:G}) in assumption (G.2). 

\bigskip
\begin{lemma} \label{L:Galerkin_estimates }
The processes $({u}_{n} {)}_{n \in \nat }$ satisfy the following estimates.
\begin{itemize}
\item[(i) ]
For every $p\in [1,4+\gamma] $  there exists a  positive constant ${C}_{1}(p)$  such that
\begin{equation} \label{E:H_estimate}
 \sup_{n \ge 1 } \e \bigl( \sup_{0 \le s \le T } |\un (s){|}_{H}^{p} \bigr) \le {C}_{1}(p) .
\end{equation}
\item[(ii)] There exists a positive constant ${C}_{2}$ such that
\begin{equation} \label{E:V_estimate}
  \sup_{n \ge 1 }  \e \bigl[ \int_{0}^{T} \norm{ \un (s)}{V}{2} \, ds \bigr] \le {C}_{2}.
\end{equation}
\end{itemize}
\end{lemma}
\noindent
Let us recall that $\gamma >0$ is defined in assumption (F.2).

\proof 
For all $n \in \nat $ and all $R>0$ let us define
\begin{equation} \label{E:stopping_time}
  {\taun }_{}(R):= \inf \{ t \ge 0: \, \, |{u}_{n}(t){|}_{H} \ge R  \} \wedge T. 
\end{equation}
Since the process $\bigl( {u}_{n}(t) {\bigr) }_{t \in [0,T]} $ is $\fmath $-adapted and right-continuous, 
$\taun (R)$ is a stopping time. Moreover, since  the process $({u}_{n})$ is c\`{a}dl\`{ag} on $[0,T]$,  the trajectories $t \mapsto {u}_{n}(t)$ are bounded on $[0,T]$, $\p $-a.s. Thus
$\taun (R)\uparrow T$, $\p $-a.s., as $R\uparrow \infty $.

\bigskip  \noindent
Assume first that $p =2$ or $p=4+\gamma $.
Using the It\^{o} formula to the function  $\phi(x):=|x{|}^{p}:= |x{|}_{H}^{p}$, $x \in H$,  we obtain for all $t \in [0,T]$
\begin{eqnarray*}
 &\bigl|  {u}_{n} & ({t\wedge \taun (R)}) {\bigr| }^{p} =  |\Pn {u}_{0}{|}^{p} \\
&+ &\int_{0}^{t\wedge \taun (R)} \bigl\{ p \bigl| \un (s){\bigr| }^{p-2} \dual{\un (s)}{- \Pn \acal \un (s) - {B}_{n}(\un (s)) + \Pn f (s)}{} \bigr\}  \, ds \\
 & +& \int_{0}^{t\wedge \taun (R)} \int_{{Y}} \big\{ \phi \bigl( \un ({s}^{-})  + \Pn F(s,\un ({s}^{-});y) \bigr) -
   \phi \bigl( \un ({s}^{-}) \bigr)  \bigr\} \, \tilde{\eta}(ds,dy)   \\
 & + &\int_{0}^{t\wedge \taun (R)} \int_{{Y}} \big\{ \phi \bigl( \un ({s}^{-})  + \Pn F(s,\un ({s}^{-});y) \bigr) -
   \phi \bigl( \un ({s}^{-}) \bigr)  \\
 &  - &\dual{\phi '(\un ({s}^{-}) ) }{\Pn F(s,\un ({s}^{-});y)}{} \bigr\} \, \nu (ds,dy)  \\ 
  &  +& \frac{1}{2} \int_{0}^{t\wedge \taun (R)} 
\tr \bigl[ \Pn G(s,\un (s)) \frac{{\partial }^{2} \phi }{\partial  {x}^{2}}
     { \bigl( \Pn G(s,\un (s)) \bigr) }^{*} \bigr]  \, ds \\
 & +&\int_{0}^{t\wedge \taun (R)}   p \,  {|\un (s)|}^{p-2} 
  \dual{\un (s)}{ \Pn G(s, \un (s) ) \, d W(s) }{} .
\end{eqnarray*} 
By (\ref{E:op_Acal}) and (\ref{E:wirowosc_B}) we obtain for all $t \in [0,T]$
\begin{eqnarray*} 
& \bigl| {u}_{n} & ({t\wedge \taun (R)}) {\bigr| }^{p} =   
  |\Pn {u}_{0}{|}^{p} \nonumber \\
&+& \int_{0}^{t\wedge \taun (R)} \bigl\{ -p \bigl| \un (s){\bigr| }^{p-2} \norm{\un (s)}{}{2} 
  + p \bigl| \un (s){\bigr| }^{p-2} \dual{\un (s)}{  f (s)}{} \bigr\}  \, ds  \nonumber \\  
 & + &\int_{0}^{t\wedge \taun (R)} 
   \int_{{Y}} \big\{  \bigl| \un ({s}^{-})  + \Pn F(s,\un ({s}^{-});y) {\bigr| }^{p} -
    \bigl| \un ({s}^{-}) {\bigr|}^{p}  \bigr\} \, \tilde{\eta}(ds,dy)  \nonumber  \\
 & +& \int_{0}^{t\wedge \taun (R)} \int_{{Y}} 
   \big\{  \bigl| \un ({s}^{-})  + \Pn F(s,\un ({s}^{-});y) {\bigr| }^{p}-
   \bigl| \un ({s}^{-} ) {\bigr| }^{p}  \nonumber \\
 &  - &p \bigl| \un ({s}^{-}) {\bigr| }^{p-2} \dual{\un ({s}^{-})  }{\Pn F(s,\un ({s}^{-});y)}{} \bigr\} \, 
 \nu (dy) ds   \nonumber \\
  &  + &\frac{1}{2} \int_{0}^{t\wedge \taun (R)} 
\tr \bigl[ \Pn G(s,\un (s)) \frac{{\partial }^{2} \phi }{\partial  {x}^{2}}
     { \bigl( \Pn G(s,\un (s)) \bigr) }^{*} \bigr]  \, ds \nonumber \\
 & +&\int_{0}^{t\wedge \taun (R)}   p \,  {|\un (s)|}^{p-2} 
  \dual{\un (s)}{ G(s, \un (s) ) \, d W(s) }{}.
\end{eqnarray*}
Let us recall that according to (\ref{E:op_Acal}) we have  $\dual{\acal u}{u}{} = \dirilsk{u}{u}{}$ and thus   
$$2\dual{\acal u}{u}{}-a \norm{u}{}{2}=(2-a ) \norm{u}{}{2}.$$
Hence  inequality (\ref{E:G}) in assumption (G.2) can be written equivalently in the following form
$$ \label{A:G'}
 \norm{G(s,u )}{\lhs ({Y}_{W},H)}{2}
  \le (2-a ) \norm{u}{}{2}+{\lambda }_{}{|u|}_{H}^{2}+ \kappa  , \qquad u \in V , \quad
   s \in [0,T]. 
$$
Hence
\begin{eqnarray*}
& &\frac{1}{2}  \int_{0}^{t\wedge \taun (R)} 
\tr \bigl[ \Pn G(s,\un (s)) \frac{{\partial }^{2} \phi }{\partial  {x}^{2}}
     { \bigl( \Pn G(s,\un (s)) \bigr) }^{*} \bigr]  \, ds  \\
& &  \le  \frac{p(p-1)}{2} \int_{0}^{t\wedge \taun (R)}  |\un (s){|}^{p-2} 
   \bigl[  (2-a ) \norm{\un (s)}{}{2}+{\lambda }_{}{|\un (s)|}_{}^{2}+ \kappa \bigr] \, ds  .  
\end{eqnarray*}
Moreover, by assumption (A.1), (\ref{E:norm_V})  and the Schwarz inequality, we obtain for every 
$\eps >0$ and  for all $s \in [0,T] $
\begin{eqnarray*}
& & \,  \dual{f(s)}{\un (s)}{}
  \le  \, {|f(s)|}_{V'} \cdot \norm{\un (s)}{V}{}
    \\
& &\le {|f(s)|}_{V'} \, |\un (s){|}_{}  +   \frac{1}{ 4\eps  } {|f(s)|}_{V'}^{2} +  \eps  \, \norm{\un (s)}{}{2}  
\end{eqnarray*}
and hence by the Young inequality  
\footnote{$ab \le \frac{1}{{q}_{1}} {a}^{{q}_{1}} 
 + \frac{1}{{q}_{2}}  {b}^{{q}_{2}}$ if $a,b >0 $, ${q}_{1},{q}_{2} \in (1,\infty )$ and $\frac{1}{{q}_{1}} + \frac{1}{{q}_{2}}=1$.} 
\begin{eqnarray*}
& & |\Pn {u}_{0}{|}^{p} + p {|\un (s)|}^{p-2} 
\Bigl(  {|f(s)|}_{V'}  |\un (s){|}_{}  +\frac{1}{ 4\eps  } {|f(s)|}_{V'}^{2}  \Bigr) \\
& &  + \frac{p(p-1)}{2} {|\un (s)|}^{p-2}  \bigl[
  {\lambda }_{} {|\un (s)|}^{2} + \kappa   \bigr]     
 =  \frac{p(p-1)\lambda }{2}  \, {|\un (s)|}^{p} \\
& & +  |\Pn {u}_{0}{|}^{p} +p {|f(s)|}_{V'} |\un (s){|}^{p-1} 
+ p\Bigl(  \frac{1}{ 4\eps  } {|f(s)|}_{V'}^{2} + \frac{(p-1)\kappa }{2}\Bigr) {|\un (s)|}^{p-2}  \\
& & \le c + {c}_{1} {|\un (s)|}^{p} 
\end{eqnarray*}
for some constants $c,{c}_{1}>0$.
Thus
\begin{eqnarray} 
& \bigl| {u}_{n} &  ({t\! \wedge \! \taun (R)}) {\bigr| }^{p}
  + \bigl[ p - p \eps - \frac{1}{2} p(p-1) (2- a)   \bigr]
 \int_{0}^{t\wedge \taun (R)} \! \! \! \bigl| \un (s) {\bigr| }^{p-2} \norm{\un (s)}{}{2}  ds \nonumber \\
& \le & c 
+ {c}_{1}  \! \int_{0}^{t\wedge \! \taun \! (R)} \bigl| \un (s){\bigr| }^{p}  \, ds  \nonumber \\
& & + \int_{0}^{t\wedge \taun (R)} 
   \int_{{Y}} \big\{  \bigl| \un ({s}^{-})  +\Pn F(s,\un ({s}^{-});y) {\bigr| }^{p} -
    \bigl| \un ({s}^{-}) {\bigr|}^{p}  \bigr\} \, \tilde{\eta}(ds,dy)  \nonumber  \\
& & + \int_{0}^{t\wedge \taun (R)} \int_{{Y}} 
   \big\{  \bigl| \un ({s}^{-})  + \Pn F(s,\un ({s}^{-});y) {\bigr| }^{p}-
   \bigl| \un ({s}^{-} ) {\bigr| }^{p}  \nonumber  \\
& &  - p \bigl| \un ({s}^{-}) {\bigr| }^{p-2} \ilsk{\un ({s}^{-})  }{\Pn F(s,\un ({s}^{-});y)}{H} \bigr\} \, 
\nu (dy) ds
\nonumber \\
& & +\int_{0}^{t\wedge \taun (R)}   p \,  {|\un (s)|}^{p-2} 
  \dual{\un (s)}{ G(s, \un (s) ) \, d W(s) }{} . \label{E:Ito_formula}
\end{eqnarray} 
Let us choose $\eps > 0 $ such that  $ p - p \eps - \frac{1}{2} p(p-1) (2- a ) >0 $,
or equivalently,
$$
  \eps < 1 - \frac{1}{2} (p-1) (2-a ).
$$
Note that since by assumption (G.2) $a\in \bigl( 2-\frac{2}{3+\gamma } , 2] $,  such an $\eps $ exists.

\bigskip  \noindent
From the Taylor formula, it follows that for each $p\ge 2 $ there exists a positive constant ${c}_{p}>0$ such that for all $x,h \in H$ the following inequality holds
\begin{eqnarray}
 &  &\bigl| |x+h{|}_{H}^{p} -|x{|}_{H}^{p}- p |x{|}_{H}^{p-2} \ilsk{x}{h}{H} \bigr| 
  \le {c}_{p} ( {|x|}_{H}^{p-2} +  {|h|}_{H}^{p-2} ) \, |h{|}_{H}^{2} .  \label{E:Taylor_II}
\end{eqnarray}
By (\ref{E:Taylor_II}), (\ref{E:F_linear_growth}) and (\ref{E:stopping_time}), the process 
$\bigl( {M}_{n}(t\wedge \taun (R)) {\bigr) }_{t \in [0,T]}$, where
$$
 {M}_{n}(t):=  \int_{0}^{t} 
   \int_{{Y}} \big\{  \bigl| \un ({s}^{-})  + \Pn F(s,\un ({s}^{-});y) {\bigr| }^{p} -
    \bigl| \un ({s}^{-}) {\bigr|}^{p}  \bigr\} \, \tilde{\eta}(ds,dy) , 
$$
$ t\in [0,T]$,
is an integrable martingale. Hence $\e [{M}_{n}(t\wedge {\tau }_{n}(R))] = 0 $ for all $t\in [0,T]$.
By (\ref{E:G}) and (\ref{E:stopping_time}), the process 
$\bigl( {N}_{n}(t\wedge \taun (R)) {\bigr) }_{t \in [0,T]}$, where
$$
{N}_{n}(t):= \int_{0}^{t} {|\un (s)|}^{p-2} \dual{\un (s)}{G(s,\un (s)) \, d W(s) }{} ,
 \qquad t \in [0,T]
$$
is an integrable martingale and thus $\e [{N}_{n} (t\wedge {\tau }_{n}(R)) ] = 0 $ for all $t\in [0,T]$.

\bigskip  \noindent
Let us denote
\begin{eqnarray}  
& & {I}_{n}(t):= \int_{0}^{t} \int_{{Y}} 
   \big\{  \bigl| \un ({s}^{-})  + \Pn F(s,\un ({s}^{-});y) {\bigr| }^{p}
   -\bigl| \un ({s}^{-}) {\bigr| }^{p}  \nonumber  \\
& &  - p \bigl| \un ({s}^{-}) {\bigr| }^{p-2} \ilsk{\un ({s}^{-})  }{\Pn F(s,\un ({s}^{-});y)}{H} \bigr\} \, 
\nu (dy) ds  , \quad t\in [0,T].  \label{E:I_n}
\end{eqnarray}
By (\ref{E:Taylor_II})
and (\ref{E:F_linear_growth}) we obtain the following inequalities
\begin{eqnarray*}
& &\bigl| {I}_{n}(t) \bigr| \nonumber \\
& &\le  {c}_{p}\int_{0}^{t} \int_{{Y}} \bigl| \Pn F(s,\un ({s}^{-});y) {\bigr| }_{H}^{2}
 \bigl\{ \bigl| \un ({s}^{-}) {\bigr| }_{H}^{p-2} + \bigl| \Pn F(s,\un ({s}^{-});y) {\bigr| }_{H}^{p-2} \bigr\}
  \nu (dy) ds \nonumber \\
& &\le  {c}_{p}\int_{0}^{t} \bigl\{ 
{C}_{2} \bigl| \un ({s}^{}) {\bigr| }_{H}^{p-2} \bigl( 1+ \bigl| \un ({s}^{}) {\bigr| }_{H}^{2}\bigr) 
 + {C}_{p}\bigl( 1+ \bigl| \un ({s}^{}) {\bigr| }_{H}^{p}\bigr)   \bigr\} \, ds \nonumber \\
& &\le  {\tilde{c}}_{p} \int_{0}^{t} \bigl\{ 1+ \bigl| \un ({s}^{}) {\bigr| }_{H}^{p} \bigr\} \, ds
= {\tilde{c}}_{p}t + {\tilde{c}}_{p}\int_{0}^{t} \bigl| \un ({s}^{}) {\bigr| }_{H}^{p} \, ds , \qquad t \in [0,T] 
\end{eqnarray*}
for some constant ${\tilde{c}}_{p}>0$.
Thus by the Fubini Theorem, we obtain the following inequality
\begin{equation} \label{E:E_I_n(t)}
 \e \bigl[ \bigl| {I}_{n}(t) \bigr| \bigr] \le  
 {\tilde{c}}_{p}t + {\tilde{c}}_{p} \int_{0}^{t} \e \bigl[ \bigl| \un (s) {\bigr| }_{H}^{p} \bigr] \, ds , \qquad t \in [0,T]. 
\end{equation}
By (\ref{E:Ito_formula}) and (\ref{E:E_I_n(t)}), we have for all $t\in [0,T]$
\begin{eqnarray} 
& &\e \bigl[  \bigl| {u}_{n}  ({t\wedge \taun (R)}) {\bigr| }_{H}^{p} \bigr] \nonumber  \\
 & &+\bigl[ p - p \eps - \frac{1}{2} p(p-1) (2- a) \,   \bigr]
 \e \Bigl[  \int_{0}^{T\wedge \taun (R)}  \bigl| \un (s){\bigr| }_{H}^{p-2} \norm{\un (s)}{}{2} \, ds \Bigr] 
 \nonumber \\
\qquad & &\le c+ {\tilde{c}}_{p}T + ({c}_{1} +{\tilde{c}}_{p} ) \int_{0}^{t\wedge \taun (R)} \e \bigl[ \bigl| \un (s) {\bigr| }_{H}^{p} \bigr] \, ds .
\label{E:Ito_formula_est}
\end{eqnarray} 
In particular,
$$ 
\e \bigl[  \bigl| {u}_{n}  ({t\wedge \taun (R)}) {\bigr| }_{H}^{p} \bigr]
 \le  c+{\tilde{c}}_{p}T + ({c}_{1}+{\tilde{c}}_{p}) \int_{0}^{t\wedge \taun (R)} \e \bigl[ \bigl| \un (s) {\bigr| }_{H}^{p} \bigr] \, ds .
$$
By the Gronwall Lemma we infer that for all $t \in [0,T]$: $
\e \bigl[  \bigl| {u}_{n}  ({t\wedge \taun (R)}) {\bigr| }^{p} \bigr]
 \le   {\tilde{\tilde{C_p}}} $
for some constant ${\tilde{\tilde{C_p}}}$ independent of $t\in [0,T]$, $R>0$ and $n \in \nat $, i.e.
$$ \label{E:E(un(t))_R_est}
\sup_{n \ge 1} \sup_{t\in [0,T]} \e \bigl[  \bigl| {u}_{n} (t\wedge \taun (R)) {\bigr| }_{H}^{p} \bigr]
 \le   {\tilde{\tilde{C_p}}} .
$$
Hence, in particular,
$$ \label{E:E(int_un(t))_R_est}
\sup_{n \ge 1} \e \Bigl[ \int_{0}^{{T\wedge \taun (R)}}  \bigl| {u}_{n} (s) {\bigr| }_{H}^{p} \, ds \Bigr]  \le   \tilde{{C}_{p}} 
$$
for some constant $\tilde{{C}_{p}}>0$.
Passing to the limit as $R \uparrow \infty $, by the Fatou Lemma we infer that
\begin{equation} \label{E:E(int_un(t))_est}
\sup_{n \ge 1} \e \Bigl[ \int_{0}^{T}  \bigl| {u}_{n} (s) {\bigr| }_{H}^{p} \, ds \Bigr]  
 \le   \tilde{{C}_{p}} .
\end{equation} 
By (\ref{E:Ito_formula_est}) and (\ref{E:E(int_un(t))_est}), we infer that
\begin{equation}  \label{E:HV_R_estimate}
\sup_{n \ge 1} \e \Bigl[  \int_{0}^{T\wedge \taun (R)}  \bigl| \un (s){\bigr| }_{H}^{p-2} \norm{\un (s)}{}{2} \, ds \Bigr]  \le {C}_{p}
\end{equation}
for some positive constant ${C}_{p}$. 
Passing to the limit as $R \uparrow \infty $ and using again the Fatou Lemma we infer that
\begin{equation}  \label{E:HV_estimate}
\sup_{n \ge 1} \e \Bigl[  \int_{0}^{T}  \bigl| \un (s){\bigr| }_{H}^{p-2} \norm{\un (s)}{}{2} \, ds \Bigr] 
\le {C}_{p} .
\end{equation} 
In particular, putting $p:=2$ by (\ref{E:norm_V}), (\ref{E:HV_estimate}) and (\ref{E:E(int_un(t))_est})   we obtain assertion (\ref{E:V_estimate}). 

\bigskip  \noindent
Let us move to the proof of inequality (\ref{E:H_estimate}). By the Burkholder-Davis-Gundy inequality we obtain
\begin{eqnarray} 
& &\e \Bigl[ \sup_{r\in [0,t ]} |{M}_{n}(r\wedge \taun (R))| \Bigr] \nonumber \\
& &=\e \Bigl[ \sup_{r\in [0,t ]} 
 \Bigl| \! \int_{0}^{r\wedge \taun (R)} \! \!
   \int_{Y} \! \big\{  \bigl| \un ({s}^{-}\! )+\Pn F(s,\un ({s}^{-}\! );y) {\bigr| }_{H}^{p} -
    \bigl| \un ({s}^{-}\! ) {\bigr| }_{H}^{p}  \bigr\}  \tilde{\eta}(ds,dy) \Bigr|  \Bigr] \nonumber \\
& &\le {\tilde{K}}_{p} \e \Bigl[ \!
 \Bigl( \! \int_{0}^{t\wedge \taun (R)} \! \! \! \!
   \int_{Y} \bigl(  \bigl| \un ({s}^{-}\! )+\Pn F(s,\un ({s}^{-}\! );y) {\bigr| }_{H}^{p} -
    \bigl| \un ({s}^{-}) {\bigr|}_{H}^{p}  {\bigr) \! }^{2}  \nu (dy) ds {\Bigr) \! }^{\frac{1}{2}} 
 \! \Bigr]   \label{E:BDG_R}
\end{eqnarray}
for some constant ${\tilde{K}}_{p}>0$.
By (\ref{E:Taylor_II}) and the Schwarz inequality we obtain the following inequalities for all $x,h \in H$
\begin{eqnarray*}
& & \bigl( |x+h{|}_{H}^{p} - |x{|}_{H}^{p} {\bigr) }^{2} 
 \le 2\bigl\{ {p}^{2}|x{|}_{H}^{2p-2} |h{|}_{H}^{2} + {c}_{p}^{2} 
   \bigl( |x{|}_{H}^{p-2} + |h{|}_{H}^{p-2} {\bigr) }^{2}   |h{|}_{H}^{4} \bigr\}  \nonumber \\
& & \le 2  {p}^{2}|x{|}_{H}^{2p-2} |h{|}_{H}^{2} + 4{c}_{p}^{2} |x{|}_{H}^{2p-4} |h{|}_{H}^{4} 
+ 4{c}_{p}^{2}  |h{|}_{H}^{2p}.  
\end{eqnarray*}
Hence by inequality (\ref{E:F_linear_growth}) in assumption (F.2) we obtain for all $s \in [0,T] $
\begin{eqnarray} 
& &\int_{Y}\bigl( \big| \un ({s}^{-})  + \Pn F(s,\un ({s}^{-});y) {\bigr| }_{H}^{p} -
    \bigl| \un ({s}^{-}) {\bigr|}_{H}^{p}  {\bigr) }^{2} \nu (dy) \nonumber \\
& & \le   2{p}^{2} \bigl| \un ({s}^{-}) {\bigr|}_{H}^{2p-2} 
   \int_{Y} \bigl|  F(s,\un ({s}^{-});y) {\bigr| }_{H}^{2} \, \nu (dy) \nonumber \\ 
& &  + 4 {c}_{p}^{2} \bigl| \un ({s}^{-}) {\bigr|}_{H}^{2p-4} 
   \int_{Y} \bigl|  F(s,\un ({s}^{-});y) {\bigr| }_{H}^{4} \, \nu (dy)
 + 4 {c}_{p}^{2} \int_{Y}\bigl|  F(s,\un ({s}^{-});y) {\bigr| }_{H}^{2p} \nu (dy) \nonumber \\
& &  \le {C}_{1} + {C}_{2} \bigl| \un ({s}^{-}) {\bigr|}_{H}^{2p-4} 
  + {C}_{3} \bigl| \un ({s}^{-}) {\bigr|}_{H}^{2p-2} 
  + {C}_{4} \bigl| \un ({s}^{-}) {\bigr|}_{H}^{2p}
\label{E:BDG_R_Taylor_II'}  
\end{eqnarray}
for some positive constants ${C}_{i}$, $i=1,...,4$.
By (\ref{E:BDG_R_Taylor_II'}) and the Young inequality we infer that
$$
\int_{Y}\bigl( \big| \un ({s}^{-})  + \Pn F(s,\un ({s}^{-});y) {\bigr| }_{H}^{p} -
    \bigl| \un ({s}^{-}) {\bigr|}_{H}^{p}  {\bigr) }^{2} \nu (dy) 
    \le {K}_{1} + {K}_{2} \bigl| \un ({s}^{-}) {\bigr|}_{H}^{2p}
$$     
for some positive constants ${K}_{1}$ and ${K}_{2}$. Thus
\begin{eqnarray}
& & 
 \Bigl( \int_{0}^{t\wedge \taun (R)} 
   \int_{Y} \bigl( \big| \un ({s}^{-})  + \Pn F(s,\un ({s}^{-});y) {\bigr| }_{H}^{p} -
    \bigl| \un ({s}^{-}) {\bigr|}_{H}^{p}  {\bigr) }^{2}  \, \nu (dy) ds {\Bigr) }^{\frac{1}{2}} 
   \nonumber  \\
 & & \le \sqrt{T{K}_{1}} + \sqrt{{K}_{2}}  \Bigl( 
  \int_{0}^{t\wedge \taun (R)}  \bigl| \un ({s}^{}) {\bigr|}_{H}^{2p} \, ds
{\Bigr) }^{\frac{1}{2}}  . \label{E:BDG_R_Taylor_II} 
\end{eqnarray}
By (\ref{E:BDG_R}), (\ref{E:BDG_R_Taylor_II}) and (\ref{E:E(int_un(t))_est})
we obtain the following inequalities
\begin{eqnarray} 
& &\e \Bigl[ \sup_{r\in [0,t]} |{M}_{n}(r \wedge \taun (R))| \Bigr] \nonumber \\
& &\le {\tilde{K}}_{p}\sqrt{T{K}_{1}} + {\tilde{K}}_{p}\sqrt{{K}_{2}} \e \Bigl[ \Bigl( 
  \int_{0}^{t\wedge \taun (R)}  \bigl| \un ({s}^{}) {\bigr|}_{H}^{2p} \, ds
{\Bigr) }^{\frac{1}{2}} \Bigr]  \nonumber \\
& &\le {\tilde{K}}_{p}\sqrt{T{K}_{1}} + {\tilde{K}}_{p}\sqrt{{K}_{2}} \e \Bigl[ 
 \Bigl( \sup_{s\in [0,t]}  \bigl| \un ({s\wedge \taun (R)}^{}) {\bigr|}_{H}^{p} {\Bigr) }^{\frac{1}{2}}
 \Bigl( \int_{0}^{t\wedge \taun (R)}  \bigl| \un ({s}^{}) {\bigr|}_{H}^{p} \, ds 
{\Bigr) }^{\frac{1}{2}} \Bigr]   \nonumber \\
& &\le {\tilde{K}}_{p}\sqrt{T{K}_{1}} + \frac{1}{4} \e \Bigl[ 
  \sup_{s\in [0,t]}  \bigl| \un ({s\wedge \taun (R)}^{}) {\bigr|}_{H}^{p} \Bigr]
 + {\tilde{K}}_{p}^{2} {K}_{2} \e \Bigl[ \int_{0}^{t\wedge \taun (R)}  \bigl| \un ({s}^{}) {\bigr|}_{H}^{p} \, ds   \Bigr] \nonumber \\
& & \le \frac{1}{4} \e \Bigl[ 
  \sup_{s\in [0,t]}  \bigl| \un ({s\wedge \taun (R)}^{}) {\bigr|}_{H}^{p} \Bigr] + \tilde{K} ,
 \label{E:BDG_R_cont.}
\end{eqnarray}
where $\tilde{K} = {\tilde{K}}_{p}\sqrt{T{K}_{1}} + {\tilde{K}}_{p}^{2} {K}_{2} {\tilde{C}}_{p}$. (The constant ${\tilde{C}}_{p}$ is the same as in (\ref{E:E(int_un(t))_est})).

\bigskip  \noindent
Similarly, by the Burkholder-Davis-Gundy inequality  we obtain
\begin{eqnarray*}
& &\e \Bigl[ \sup_{r\in [0,t ]} |{N}_{n}(r\wedge \taun (R))| \Bigr] \\
& & =\e \Bigl[ \sup_{r\in [0,t ]}
 \Bigl| \int_{0}^{r\wedge \taun (R)} p \, {|\un (s )|}^{p-2} \dual{\un (s )}{\Pn G(s , \un (s ) ) \, d W(s ) }{}
  \Bigr| \Bigr] \\
& & \le C \, p \cdot \e \Bigl[
 {\Bigl( \int_{0}^{t \wedge \taun (R)} \, { |\un (s )|}^{2p-2} \cdot
 \norm{  G(s , \un (s ) ) }{\lhs (Y,H)}{2}
   \, ds
  \Bigr) }^{\frac{1}{2}} \Bigr]   \\
& & \le C  p  \e \Bigl[ \Bigl( \sup_{ s\in  [0, t] } {|\un (s\wedge {\tau }_{n}(R))|}^{p} 
{\Bigr) }^{\frac{1}{2}}
 {\Bigl( \int_{0}^{t \wedge \taun (R)}   { |\un (s )|}^{p-2} 
 \norm{  G(s , \un (s ) ) }{\lhs (Y,H)}{2}
    ds
  \Bigr) \!  }^{\frac{1}{2}} \Bigr]  .
\end{eqnarray*} 
By inequality (\ref{E:G}) in assumption (G.2) and
estimates (\ref{E:HV_estimate}), (\ref{E:E(int_un(t))_est})  we have the following inequalities
\begin{eqnarray}
& &\e \Bigl[ \sup_{r\in [0,t ]} |{N}_{n}(r\wedge \taun (R))| \Bigr] 
 \le C \, p \cdot \e \Bigl[ 
\Bigl( \sup_{ s\in  [0, t] } {|\un (s\wedge {\tau }_{n}(R) )|}^{p} 
{\Bigr) }^{\frac{1}{2}}
\nonumber \\
& & \cdot \Bigl( \int_{0}^{t\wedge \taun (R)} \,  { |\un (s )|}^{p-2} 
  \cdot
 \bigl[ {\lambda }_{} \,  {|\un (s )|}^{2} + \kappa   
+ (2- a ) \norm{\un (s )}{}{2} \bigr]
   \, ds
  {\Bigr) }^{\frac{1}{2}} \Bigr]  \nonumber  \\
& &\le \frac{1}{4} \e \bigl[ \sup_{ r \in[0, t] } {|\un (r\wedge {\tau }_{n}(R))|}^{p} \bigr] 
  \nonumber \\
& &+ {C}^{2}{p}^{2}
 \e \Bigl[ \int_{0}^{t\wedge \taun (R)} \bigl[ {\lambda }_{}  { |\un (s )|}^{p} 
 + \kappa { |\un (s )|}^{p-2}   + (2- a) { |\un (s )|}^{p-2} \norm{\un (s )}{}{2} \bigr]
   \, ds
 \Bigr]  \nonumber \\
& &\le \frac{1}{4} \e \bigl[ \sup_{ r\in [0, t] } {|\un (r\wedge {\tau }_{n}(R) )|}^{p} \bigr]
 + \tilde{\tilde{K}} , 
    \label{E:BDG_Nn_R}
\end{eqnarray}
where $\tilde{\tilde{K}} = {C}^{2}{p}^{2} [\lambda {\tilde{C}}_{p} + \kappa {\tilde{C}}_{p-2} + (2-a){C}_{2}]$.
(The constants ${\tilde{C}}_{p}, {\tilde{C}}_{p-2}$ are the same as in (\ref{E:E(int_un(t))_est}) and ${C}_{2}$ is the same as in (\ref{E:HV_estimate}).)
Therefore by (\ref{E:Ito_formula}) for all $t\in [0,T]$
\begin{eqnarray} 
& &\bigl| {u}_{n}  ({t\wedge \taun (R)}) {\bigr| }^{p}
 \le  c
+ {c}_{1}\int_{0}^{T} \bigl| \un (s){\bigr| }^{p}  \, ds 
 + \sup_{r \in [0,T] } |{M}_{n}(r\wedge \taun (R))|
   \nonumber  \\
& &  + |{I}_{n}(T\wedge \taun (R)) | + \sup_{r \in [0,T] } |{N}_{n}(r\wedge \taun (R))|,
 \label{E:Ito_BDG_R}
\end{eqnarray} 
where ${I}_{n}$ is defined by (\ref{E:I_n}). Since inequality (\ref{E:Ito_BDG_R}) holds for all $t\in [0,T]$ and the right-hand side of (\ref{E:Ito_BDG_R}) in independent of $t$, we infer that
\begin{eqnarray} 
& &\e \Bigl[ \sup_{t \in [0,T]}\bigl| {u}_{n}  (t\wedge \taun (R)) {\bigr| }^{p} \Bigr]
 \le  c  +
 + {c}_{1}\e \Bigl[ \int_{0}^{T} \bigl| \un (s){\bigr| }^{p}  \, ds  \Bigr] \nonumber \\
& & +  \e \Bigl[ \sup_{r \in [0,T] } |{M}_{n}(r\wedge \taun (R))| \Bigr] 
 +  \e \Bigl[ |{I}_{n}(T\wedge \taun (R)) |  \Bigr]
 +  \e \Bigl[ \sup_{r \in [0,T] } |{N}_{n}(r\wedge \taun (R))|  \Bigr] . \nonumber \\
    \label{E:Ito_BDG_R_E}
\end{eqnarray}
Using inequalities (\ref{E:E(int_un(t))_est}), (\ref{E:BDG_R_cont.}), (\ref{E:E_I_n(t)}) and (\ref{E:BDG_Nn_R}) in (\ref{E:Ito_BDG_R_E})  we infer that
$$ \label{E:H_est_R}
 \e \Bigl[ \sup_{t \in T}\bigl| {u}_{n}  (t\wedge \taun (R)) {\bigr| }^{p} \Bigr]
 \le {C}_{1}(p)
$$
for some constant ${C}_{1}(p)$ independent of $n \in \nat $ and $R>0$. Passing to the limit as $R \to \infty $, we obtain inequality  (\ref{E:H_estimate}). 
Thus the Lemma holds for $p\in \{ 2,4+\gamma  \} $.

\bigskip  \noindent
Let now $p \in [1,4+\gamma ) \setminus \{ 2\} $. Let us fix $n \in \nat $. Then
$$
   |{u}_{n}(t){|}_{H}^{p} =  \bigl( |{u}_{n}(t){|}_{H}^{4+\gamma } {\bigr) }^{\frac{p}{4+\gamma }}
   \le \bigl( \sup_{t\in [0,T]}|{u}_{n}(t){|}_{H}^{4+\gamma } {\bigr) }^{\frac{p}{4+\gamma }} ,
   \qquad t \in [0,T].
$$
Thus 
$$
  \sup_{t\in [0,T]} |{u}_{n}(t){|}_{H}^{p}
   \le \bigl( \sup_{t\in [0,T]}|{u}_{n}(t){|}_{H}^{4+\gamma } {\bigr) }^{\frac{p}{4+\gamma }}
$$   
and by the H\H{o}lder inequality 
\begin{eqnarray*}
& \e \Bigl[ \sup_{t\in [0,T]} |{u}_{n}(t){|}_{H}^{p} \Bigr] 
   &\le \e \Bigl[ \bigl( \sup_{t\in [0,T]}|{u}_{n}(t){|}_{H}^{4+\gamma } {\bigr) }^{\frac{p}{4+\gamma }}
    \Bigr]  \\
&  &  \le \Bigl( \e \Bigl[ \sup_{t\in [0,T]}|{u}_{n}(t){|}_{H}^{4+\gamma }  \Bigr] 
   {\Bigr) }^{\frac{p}{4+\gamma }} 
\le \bigl[ {C}_{1} (4+\gamma ) {\bigr] }^{\frac{p}{4+\gamma }} .
\end{eqnarray*}
Since $n \in \nat $ was chosen in an arbitray way, we infer that
$$
\sup_{n \in \nat } \e \Bigl[ \sup_{t\in [0,T]} |{u}_{n}(t){|}_{H}^{p} \Bigr] 
   \le {C}_{1}(p) ,
$$
where ${C}_{1}(p)=\bigl[ {C}_{1} (4+\gamma ) {\bigr] }^{\frac{p}{4+\gamma }}$.
The proof of Lemma is thus complete. \qed

\subsection{Tightness}

Let  $m > \frac{d}{2}+1$ be  fixed and let $U$ be the space defined by (\ref{E:U_comp_V_m}).
We will apply Corollary \ref{C:tigthness_criterion_cadlag_unbound} with $q:=2$.
So, let us consider the space 
\begin{equation} \label{E:Z_cadlag_NS} 
  \zcal  : = {L}_{w}^{2}(0,T;V) \cap {L}^{2}(0,T;{H}_{loc}) \cap \dmath ([0,T];U')
  \cap \dmath ([0,T];{H}_{w}).
\end{equation}
For each $n \in \nat $, the solution ${u}_{n} $ of the Galerkin equation defines a measure
$\lcal ({u}_{n})$ on $(\zcal , \tcal )$. 
Using Corollary \ref{C:tigthness_criterion_cadlag_unbound} we will prove that the set of measures 
$\bigl\{ \lcal ({u}_{n} ) , n \in \nat  \bigr\} $ is tight on $(\zcal , \tcal )$. The inequalities 
(\ref{E:H_estimate}) and (\ref{E:V_estimate}) in Lemma \ref{L:Galerkin_estimates } are of crucial importance. However, to prove tightness it is sufficient to use inequality (\ref{E:H_estimate}) only with 
$p=2$.

\begin{lemma} \label{L:comp_Galerkin}
The set of measures $\bigl\{ \lcal ({u}_{n} ) , n \in \nat  \bigr\} $ is tight on $(\zcal , \tcal )$.
\end{lemma}

\proof
We will apply Corollary \ref{C:tigthness_criterion_cadlag_unbound}.
By estimates  (\ref{E:H_estimate}) and (\ref{E:V_estimate}), conditions (a), (b) are satisfied.
Thus, it is sufficient to  prove that the sequence $(\un {)}_{n \in \nat }$ satisfies the Aldous condition \bf [A] \rm in the space $U'$. We will use Lemma \ref{L:Aldous_criterion} in Appendix A.
Let ${(\taun )}_{n \in \nat} $ be a sequence of stopping times such that $0 \le \taun \le T$.
By (\ref{E:Galerkin}), we have
\begin{eqnarray*}
& {u}_{n} (t) \,
 =  & \Pn {u}_{0}  - \int_{0}^{t} \Pn \acal  {u}_{n} (s) \, ds
  - \int_{0}^{t} \Bn \bigl( {u}_{n} (s) \bigr) \, ds 
  + \int_{0}^{t} \Pn f(s) \, ds  \\
& & + \int_{0}^{t} \int_{Y}\Pn F (s,{u}_{n}({s}^{-}),y) \tilde{\eta}(ds,dy) 
  + \int_{0}^{t} \Pn G(s,\un (s)) \, dW(s) \nonumber \\
& & =:    \Jn{1} + \Jn{2}(t) + \Jn{3}(t) + \Jn{4}(t) + \Jn{5}(t) + \Jn{6}(t), \qquad t \in [0,T].
\end{eqnarray*}
Let $\theta  >0 $. 
We will check that each term $\Jn{i}$, i=1,...,6, satisfies condition (\ref{E:Aldous_est})
in Lemma \ref{L:Aldous_criterion}.

\bigskip \noindent
Since $\acal :V \to V'$ and ${|\acal (u)|}_{V'} \le \norm{u}{}{}$
and the embedding  $V' \hookrightarrow U'$ is continuous,  by the H\H{o}lder inequality and (\ref{E:V_estimate}), we have the following estimates
\begin{eqnarray*}
& &\e \bigl[ \bigl| \Jn{2} (\taun + \theta ) - \Jn{2}(\taun ) {\bigr| }_{U'}  \bigr]
 = \e \Bigl[ {\Bigl| \int_{\taun }^{\taun + \theta } \Pn \acal {u}_{n} (s) \, ds \Bigr| }_{U'} \Bigr]
 \nonumber \\
& & \le c \, \e \Bigl[  \int_{\taun }^{\taun + \theta }
{\bigl|  \acal  {u}_{n} (s) \bigr| }_{V'}  \, ds \Bigr]
 \le c \, \e \Bigl[  \int_{\taun }^{\taun + \theta } \nonumber
 \norm{ {u}_{n} (s) }{}{}  \, ds \Bigr] \\
& & \le c \, \e \Bigl[ {\theta }^{\frac{1}{2}} 
 \Bigl( \int_{0 }^{T } \norm{{u}_{n} (s) }{}{2}\, ds {\Bigr) }^{\frac{1}{2}} \Bigr]
\le c \sqrt{{C}_{2}} \cdot {\theta }^{\frac{1}{2}}=: {c}_{2} \cdot {\theta }^{\frac{1}{2}}.  \label{E:Jn2}
\end{eqnarray*}
Thus $\Jn{2}$ satifies condition (\ref{E:Aldous_est}) with $\alpha =1$ and $\beta =\frac{1}{2}$.

\bigskip \noindent
Let us consider the term ${\Jn{3}}$.
Since $m > \frac{d}{2} +1 $ and $U\hookrightarrow {V}_{m}$,  by (\ref{E:estimate_B_ext}) 
and (\ref{E:H_estimate}) we have the following inequalities
\begin{eqnarray*}
& &\e \bigl[ \bigl| \Jn{3} (\taun + \theta )  - \Jn{3}(\taun) {\bigr| }_{U'}  \bigr]
 = \e \Bigl[ { \Bigl| \int_{\taun }^{\taun + \theta }
  \Bn \bigl({u}_{n} (s) \bigr) \, ds \Bigr| }_{U'} \Bigr] \nonumber \\
& &\le c\e \Bigl[  \int_{\taun }^{\taun + \theta }
{ \bigl|  B\bigl( {u}_{n} (s)  \bigr)  \bigr| }_{{V}_{m}'} \, ds \Bigr] 
 \le c\e \Bigl[  \int_{\taun }^{\taun + \theta }
  \| B \| \cdot \bigl|  {u}_{n} (s) { \bigr| }_{H}^{2}   \, ds \Bigr] \nonumber \\
& & \le c\| B \|  \cdot  \e \bigl[ \sup_{s \in [0,T]} \bigl| {u}_{n} (s) { \bigr| }_{H}^{2}\bigr] \cdot \theta
 \le c\| B \| \,  {C}_{1}(2)  \cdot \theta =: {c}_{3} \cdot \theta , \label{E:Jn3}
\end{eqnarray*}
where $\| B \| $ stands for the norm of $B: H \times H \to {V}_{m}'$. This means that ${\Jn{3}}$ satisfies condition (\ref{E:Aldous_est})  with $\alpha =\beta  =1$.

\bigskip \noindent
Let us move to the term ${\Jn{4}}$. By the H\H{o}lder inequality, we have
\begin{eqnarray*}
& &\e \bigl[ \bigl| \Jn{4} (\taun + \theta ) - \Jn{4}(\taun) {\bigr| }_{U'}  \bigr]
\le c  \e \Bigl[ {\Bigl| \int_{\taun }^{\taun + \theta } \!
  \Pn f (s)  ds \Bigr| }_{V'} \Bigr]  \nonumber \\
& & \le c  \cdot {\theta }^{\frac{1}{2}} \cdot   \e \Bigl[  
\Bigl( \int_{0 }^{T } {\bigl|  f(s) \bigr| }_{V'}^{2}  \, ds {\Bigr) }^{\frac{1}{2}} \Bigr]
= c  \cdot {\theta }^{\frac{1}{2}} \cdot \norm{f}{{L}^{2}(0,T;V')}{} =: {c}_{4} \cdot {\theta }^{\frac{1}{2}}. \label{E:Jn4}
\end{eqnarray*}
Hence condition (\ref{E:Aldous_est}) holds with $\alpha =1$ and $\beta =\frac{1}{2}$.

\bigskip \noindent
Let us consider the term ${\Jn{5}}$.  Since  $ H \hookrightarrow U' $, 
by (\ref{E:isometry}),
condition (\ref{E:F_linear_growth}) with $p=2$ in Assumption (F.2) and by (\ref{E:H_estimate}), we obtain the following inequalities
\begin{eqnarray*}
& & \e \bigl[ \bigl| \Jn{5} (\taun + \theta ) - \Jn{5}(\taun) {\bigr| }_{U'}^{2}  \bigr] \nonumber 
  =  \e \Bigl[ \Bigl| \int_{\taun }^{\taun + \theta } \int_{Y} \Pn F(s, {u}_{n} (s);y ) \, \tilde{\eta}(ds,dy)  {\Bigr| }_{U'}^{2} 
\Bigr] \nonumber \\
& &  \le c\e \Bigl[ \Bigl| \int_{\taun }^{\taun + \theta } \int_{Y} \Pn F(s, {u}_{n} (s);y ) \, \tilde{\eta}(ds,dy)  {\Bigr| }_{H}^{2} 
\Bigr] \nonumber \\
& & = c\e \Bigl[  \int_{\taun }^{\taun + \theta } \int_{Y} {\bigl| \Pn F(s, {u}_{n} (s);y ) \bigr| }_{H}^{2} \, \nu (dy)ds   
\Bigr]  \le C  \e \Bigl[  \int_{\taun }^{\taun + \theta }
  ( 1 + | {u}_{n} (s) {|}_{H}^{2} ) ds  \Bigr]  \nonumber \\
&  & \le C \cdot  \theta \cdot \bigl( 1+
  \e \bigl[ \sup_{s \in [0,T]} \bigl| {u}_{n} (s) { \bigr| }_{H}^{2} \bigr] \bigr)
\le C \cdot (1+ {C}_{1}(2)) \cdot \theta =: {c}_{5} \cdot \theta . \label{E:Jn5}
\end{eqnarray*}
Thus ${\Jn{5}}$ satisfies condition (\ref{E:Aldous_est})  with $\alpha =2$ and $\beta  =1$. 

\bigskip \noindent
Let us consider the term ${\Jn{6}}$. 
By the It\^{o} isometry,
condition (\ref{E:G*}) in assumption (G.3), continuity of the embedding $V'\hookrightarrow U'$ and inequality (\ref{E:H_estimate}), we have
\begin{eqnarray*}
& &\e \bigl[ \bigl| \Jn{6} (\taun +\theta ) - \Jn{6}(\taun) {\bigr| }_{U '}^{2}  \bigr] 
  = \e \Bigl[ \Bigl| \int_{\taun }^{\taun +\theta } 
    \Pn G(s,\un (s))\, dW(s) {\Bigr| }_{U '}^{2} \Bigr]  \nonumber \\
& & \le c  \e \Bigl[ \int_{\taun }^{\taun +\theta }(1+|\un (s){|}_{H}^{2} ) ds  \Bigr] 
 \le c   \theta \bigl( 1+ \e \bigl[ \sup_{s \in [0,T]} \bigl| \un (s) { \bigr| }_{H}^{2}\bigr] \bigr) 
\le c (1+ {C}_{1}(2) ) \theta  .  \label{E:Jn6}
\end{eqnarray*}
Thus ${\Jn{6}}$ satisfies condition (\ref{E:Aldous_est})  with $\alpha =2 $ and $ \beta =1$.

\bigskip  \noindent
By Lemma \ref{L:Aldous_criterion} the sequence $({u}_{n}{)}_{n\in \nat }$ satisfies the Aldous condition in the space $U'$.
This completes the proof of Lemma. \qed

\bigskip  \noindent
We will now  move to the proof of the main Theorem of existence of a martingale solution.
 The main difficulties occur in the term containing the nonlinearity $B$ and in the noise terms $F$ and $G$. To deal with the nonlinear term, we need inequality (\ref{E:H_estimate}) 
 for some $p>4$. Moreover, we will see that the sequence $({\bar{u}}_{n})$ of approximate solutions is convergent in the Fr\'{e}chet space ${L}^{2}(0,T;{H}_{loc })$. 
So, we will use the  property of the mapping $B$ contained in Lemma \ref{L:B_conv_aux} below.
Analogous problems appear in the noise terms, where assumptions (F.3) and (G.3) will be needed in the case when the domain $\ocal $ is unbounded. 
For simplicity  we  assume that $\dim {Y}_{W}=1$, i.e. we consider one-dimensional cylindrical Wiener process $W(t)$, $t \in [0,T]$. 
Construction of a martingale solution is based on the Skorokhod Theorem for nonmetric spaces.
The method is closely related to the approach due to Brze\'{z}niak and Hausenblas  
\cite{Brzezniak_Hausenblas_2010}.

\subsection{Proof of Theorem \ref{T:existence} }

\noindent
By Lemma \ref{L:comp_Galerkin} the set of measures $\bigl\{ \lcal ({u}_{n} ) , n \in \nat  \bigr\} $ is tight on the space $(\zcal , \tcal )$. Let ${\eta }_{n}:= \eta $, $n \in \nat $. The set of measures $\bigl\{ \lcal ({\eta }_{n} ) , n \in \nat  \bigr\} $ is tight on the space 
${M}_{\bar{\nat }}([0,T]\times Y)$.
Let ${W}_{n}:= W $, $n \in \nat $. The set $\bigl\{ \lcal ({W}_{n}) , n \in \nat  \bigr\} $ is tight 
on the space $\ccal ([0,T];\rzecz )$ of continuous function from $[0,T]$ to $\rzecz $ with the standard supremum-norm .
 Thus the set $ \bigl\{ \lcal ({u}_{n}, {\eta }_{n},{W}_{n}) , n \in \nat  \bigr\}$ is tight on  $\zcal \times {M}_{\bar{\nat }}([0,T]\times Y) \times \ccal ([0,T];\rzecz )$.
By Corollary \ref{C:Skorokhod_J,B,H} and Remark \ref{R:separating_maps}, see Appendix B, there exists a subsequence $({n}_{k}{)}_{k\in \nat }$, a probability space 
$\bigl( \bar{\Omega }, \bar{\fcal },\bar{\p }  \bigr) $ and, on this space, 
$\zcal \times {M}_{\bar{\nat }}([0,T]\times Y) \times \ccal ([0,T];\rzecz )$-valued random variables $({u}_{*},{\eta }_{*},{W}_{*})$, 
$({\bar{u}}_{k}, {\bar{\eta }}_{k},{\bar{W}}_{k})$, $k \in \nat $ such that
\begin{itemize}
\item[(i) ] $\lcal \bigl( ({\bar{u}}_{k}, {\bar{\eta }}_{k},{\bar{W}}_{k}) \bigr) = \lcal \bigl( ({u}_{{n}_{k}}, {\eta }_{{n}_{k}},{W}_{{n}_{k}}) \bigr) $ for all $ k \in \nat $;
\item[(ii) ] ${(\bar{u}}_{k}, {\bar{\eta }}_{k},{\bar{W}}_{k}) \to ({u}_{*},{\eta }_{*},{W}_{*})$ in $\zcal \times {M}_{\bar{\nat }}([0,T]\times Y)  \times \ccal ([0,T];\rzecz )$ with probability $1$ on $\bigl( \bar{\Omega }, \bar{\fcal },\bar{\p }  \bigr) $ as $k \to \infty $;
\item[(iii) ] $ ({\bar{\eta }}_{k} (\bar{\omega }), {\bar{W}}_{k}(\bar{\omega }) )
 =  ({\eta }_{*}(\bar{\omega }),{W}_{*}(\bar{\omega }))$ for all $\bar{\omega } \in \bar{\Omega }$.
\end{itemize}
We will denote this sequences again by $\bigl( \! ({u}_{n}, {\eta }_{n},{W}_{n}) \! {\bigr) }_{\! n\in \nat }$ and 
$\bigl( \! ({\bar{u}}_{n}, {\bar{\eta }}_{n},{\bar{W}}_{n}) \! {\bigr) }_{\! n\in \nat } $.
Moreover, ${\bar{\eta }}_{n}$, $n \in \nat $, and ${\eta }_{*}$ are time homogeneous Poisson random measures on $(Y, \ycal )$ with intensity measure $\nu $ and ${\bar{W}}_{n}$, $n \in \nat $, and ${W}_{*}$ are cylindrical Wiener processes, see \cite[Section 9]{Brzezniak_Hausenblas_2010}.
Using the definition of the space $\zcal $, see (\ref{E:Z_cadlag_NS} ),
in particular, we have
\begin{equation} \label{E:Z_cadlag_NS_bar_un_conv} 
  {\bar{u}}_{n} \to {u}_{*} \quad \mbox{in} \quad  {L}_{w}^{2}(0,T;V) \cap {L}^{2}(0,T;{H}_{loc}) \cap \dmath ([0,T];U')
  \cap \dmath ([0,T];{H}_{w}) \quad \bar{\p } \mbox{-a.s.}
\end{equation}
Since the random variables ${\bar{u}}_{n}$ and ${u}_{n}$ are identically distributed, we have the following inequalities.
For every $p\in [1,4+\gamma ] $  
\begin{equation} \label{E:H_estimate_bar_u_n}
 \sup_{n \ge 1 } \bar{\e } \bigl( \sup_{0 \le s \le T } |\bun{s}{|}_{H}^{p} \bigr) \le {C}_{1}(p) .
\end{equation}
and
\begin{equation} \label{E:V_estimate_bar_u_n}
  \sup_{n \ge 1 }  \bar{\e } \bigl[ \int_{0}^{T} \norm{ \bun{s}}{V}{2} \, ds \bigr] \le {C}_{2}.
\end{equation}
Let us fix $v \in U $. Analogously to \cite{Brzezniak_Hausenblas_2010}, let us denote
\begin{eqnarray} 
& &  {\kcal }_{n}({\bar{u}}_{n}, {\bar{\eta }}_{n},{\bar{W}}_{n}, v) (t)
:=  \ilsk{{\bar{u}}_{n}(0)}{v}{H} \nonumber \\
& & + \int_{0}^{t} \dual{\Pn \acal {\bar{u}}_{n}(s)}{v}{}  ds   
+ \int_{0}^{t} \dual{\Bn ({\bar{u}}_{n}(s))}{v}{}  ds \nonumber \\
& & + \int_{0}^{t} \dual{\Pn f(s)}{v}{}\, ds 
 + \int_{0}^{t} \int_{Y} \ilsk{\Pn F(s,{\bar{u}}_{n}({s}^{-});y)}{v}{H} \,  {\tilde{\bar{\eta }}}_{n} (ds,dy)
  \nonumber \\
& & + \Dual{\int_{0}^{t}\Pn G(s,{\bar{u}}_{n}(s))\,  d{\bar{W}}_{n}(s)}{v}{} 
\label{E:K_n_bar_u_n}
\end{eqnarray}
and 
\begin{eqnarray} 
& &  \kcal ({u}_{*}, {\eta }_{*},{W}_{*}, v) (t)
:=  \ilsk{{u}_{*}(0)}{v}{H} 
   + \int_{0}^{t} \dual{ \acal {u}_{*}(s)}{v}{}  ds 
 + \int_{0}^{t} \dual{ B({u}_{*}(s))}{v}{}  ds \nonumber \\
& & + \int_{0}^{t} \dual{ f(s)}{v}{}\, ds 
 + \int_{0}^{t} \int_{Y} \ilsk{ F(s,{u}_{*}({s}^{-});y)}{v}{H} \,  \tilde{{\eta }}_{*} (ds,dy) \nonumber \\
& & +  \Dual{\int_{0}^{t}G(s,{u}_{*}(s))\,  d{W}_{*}(s)}{v}{}  ,
\quad t \in [0,T] .
\label{E:K_u*}
\end{eqnarray}
\bf Step ${1}^{0}$. \rm  We will prove that 
\begin{equation} \label{E:bar_u_n_convergence}
 \lim_{n \to \infty } \norm{\ilsk{{\bar{u}}_{n}(\cdot )}{v}{H} 
 - \ilsk{{u}_{*}(\cdot )}{v}{H} }{{L}^{2}([0,T]\times \bar{\Omega })}{} =0 
\end{equation}
and
\begin{equation} \label{E:K_n_bar_u_n_convergence}
     \lim_{n \to \infty } \norm{  {\kcal }_{n}({\bar{u}}_{n}, {\bar{\eta }}_{n},{\bar{W}}_{n},v)
 - \kcal ({u}_{*}, {\eta }_{*},{W}_{*}, v) }{{L}^{2}([0,T]\times \bar{\Omega })}{} =0 .
\end{equation}

\bigskip  \noindent
To prove (\ref{E:bar_u_n_convergence}) let us write
\begin{eqnarray*}
& & \norm{\ilsk{{\bar{u}}_{n}(\cdot )}{v}{H} 
 - \ilsk{{u}_{*}(\cdot )}{v}{H} }{{L}^{2}([0,T]\times \bar{\Omega })}{2} \\
& &  = \int_{\bar{\Omega }}\int_{0}^{T} \bigl| \ilsk{{\bar{u}}_{n}(t) -{u}_{*}(t)}{v}{H} {\bigr| }^{2} \, dt \bar{\p }(d\omega ) 
=  \bar{\e } \Bigl[ \int_{0}^{T} \bigl| \ilsk{{\bar{u}}_{n}(t) -{u}_{*}(t)}{v}{H} {\bigr| }^{2} \, dt \Bigr] .
\end{eqnarray*}
Moreover,  
\begin{eqnarray*}
& &  \int_{0}^{T} \bigl| \ilsk{{\bar{u}}_{n}(t) -{u}_{*}(t)}{v}{H} {\bigr| }^{2} \, dt
= \int_{0}^{T} \bigl| \dual{{\bar{u}}_{n}(t) -{u}_{*}(t)}{v}{U',U} {\bigr| }^{2} \, dt \\
& & \le \norm{v}{U}{2}\int_{0}^{T} | {\bar{u}}_{n}(t) -{u}_{*}(t){|}_{U'}^{2} \, dt.
\end{eqnarray*}
Since by (\ref{E:Z_cadlag_NS_bar_un_conv})
$ {\bar{u}}_{n} \to {u}_{*} $ in $\dmath ([0,T];U')$ and by (\ref{E:H_estimate_bar_u_n}) 
 $\sup_{t\in [0,T]}| {\bar{u}}_{n}(t){|}_{H}^{2} < \infty $, $\bar{\p }$-a.s. and the embedding $H \hookrightarrow U'$ is continuous, by the Dominated Convergence Theorem we infer that $\bar{\p }$-a.s.,
$ {\bar{u}}_{n} \to {u}_{*} $ in $ {L}^{2}(0,T;U') $. Then 
\begin{equation} \label{E:Vitali_{u}_{n}_conv}
  \lim_{n \to \infty } \int_{0}^{T} \bigl| \ilsk{{\bar{u}}_{n}(t) -{u}_{*}(t)}{v}{H} {\bigr| }^{2} \, dt
  =0 .
\end{equation} 
Moreover,  by the H\H{o}lder inequality and (\ref{E:H_estimate_bar_u_n})
for every $n \in \nat$ and every $r\in \bigl( 1, 2+ \frac{\gamma }{2}\bigr] $
\begin{eqnarray}
& & \bar{\e }\Bigl[ \Bigl|  \int_{0}^{T} \bigl| {\bar{u}}_{n}(t) -{u}_{*}(t) {\bigr| }_{H}^{2} \, dt 
 {\Bigr| }^{r}\Bigr] 
  \le c \bar{\e }\Bigl[  \int_{0}^{T} \bigl( \bigl| {\bar{u}}_{n}(t)  {\bigr| }_{H}^{2r}
+ \bigl| {u}_{*}(t) {\bigr| }_{H}^{2r} \bigr) \, dt \Bigr]  \nonumber \\
& & \le \tilde{c} \bar{\e }\bigl[  \sup_{t\in [0,T]}  \bigl| {\bar{u}}_{n}(t)  {\bigr| }_{H}^{2r}  \bigr]  \le \tilde{c} {C}_{1}(2r)    \label{E:Vitali_{u}_{n}_est}
\end{eqnarray} 
for some constants $c, \tilde{c}>0$.
By (\ref{E:Vitali_{u}_{n}_conv}), (\ref{E:Vitali_{u}_{n}_est}) and the Vitali Theorem we infer that
\begin{eqnarray*}
\lim_{n\to \infty }
\bar{\e } \Bigl[ \int_{0}^{T} \bigl| \ilsk{{\bar{u}}_{n}(t) -{u}_{*}(t)}{v}{H} {\bigr| }^{2} \, dt \Bigr]
=0,
\end{eqnarray*}
i.e. (\ref{E:bar_u_n_convergence}) holds.

\bigskip \noindent 
Let us move to the proof of (\ref{E:K_n_bar_u_n_convergence}).
Note that by the Fubini Theorem, we have
\begin{eqnarray*}
& & \norm{  {\kcal }_{n}({\bar{u}}_{n}, {\bar{\eta }}_{n},{\bar{W}}_{n},v)
 - \kcal ({u}_{*}, {\eta }_{*},{W}_{*}, v) }{{L}^{2}([0,T]\times \bar{\Omega })}{2} \\
& &= \int_{0}^{T} \int_{\bar{\Omega }} \bigl| {\kcal }_{n}({\bar{u}}_{n}, {\bar{\eta }}_{n},{\bar{W}}_{n},v)(t)
 - \kcal ({u}_{*}, {\eta }_{*},{W}_{*}, v)(t) {\bigr| }^{2}\,  d\bar{\p } (\omega )dt \\
& &= \int_{0}^{T} \bar{\e } \bigl[ \bigl| {\kcal }_{n}({\bar{u}}_{n}, {\bar{\eta }}_{n},{\bar{W}}_{n},v)(t)
 - \kcal ({u}_{*}, {\eta }_{*},{W}_{*} ,v)(t) {\bigr| }^{2}\, \bigr] dt .
 \end{eqnarray*}
We will prove that each term on the right hand side of  
(\ref{E:K_n_bar_u_n}) tends in ${L}^{2}([0,T]$ $\times \bar{\Omega })$ to the corresponding term in (\ref{E:K_u*}).

\bigskip  \noindent
Since by (\ref{E:Z_cadlag_NS_bar_un_conv}) ${\bar{u}}_{n} \to {u}_{*}$ in $\dmath (0,T;{H}_{w})$ $\bar{\p }$-a.s. and ${u}_{*}$ is continuous at $t=0$, we infer that
$ \ilsk{{\bar{u}}_{n}(0)}{v}{H} \to \ilsk{{u}_{*}(0)}{v}{H} $ $\bar{\p }$-a.s.
By (\ref{E:H_estimate_bar_u_n}) and the Vitali Theorem, we have
$$
  \lim_{n\to \infty } \bar{\e } \Bigl[ \bigl| \ilsk{{\bar{u}}_{n}(0)-{u}_{*}(0)}{v}{H}  {\bigr| }^{2} \Bigr] =0.
$$
Hence 
\begin{eqnarray} \label{E:u_n_0_conv}
 \lim_{n \to \infty } \norm{\ilsk{{\bar{u}}_{n}(0)-{u}_{*}(0)}{v}{H} }{{L}^{2}([0,T]\times \bar{\Omega })}{2}   =0 .
\end{eqnarray}

\bigskip  \noindent
By (\ref{E:Z_cadlag_NS_bar_un_conv}) ${\bar{u}}_{n} \to {u}_{*}$ in ${L}_{w}^{2}(0,T;V)$, $\bar{\p }$-a.s.
Moreover, since $v \in U$,  $\Pn v \to v $ in $V$, see Section \ref{S:Some_operators}.
Thus by  relation (\ref{E:op_Acal}) we infer that
\begin{eqnarray} 
& &\lim_{n\to \infty } \int_{0}^{t} \dual{\Pn \acal {\bar{u}}_{n}(s)}{v}{} \, ds 
 = \lim_{n \to \infty } \int_{0}^{t} \dirilsk{{\bar{u}}_{n}(s)}{\Pn v}{} \, ds \nonumber \\
& &= \int_{0}^{t} \dirilsk{ {u}_{*}(s)}{v}{} \, ds 
 = \int_{0}^{t} \dual{ {u}_{*}(s)}{v}{} \, ds . 
 \label{E:Vitali_{A}_{n}_conv}
\end{eqnarray}
By (\ref{E:op_Acal}), Lemma \ref{L:A_acal_rel}, the H\H{o}lder inequality and (\ref{E:H_estimate_bar_u_n})
 for all $t \in [0,T]$,  $r \in (0,2+\gamma ] $ and  $n \in \nat $
\begin{eqnarray} 
& &  \bar{\e }\Bigl[ \Bigl | \int_{0}^{t}  \dual{ \Pn \acal {\bar{u}}_{n}(s)}{v}{} \, ds {\Bigr| }^{2+r} \Bigr]
  = \bar{\e } \Bigl[ \Bigl | \int_{0}^{t}  \ilsk{  {\bar{u}}_{n}(s)}{(A-I) \Pn v}{H}  \, ds {\Bigr| }^{2+r} \Bigr] \nonumber \\
& & \le c\| { v}{\| }_{U}^{2+r}  \bar{\e } \Bigl[  \int_{0}^{T} \bigl| {\bar{u}}_{n}(s){\bigr| }_{H}^{2+r} 
    \, ds  \Bigr] 
  \le \tilde{c} \bar{\e } \bigl[  \sup_{0 \le s \le T } |\bun{s}{|}_{H}^{2+r}  \bigr]
  \le \tilde{c}{C}_{1}(2+r) 
  \label{E:Vitali_{A}_{n}_est}
\end{eqnarray}
for some constants $c, \tilde{c}>0$. 
Therefore by (\ref{E:Vitali_{A}_{n}_conv}), (\ref{E:Vitali_{A}_{n}_est}) and the Vitali Theorem we infer that for all $t \in [0,T]$
$$
  \lim_{n\to \infty } \bar{\e } \Bigl[ \Bigl| \int_{0}^{t} 
\dual{\Pn \acal {\bar{u}}_{n}(s)-\acal {u}_{*}(s)}{v}{} \, ds 
 {\Bigr| }^{2} \Bigr] =0.
$$
Hence by (\ref{E:H_estimate_bar_u_n}) and the Dominated Convergence Theorem
\begin{eqnarray} \label{E:{A}_{n}_conv}
  \lim_{n\to \infty } \int_{0}^{T}\bar{\e } \Bigl[ \Bigl| \int_{0}^{t} 
\dual{\Pn \acal {\bar{u}}_{n}(s)-\acal {u}_{*}(s)}{v}{} \, ds 
 {\Bigr| }^{2} \Bigr] \, dt =0.
\end{eqnarray}
Let us move to the nonlinear term. We will use the following auxilliary result proven in 
\cite{Brzezniak_Motyl_NS}. (We recall the proof in Appendix D.)

\bigskip  \noindent
\begin{lemma}  \label{L:B_conv_aux}
\rm (Lemma B.1 in \cite{Brzezniak_Motyl_NS}) \it 
Let $u \in {L}^{2}(0,T;H)$ and let  ${(\un )}_{n} $ be a bounded sequence in ${L}^{2}(0,T;H)$ such that $\un \to u $ in ${L}^{2}(0,T;{H}_{loc})$. Let $m > \frac{d}{2}+1$.
Then for all $t \in [0,T]$ and all $\psi \in {V}_{m}$:
$$
  \nlim \int_{0}^{t}  \dual{ B ( \un (s ))   }{\psi }{} \, ds
  = \int_{0}^{t} \dual{  B ( u (s ))   }{\psi }{} \, ds .  
$$
(Here $\dual{\cdot }{\cdot }{}$ denotes the dual pairing between the space ${V}_{m}$ and ${V}_{m}'$.) \qed
\end{lemma}

\bigskip  \noindent
Let us fix $m > \frac{d}{2}+1$.
Since by (\ref{E:V_estimate_bar_u_n}) and (\ref{E:norm_V}) the sequence $({\bar{u}}_{n})$ is bounded in ${L}^{2}(0,T;H)$ and 
by (\ref{E:Z_cadlag_NS_bar_un_conv}) ${\bar{u}}_{n} \to {u}_{*}$ in 
${L}^{2}(0,T;{H}_{loc})$ $\bar{\p }$-a.s.,
by Lemma \ref{L:B_conv_aux} we infer that $\bar{\p }$-a.s.
for all $t \in [0,T]$ and all $v \in {V}_{m}$
$$
 \lim_{n\to \infty } \int_{0}^{t} \dual{ B({\bar{u}}_{n}(s)) - B({u}_{*}(s))}{v}{} \, ds =0 .
$$
It is easy to see that for sufficiently large $n \in \nat $ 
$$
     \Bn ({\bar{u}}_{n}(s)) = \Pn B({\bar{u}}_{n}(s)) , \qquad s \in [0,T].
$$
Moreover, if $v \in U$ then $\Pn v \to v $ in ${V}_{m}$, see Section \ref{S:Some_operators}. Since $U \subset {V}_{m}$, we infer that for all $v \in U$ and all $t \in [0,T]$
\begin{equation} \label{E:Vitali_{B}_{n}_conv}
 \lim_{n\to \infty } \int_{0}^{t} \dual{ \Bn ({\bar{u}}_{n}(s)) - B({u}_{*}(s))}{v}{} \, ds =0 
 \qquad \bar{\p }\mbox{a.s}.
\end{equation}
By the H\H{o}lder inequality, (\ref{E:estimate_B_ext}) and  (\ref{E:H_estimate_bar_u_n}) we obtain for all 
 $t \in [0,T]$,  $r\in \bigl(0,\frac{\gamma }{2}\bigr] $ and  $n \in \nat $
\begin{eqnarray} 
& & \bar{\e } \Bigl[ \Bigl| \int_{0}^{t} \dual{\Bn ({\bar{u}}_{n}(s)) }{v}{} \, ds {\Bigr| }^{2+r}  \Bigr] 
    \le \bar{\e } \Bigl[ {t}^{1+r}{\| v\| }_{{V}_{m}}^{2+r} \int_{0}^{t} \bigl| \Bn ({\bar{u}}_{n}(s)) {\bigr| }_{{V}_{m}'}^{2+r}  \, ds   \Bigr] \nonumber \\
&  & \le C \bar{\e } \Bigl[ \! \int_{0}^{t} \! \bigl| {\bar{u}}_{n}(s) {\bigr| }_{H}^{2(2+r)} ds \Bigr]  
  \le \tilde{C} \bar{\e } \bigl[ \! \sup_{s\in [0,T]} \bigl| {\bar{u}}_{n}(s) {\bigr| }_{H}^{2(2+r)}   \bigr] 
 \le \tilde{C} {C}_{1}(4+2r) .  \label{E:Vitali_{B}_{n}_est}
\end{eqnarray}
In view of (\ref{E:Vitali_{B}_{n}_conv}) and (\ref{E:Vitali_{B}_{n}_est}), by the Vitali Theorem we 
obtain for all $t\in [0,T]$
\begin{equation} \label{E:Lebesque_{B}_{n}_conv}
  \lim_{n \to \infty } \bar{\e } \Bigl[ \Bigl| \int_{0}^{t} \dual{\Bn ({\bar{u}}_{n}(s)) - B({u}_{*}(s))}{v}{} \, ds {\Bigr| }^{2}  \Bigr] =0 .
\end{equation}
Since by (\ref{E:H_estimate_bar_u_n}) for all $t \in [0,T]$ and all $n \in \nat $
$$
   \bar{\e } \Bigl[ \Bigl| \int_{0}^{t} \dual{\Bn ({\bar{u}}_{n}(s)) }{v}{} \, ds {\Bigr| }^{2} \Bigr] \le c \bar{\e } \bigl[ \sup_{s\in [0,T]} \bigl| {\bar{u}}_{n}(s)) {\bigr| }_{H}^{4}   \bigr]
 \le c {C}_{1}(4)
$$
for some $c>0$,
by (\ref{E:Lebesque_{B}_{n}_conv}) and the Dominated Convergence Theorem, we infer that 
\begin{eqnarray} \label{E:{B}_{n}_conv}
  \lim_{n \to \infty } \int_{0}^{T} \bar{\e } \Bigl[ \Bigl| \int_{0}^{t} \dual{\Bn ({\bar{u}}_{n}(s)) - B({u}_{*}(s))}{v}{} \, ds {\Bigr| }^{2}  \Bigr] \, dt  =0 .
\end{eqnarray}

\bigskip  \noindent
Let us move to the noise terms. Let us assume first that $v \in \vcal $. For all $t \in [0,T]$ we have
\begin{eqnarray*}
& & \int_{0}^{t} \int_{Y} \bigl| \dual{F(s, {\bar{u}}_{n}({s}^{-});y)- F(s,{u}_{*}({s}^{-});y)}{v}{} {\bigr| }^{2} \, d\nu (y)ds  \nonumber \\
& & \int_{0}^{t} \int_{Y} \bigl| \ilsk{F(s, {\bar{u}}_{n}({s}^{-});y)- F(s,{u}_{*}({s}^{-});y)}{v}{H} {\bigr| }^{2} \, d\nu (y)ds  \nonumber \\
& &= \int_{0}^{t} \int_{Y} \bigl|  {\tilde{F}}_{v}({\bar{u}}_{n})(s,y) - {\tilde{F}}_{v}({u}_{*})(s,y)
    {\bigr| }^{2} \, d\nu (y)ds  \nonumber \\
& &\le \int_{0}^{T} \int_{Y} \bigl|  {\tilde{F}}_{v}({\bar{u}}_{n})(s,y) - {\tilde{F}}_{v}({u}_{*})(s,y)
    {\bigr| }^{2} \, d\nu (y)ds  \nonumber \\
& &= \norm{{\tilde{F}}_{v}({\bar{u}}_{n}) - {\tilde{F}}_{v}({u}_{*})}{{L}^{2}([0,T]\times Y;\rzecz )}{2}, 
\end{eqnarray*}
where ${\tilde{F}}_{v}$ is the mapping defined by (\ref{E:F**}).  
Since by (\ref{E:Z_cadlag_NS_bar_un_conv})
${\bar{u}}_{n} \to {u}_{*}$ in ${L}^{2}(0,T;{H}_{loc})$, $\bar{\p }$-a.s.,
by  assumption (F.3) we infer that for all $t\in [0,T]$ 
\begin{eqnarray}
& & \lim_{n \to \infty } \int_{0}^{t} \int_{Y} \bigl| \ilsk{F(s, {\bar{u}}_{n}({s}^{-});y)- F(s,{u}_{*}({s}^{-});y)}{v}{H} {\bigr| }^{2} \, d\nu (y)ds  =0. 
\label{E:Vitali_{noise}_{n}_conv}
\end{eqnarray}
Moreover, by inequality (\ref{E:F_linear_growth}) in assumption (F.2) and by (\ref{E:H_estimate_bar_u_n})
for every $t \in [0,T]$ every $r \in \bigl( 1,2+ \frac{\gamma }{2} \bigr] $ and every $n \in \nat $ the following inequalities hold
\begin{eqnarray} 
& & \bar{\e } \Bigl[ \Bigl| \int_{0}^{t} \int_{Y} \bigl| \ilsk{F(s, {\bar{u}}_{n}({s}^{-});y)- F(s,{u}_{*}({s}^{-});y)}{v}{H} {\bigr| }^{2} \, d\nu (y)ds {\Bigr| }^{r} \Bigr] \nonumber \\
& & \le {2}^{r} {|v|}_{H}^{2r}\bar{\e } \Bigl[ \Bigl| \int_{0}^{t} \int_{Y} \bigl\{ 
 \bigl| F(s, {\bar{u}}_{n}({s}^{-});y) {\bigr| }_{H}^{2}
 +\bigl| F(s,{u}_{*}({s}^{-});y) {\bigr| }_{H}^{2} \bigr\} \, d\nu (y)ds {\Bigr| }^{r} \Bigr] \nonumber \\
& &\le {2}^{r} {C}_{2}^{r} {|v|}_{H}^{2r}\, \bar{\e } \Bigl[ \Bigl|  \int_{0}^{t} \bigl\{ 2+
\bigl|  {\bar{u}}_{n}(s) {\bigr| }_{H}^{2}
 + \bigl| {u}_{*}(s) {\bigr| }_{H}^{2}
\bigr\}  \, ds {\Bigr| }^{r} \Bigr] 
 \le c \bigl( 1+ \bar{\e } \bigl[ \sup_{s \in [0,T]} \bigl| {\bar{u}}_{n}(s)  {\bigr| }_{H}^{2r}\bigr] \bigr)
 \nonumber \\
& &\le c (1+ {C}_{1}(2r))  
\label{E:Vitali_{noise}_{n}_est}
\end{eqnarray} 
for some constant $c>0$.
Thus by (\ref{E:Vitali_{noise}_{n}_conv}), (\ref{E:Vitali_{noise}_{n}_est}) and the Vitali Theorem for all $t\in [0,T]$
\begin{eqnarray}
\lim_{n\to \infty } \bar{\e } \Bigl[  \int_{0}^{t} \int_{Y} \bigl| \dual{F(s, {\bar{u}}_{n}({s}^{-});y)- F(s,{u}_{*}({s}^{-});y)}{v}{} {\bigr| }^{2} \, d\nu (y)ds  \Bigr] =0,
\quad  v \in \vcal . \label{E:Lebesgue_{noise}_{n}_conv_vcal}
\end{eqnarray} 
Let now $v \in H$ and let $\eps >0$. Since $\vcal $ is dense in $H$, there exists ${v}_{\eps } \in \vcal $ such that $|v-{v}_{\eps }{|}_{H}^{2}< \eps $. By (\ref{E:F_linear_growth}) the following inequalities hold
\begin{eqnarray*}
& & \int_{0}^{t} \int_{Y} \bigl| \dual{F(s, {\bar{u}}_{n}({s}^{-});y)- F(s,{u}_{*}({s}^{-});y)}{v}{} {\bigr| }^{2} \, d\nu (y)ds \\
& & \le 2 \int_{0}^{t} \int_{Y} \bigl| \dual{F(s, {\bar{u}}_{n}({s}^{-});y)- F(s,{u}_{*}({s}^{-});y)}{v-{v}_{\eps }}{} {\bigr| }^{2} \, d\nu (y)ds \\
& & \quad + 2 \int_{0}^{t} \int_{Y} \bigl| \dual{F(s, {\bar{u}}_{n}({s}^{-});y)- F(s,{u}_{*}({s}^{-});y)}{{v}_{\eps }}{} {\bigr| }^{2} \, d\nu (y)ds \\
& & \le 4 {C}_{2} {\eps }^{2} \int_{0}^{t} \bigl\{ 2+ \bigl| {\bar{u}}_{n}(s) {\bigr| }_{H}^{2}
 + \bigl| {u}_{*}(s) {\bigr| }_{H}^{2}  \bigr\} \,ds \\
& & \quad + 2 \int_{0}^{t} \int_{Y} \bigl| \dual{F(s, {\bar{u}}_{n}({s}^{-});y)- F(s,{u}_{*}({s}^{-});y)}{{v}_{\eps }}{} {\bigr| }^{2} \, d\nu (y)ds.
\end{eqnarray*} 
Hence by (\ref{E:H_estimate_bar_u_n})
\begin{eqnarray*}
& & \bar{\e } \Bigl[\int_{0}^{t} \int_{Y} \bigl| \dual{F(s, {\bar{u}}_{n}({s}^{-});y)- F(s,{u}_{*}({s}^{-});y)}{v}{} {\bigr| }^{2} \, d\nu (y)ds \Bigr] \\
& & \le 4 {C}_{2} {\eps }^{2} \bar{\e } \Bigl[ \int_{0}^{t} \bigl\{ 2+ \bigl| {\bar{u}}_{n}(s) {\bigr| }_{H}^{2}
 + \bigl| {u}_{*}(s) {\bigr| }_{H}^{2}  \bigr\} \,ds \Bigr] \\
& & \quad + 2 \bar{\e } \Bigl[ \int_{0}^{t} \int_{Y} \bigl| \dual{F(s, {\bar{u}}_{n}({s}^{-});y)- F(s,{u}_{*}({s}^{-});y)}{{v}_{\eps }}{} {\bigr| }^{2} \, d\nu (y)ds \bigr] \\
& & \le \tilde{c}  {\eps }^{2} 
+ 2 \bar{\e } \Bigl[ \int_{0}^{t} \int_{Y} \bigl| \dual{F(s, {\bar{u}}_{n}({s}^{-});y)- F(s,{u}_{*}({s}^{-});y)}{{v}_{\eps }}{} {\bigr| }^{2} \, d\nu (y)ds \bigr] .
\end{eqnarray*} 
Passing to the upper limit as $n \to \infty $ in the above inequality, by 
(\ref{E:Lebesgue_{noise}_{n}_conv_vcal}) we obtain
\begin{eqnarray*}
& & \limsup_{n \to \infty } \bar{\e } \Bigl[\int_{0}^{t} \int_{Y} \bigl| \dual{F(s, {\bar{u}}_{n}({s}^{-});y)- F(s,{u}_{*}({s}^{-});y)}{v}{} {\bigr| }^{2} \, d\nu (y)ds \Bigr] \le \tilde{c} {\eps }^{2}. 
\end{eqnarray*}
Since $\eps > 0$ was chosen in an arbitrary way, we infer that  for all $v \in H $
\begin{eqnarray*}
& & \lim_{n \to \infty } \bar{\e } \Bigl[\int_{0}^{t} \int_{Y} \bigl| \dual{F(s, {\bar{u}}_{n}({s}^{-});y)- F(s,{u}_{*}({s}^{-});y)}{v}{} {\bigr| }^{2} \, d\nu (y)ds \Bigr] =0 . 
\end{eqnarray*}
Moreover, since the restriction of $\Pn $ to the space $H$ is the $\ilsk{\cdot }{\cdot }{H}$-projection onto ${H}_{n}$, see Section \ref{S:Some_operators}, we infer that also
\begin{eqnarray*}
\lim_{n\to \infty } \bar{\e } \Bigl[  \int_{0}^{t} \int_{Y} \bigl| \dual{\Pn F(s, {\bar{u}}_{n}({s}^{-});y)- F(s,{u}_{*}({s}^{-});y)}{v}{} {\bigr| }^{2} \, d\nu (y)ds  \Bigr] =0,
\quad  v \in H . 
\end{eqnarray*}
Hence by the properties of the integral with respect to the compensated Poisson random measure
and the fact that ${\bar{\eta }}_{n} = {\eta }_{*}$, we have
\begin{equation}  \label{E:Lebesque_{noise}_{n}_conv}
\lim_{n\to \infty } \bar{\e } \Bigl[ \Bigl|  \int_{0}^{t} \int_{Y}  \dual{ \Pn F(s, {\bar{u}}_{n}({s}^{-});y)- F(s,{u}_{*}({s}^{-});y)}{v}{}  \, {\tilde{{\eta }}}_{*} (ds,dy)
{\Bigr| }^{2}  \Bigr]  =0.
\end{equation}
Moreover, by inequality (\ref{E:F_linear_growth}) in assumption (F.2) and by (\ref{E:H_estimate_bar_u_n}) we obtain the following inequalities
\begin{eqnarray} 
& & \bar{\e } \Bigl[ \Bigl|  \int_{0}^{t} \int_{Y}  \dual{\Pn F(s, {\bar{u}}_{n}({s}^{-});y)- F(s,{u}_{*}({s}^{-});y)}{v}{}  \, {\tilde{\eta }}_{*} (ds,dy)
{\Bigr| }^{2}  \Bigr] \nonumber \\
& & =  \bar{\e } \Bigl[   \int_{0}^{t} \int_{Y} \bigl| \ilsk{\Pn F(s, {\bar{u}}_{n}({s}^{-});y)- F(s,{u}_{*}({s}^{-});y)}{v}{H} {\bigr| }^{2} \, \nu (dy)ds  \Bigr] \nonumber  \\
& & \le 2 {|v|}_{H}^{2} \bar{\e } \Bigl[  \int_{0}^{t} \int_{Y}\bigl\{
    \bigl| \Pn F(s,{\bar{u}}_{n}({s}^{-});y)  {\bigr| }_{H}^{2}
  + \bigl| F(s,{u}_{*}({s}^{-});y) {\bigr| }_{H}^{2} \bigr\} \, \nu (dy)ds \Bigr] \nonumber \\
& & \le 2{C}_{2}  {|v|}_{H}^{2} \bar{\e } \Bigl[  \int_{0}^{t} \bigl\{ 2+
    \bigl| {\bar{u}}_{n}(s)  {\bigr| }_{H}^{2} +
   \bigl| {u}_{*}(s) {\bigr| }_{H}^{2} \bigr\} \, ds \Bigr] 
\le c\bigl( 1+ \bar{\e } \bigl[ \sup_{s \in [0,T]} \bigl| {\bar{u}}_{n}(s){\bigr| }_{H}^{2}\bigr] \bigr) 
 \nonumber  \\
& & \le c ( 1+ {C}_{1}(2)) .
 \label{E:Lebesque_{noise}_{n}_est}
\end{eqnarray}
By (\ref{E:Lebesque_{noise}_{n}_conv}), (\ref{E:Lebesque_{noise}_{n}_est}) and the Dominated Convergence Theorem, we have for all $v \in H $
\begin{eqnarray}
 \lim_{n\to \infty } \int_{0}^{T}\bar{\e } \Bigl[ \Bigl|  \int_{0}^{t} \int_{Y}  \dual{\Pn F(s, {\bar{u}}_{n}({s}^{-});y)- F(s,{u}_{*}({s}^{-});y)}{v}{}  \, {\tilde{\eta }}_{*} (ds,dy) 
{\Bigr| }^{2}  \Bigr] \, dt  =0.  \label{E:{noise}_{n}_conv}
\end{eqnarray} 
Since $U \subset H$, (\ref{E:{noise}_{n}_conv}) holds for all $v \in U $, as well.

\bigskip  \noindent
Let us move to the second part of the noise. Let us assume first that $v \in \vcal $.
We have
\begin{eqnarray*} 
& & \int_{0}^{t}  \bigl\| \dual{G(s, {\bar{u}}_{n}(s))- G(s,{u}_{*}(s))}{v}{} 
  {\bigr\| }^{2}_{\lhs ({Y}_{w};\rzecz )} \, ds \\
& & =\int_{0}^{t}  \bigl\| \ilsk{G(s, {\bar{u}}_{n}(s))- G(s,{u}_{*}(s))}{v}{H} 
  {\bigr\| }^{2}_{\lhs ({Y}_{w};\rzecz )} \, ds  \nonumber  \\
& &= \int_{0}^{t}  \bigl\|  {\tilde{G}}_{v}({\bar{u}}_{n})(s) - {\tilde{G}}_{v}({u}_{*})(s)
    {\bigr\| }^{2}_{\lhs ({Y}_{w};\rzecz )} \, ds \\
& & \le \int_{0}^{T}  \bigl\|  {\tilde{G}}_{v}({\bar{u}}_{n})(s) - {\tilde{G}}_{v}({u}_{*})(s)
    {\bigr\| }^{2}_{\lhs ({Y}_{w};\rzecz )} \, ds \nonumber \\ 
& & = \norm{{\tilde{G}}_{v}({\bar{u}}_{n}) - {\tilde{G}}_{v}({u}_{*})}{{L}^{2}([0,T];\lhs ({Y}_{w};\rzecz )) }{2} ,
 \end{eqnarray*}
where   ${\tilde{G}}_{v}$ is the mapping defined by (\ref{E:G**}).
Since by (\ref{E:Z_cadlag_NS_bar_un_conv})
${\bar{u}}_{n} \to {u}_{*}$ in ${L}^{2}(0,\! T;\! {H}_{loc}\! )$, $\bar{\p }$-a.s.,
by the second part of assumption (G.3) we infer that for all $t\in [0,T]$ and all $v \in \vcal $
\begin{eqnarray} 
& & \lim_{n \to \infty } \int_{0}^{t}  \bigl\| \dual{G(s, {\bar{u}}_{n}(s))- G(s,{u}_{*}(s))}{v}{} 
  {\bigr\| }^{2}_{\lhs ({Y}_{w};\rzecz )} \, ds =0  .
   \label{E:Vitali_{Gaussian}_{n}_conv}
\end{eqnarray}
Moreover, by (\ref{E:G*}) and (\ref{E:H_estimate_bar_u_n}) we see that
for every $t \in [0,T]$ every $r \in \bigl( 1,2+ \frac{\gamma }{2} \bigr] $ and every $n \in \nat $
\begin{eqnarray} 
& &\bar{\e } \Bigl[  \Bigl| \int_{0}^{t}  \bigl\| \dual{G(s, {\bar{u}}_{n}(s))- G(s,{u}_{*}(s))}{v}{} 
  {\bigr\| }^{2}_{\lhs ({Y}_{w};\rzecz )} \, ds {\Bigr| }^{r} \Bigr]  \nonumber \\
& & \le  c \, \bar{\e } \Bigl[ \norm{v}{V}{2r} \cdot \int_{0}^{t} \bigl\{ 
  \norm{G(s, {\bar{u}}_{n}(s))}{\lhs ({Y}_{W};V')}{2r} 
  + \norm{G(s, {u}_{*}(s))}{\lhs ({Y}_{W};V')}{2r} \bigr\} \, ds   \Bigr] \nonumber \\
& & \le  {c}_{1} \, \bar{\e } \Bigl[ \int_{0}^{T} (1+ |{\bar{u}}_{n}(s){|}_{H}^{2r}
    + |{u}_{*}(s){|}_{H}^{2r}) \, ds  \Bigr] \nonumber  \\
&  & \le  \tilde{c} \Bigl\{ 1+ \bar{\e } \Bigl[ \sup_{s\in [0,T]}  |{\bar{u}}_{n}(s){|}_{H}^{2r}
    + \sup_{s\in [0,T]} |{u}_{*}(s){|}_{H}^{2r})  \Bigr] \Bigr\} \le \tilde{c} (1+2 {C}_{1}(2r))
 \label{E:Vitali_{Gaussian}_{n}_est}   
\end{eqnarray}
for some positive constants $c,{c}_{1},\tilde{c}$.
Thus by (\ref{E:Vitali_{Gaussian}_{n}_conv}), (\ref{E:Vitali_{Gaussian}_{n}_est}) and the Vitali Theorem 
\begin{eqnarray}
\lim_{n\to \infty } \bar{\e } \Bigl[  \int_{0}^{t}  
 \bigl\| \dual{G(s, {\bar{u}}_{n}(s))- G(s,{u}_{*}(s))}{v}{} {\bigr\| }^{2}_{\lhs ({Y}_{w};\rzecz )} \, ds
  \Bigr] =0  \quad \mbox{for all } v \in \vcal . \label{E:Lebesgue_{Gaussian}_{n}_conv_vcal}
\end{eqnarray}
Let now $v \in V$ and let $\eps >0$. Since $\vcal $ is dense in $V$, there exists ${v}_{\eps }\in \vcal  $ such that $\norm{v-{v}_{\eps }}{V}{} \le \eps $. We have the following inequalities
\begin{eqnarray*}
& & \int_{0}^{t}  
 \bigl\| \dual{G(s, {\bar{u}}_{n}(s))- G(s,{u}_{*}(s))}{v}{} {\bigr\| }^{2}_{\lhs ({Y}_{w};\rzecz )} \, ds \\ 
& & \le 2 \int_{0}^{t}  
 \bigl\| \dual{G(s, {\bar{u}}_{n}(s))- G(s,{u}_{*}(s))}{v-{v}_{\eps }}{} {\bigr\| }^{2}_{\lhs ({Y}_{w};\rzecz )} \, ds \\
& & \quad + 2 \int_{0}^{t}  
 \bigl\| \dual{G(s, {\bar{u}}_{n}(s))- G(s,{u}_{*}(s))}{{v}_{\eps }}{} {\bigr\| }^{2}_{\lhs ({Y}_{w};\rzecz )} \, ds .
\end{eqnarray*}
Moreover, by inequality  (\ref{E:G*}) in assumption (G.3), we obtain the following estimates
\begin{eqnarray*}
& & \int_{0}^{t}  
 \bigl\| \dual{G(s, {\bar{u}}_{n}(s))- G(s,{u}_{*}(s))}{v-{v}_{\eps }}{} {\bigr\| }^{2}_{\lhs ({Y}_{w};\rzecz )} \, ds \\
& &  \le  \norm{v-{v}_{\eps }}{V}{2} \int_{0}^{t}  
 \bigl\| G(s, {\bar{u}}_{n}(s))- G(s,{u}_{*}(s))  {\bigr\| }^{2}_{\lhs ({Y}_{w};V' )} \, ds \\
& &  \le c \bigl( 1+ \sup_{s\in [0,T]} {|{\bar{u}}_{n}(s)|}_{H}^{2} + \sup_{s\in [0,T]} {|{u}_{*}(s)|}_{H}^{2} \bigr) 
 {\eps }^{2} 
\end{eqnarray*}
for some $c>0$.
Thus by (\ref{E:H_estimate_bar_u_n}) we obtain the following inequalities
\begin{eqnarray*}
& & \bar{\e } \Bigl[ \int_{0}^{t}  
 \bigl\| \dual{G(s, {\bar{u}}_{n}(s))- G(s,{u}_{*}(s))}{v}{} {\bigr\| }^{2}_{\lhs ({Y}_{w};\rzecz )} \, ds
 \Bigr]  \\
& &  \le 2c \bigl\{ 1+ 2{C}_{1}(2) \bigr\} {\eps }^{2}  
  + 2 \bar{\e } \Bigl[\int_{0}^{t}  
 \bigl\| \dual{G(s, {\bar{u}}_{n}(s))- G(s,{u}_{*}(s))}{{v}_{\eps }}{} {\bigr\| }^{2}_{\lhs ({Y}_{w};\rzecz )} \, ds \Bigr] .
\end{eqnarray*} 
Passing to the upper limit as $n \to \infty $ by (\ref{E:Lebesgue_{Gaussian}_{n}_conv_vcal}) 
we infer that for all $v \in V$
\begin{eqnarray*}
& & \limsup_{n \to \infty } \bar{\e } \Bigl[ \int_{0}^{t}  
 \bigl\| \dual{G(s, {\bar{u}}_{n}(s))- G(s,{u}_{*}(s))}{v}{} {\bigr\| }^{2}_{\lhs ({Y}_{w};\rzecz )} \, ds
 \Bigr]   
\le C {\eps }^{2} , 
\end{eqnarray*} 
where $C=2c \bigl\{ 1+ 2{C}_{1}(2) \bigr\} $.
Since $\eps >0$ was chosen in an arbitrary way, we infer that
\begin{eqnarray}
\lim_{n\to \infty } \bar{\e } \Bigl[  \int_{0}^{t}  
 \bigl\| \dual{G(s, {\bar{u}}_{n}(s))- G(s,{u}_{*}(s))}{v}{} {\bigr\| }^{2}_{\lhs ({Y}_{w};\rzecz )} \, ds
  \Bigr] =0 \quad \mbox{for all } v \in V.  \label{E:Lebesgue_{Gaussian}_{n}_conv_V}
\end{eqnarray}  
For every $v \in V $ and every $s \in [0,T]$  we have
\begin{eqnarray*}
& &  \dual{\Pn G(s, {\bar{u}}_{n}(s))- G(s,{u}_{*}(s))}{v}{}
 = \dual{ G(s, {\bar{u}}_{n}(s))}{\Pn v}{}
 - \dual{ G(s,{u}_{*}(s))}{v}{} \\
& & = \dual{ G(s, {\bar{u}}_{n}(s))}{\Pn v - v}{}
 + \dual{ G(s, {\bar{u}}_{n}(s)) -G(s,{u}_{*}(s))}{v}{} \\
& & \le \norm{G(s, {\bar{u}}_{n}(s))}{\lhs ({Y}_{W},V')}{}  \norm{\Pn v -v}{V}{}
 + \dual{ G(s, {\bar{u}}_{n}(s)) -G(s,{u}_{*}(s))}{v}{} . 
\end{eqnarray*}
Thus by inequality  (\ref{E:G*}) in assumption (G.3) and by (\ref{E:H_estimate_bar_u_n}) we obtain
\begin{eqnarray*}
& &  \bar{\e } \Bigl[  \int_{0}^{t}  
 \bigl\| \dual{\Pn G(s, {\bar{u}}_{n}(s))- G(s,{u}_{*}(s))}{v}{} {\bigr\| }^{2}_{\lhs ({Y}_{w};\rzecz )} \, ds  \Bigr] \\
& &  \le 2 \norm{\Pn v -v}{V}{2} \, \bar{\e } \Bigl[  \int_{0}^{T} \norm{G(s, {\bar{u}}_{n}(s))}{\lhs ({Y}_{W},V')}{2} \, ds \Bigr] \\
& & \quad + 2 \bar{\e } \Bigl[  \int_{0}^{t}  
 \bigl\| \dual{ G(s, {\bar{u}}_{n}(s))- G(s,{u}_{*}(s))}{v}{} {\bigr\| }^{2}_{\lhs ({Y}_{w};\rzecz )} \, ds  \Bigr] \\
 & &  \le 2C \norm{\Pn v -v}{V}{2} \, \bar{\e } \Bigl[  \int_{0}^{T} \bigl( 1+ {|{\bar{u}}_{n}(s)|}_{H}^{2} \bigr) \, ds \Bigr] \\
& & \quad + 2 \bar{\e } \Bigl[  \int_{0}^{t}  
 \bigl\| \dual{ G(s, {\bar{u}}_{n}(s))- G(s,{u}_{*}(s))}{v}{} {\bigr\| }^{2}_{\lhs ({Y}_{w};\rzecz )} \, ds  \Bigr]  \\
 & & \le 2CT(1+ {C}_{1}(2)) \norm{\Pn v -v}{V}{2} \\
& & \quad  + 2 \bar{\e } \Bigl[  \int_{0}^{t}  
 \bigl\| \dual{ G(s, {\bar{u}}_{n}(s))- G(s,{u}_{*}(s))}{v}{} {\bigr\| }^{2}_{\lhs ({Y}_{w};\rzecz )} \, ds  \Bigr] .
\end{eqnarray*}
Since $U \subset V $ and $\norm{\Pn v -v}{V}{} \to 0 $ for all $v\in U$, see Section \ref{S:Some_operators}, by 
(\ref{E:Lebesgue_{Gaussian}_{n}_conv_V}) we infer that
\begin{eqnarray*}
 \lim_{n\to \infty } \bar{\e } \Bigl[ \int_{0}^{t}  
 \bigl\| \dual{\Pn G(s, {\bar{u}}_{n}(s))- G(s,{u}_{*}(s))}{v}{} {\bigr\| }^{2}_{\lhs ({Y}_{w};\rzecz )} \, ds  \Bigr]  =0 \quad \mbox{for all } v \in U .
\end{eqnarray*} 
Hence by the properties of the It\^{o} integral we infer that for all $t \in [0,T]$ and all $v \in U$
\begin{equation}  \label{E:Lebesque_{Gaussian}_{n}_conv}
\lim_{n\to \infty } \bar{\e } \Bigl[ \Bigl| 
 \Dual{\int_{0}^{t}\bigl[ \Pn G(s, {\bar{u}}_{n}(s))- G(s,{u}_{*}(s)) \bigr] \, d{W}_{*}(s) }{v}{}
 {\Bigr| }^{2} \Bigr] =0 .
\end{equation}
Moreover, by the It\^{o} isometry, inequality  (\ref{E:G*}) in assumption (G.3), and (\ref{E:H_estimate_bar_u_n}) we have for all $t\in [0,T]$
and all $n \in \nat $ 
\begin{eqnarray} 
& & \bar{\e } \Bigl[ \Bigl| 
 \Dual{\int_{0}^{t}\bigl[ \Pn G(s, {\bar{u}}_{n}(s))- G(s,{u}_{*}(s)) \bigr] \, d{W}_{*}(s) }{v}{}
 {\Bigr| }^{2} \Bigr] \nonumber \\
& & = \bar{\e } \Bigl[\int_{0}^{t}  \bigl\| \dual{\Pn G(s, {\bar{u}}_{n}(s))- G(s,{u}_{*}(s))}{v}{} 
  {\bigr\| }^{2}_{\lhs ({Y}_{w};\rzecz )} \, ds \Bigr] \nonumber \\
& & \le  c \Bigl\{ 1+ \bar{\e } \Bigl[ \sup_{s\in [0,T]}  |{\bar{u}}_{n}(s){|}_{H}^{2}
    + \sup_{s\in [0,T]} |{u}_{*}(s){|}_{H}^{2})  \Bigr] \Bigr\} \le c (1+2 {C}_{1}(2)) 
 \label{E:Lebesque_{Gaussian}_{n}_est}  
\end{eqnarray}
for some $c>0$.
By (\ref{E:Lebesque_{Gaussian}_{n}_conv}), (\ref{E:Lebesque_{Gaussian}_{n}_est}) and the Dominated Convergence Theorem we infer that
\begin{eqnarray}
 \lim_{n\to \infty } \int_{0}^{T}\bar{\e } \Bigl[ \Bigl| 
 \Dual{\int_{0}^{t}\bigl[ \Pn G(s, {\bar{u}}_{n}(s))- G(s,{u}_{*}(s)) \bigr] \, d{W}_{*}(s) }{v}{}
 {\Bigr| }^{2} \Bigr] =0.  \label{E:{Gaussian}_{n}_conv}
\end{eqnarray}
By (\ref{E:u_n_0_conv}), (\ref{E:{A}_{n}_conv}), (\ref{E:{B}_{n}_conv}), (\ref{E:{noise}_{n}_conv}) and  
(\ref{E:{Gaussian}_{n}_conv})
 the proof of (\ref{E:K_n_bar_u_n_convergence}) is complete.

\bigskip  \noindent
\bf Step ${2}^{0}$. \rm Since ${u}_{n}$ is a solution of the Galerkin equation, for all $t\in [0,T]$
$$
   \ilsk{{u}_{n}(t)}{v}{H} = {\kcal }_{n} ({u}_{n},{\eta }_{n},{W}_{n},v)(t) , \qquad \p \mbox{-a.s.}
$$
In particular,
$$
  \int_{0}^{T} \e \bigl[ \bigl|  \ilsk{\un (t)}{v}{H} 
 - {\kcal }_{n} ({u}_{n},{\eta }_{n},{W}_{n},v)(t){\bigr| }^{2} \, \bigr] \, dt  =0.
$$
Since $\lcal ({u}_{n},{\eta }_{n},{W}_{n}) = \lcal ({\bar{u}}_{n},{\bar{\eta }}_{n},{\bar{W}}_{n})$,
$$
  \int_{0}^{T} \bar{\e } \bigl[ \bigl|  \ilsk{\bun{t}}{v}{H} 
- {\kcal }_{n} ({\bar{u}}_{n},{\bar{\eta } }_{n},{\bar{W}}_{n},v)(t){\bigr| }^{2} \, \bigr] \, dt  =0.
$$
Moreover, by (\ref{E:bar_u_n_convergence}) and (\ref{E:K_n_bar_u_n_convergence})
$$
  \int_{0}^{T} \bar{\e } \bigl[ \bigl| 
 \ilsk{{u}_{*}(t)}{v}{H} - {\kcal }_{} ({u}_{*},{\eta }_{*},{W}_{*},v)(t){\bigr| }^{2} \, \bigr] \, dt  =0.
$$
Hence for $l$-almost all $t \in [0,T]$ and $\bar{\p }$-almost all $\omega \in \bar{\Omega }$
$$
\ilsk{{u}_{*}(t)}{v}{H} - {\kcal }_{} ({u}_{*},{\eta }_{*},{W}_{*},v)(t)=0, 
$$
i.e. for $l$-almost all $t \in [0,T]$ and $\bar{\p }$-almost all $\omega \in \bar{\Omega }$
\begin{eqnarray*}
& \ilsk{{u}_{*}(t)}{v}{H} - & \ilsk{{u}_{*}(0)}{v}{H}  +  \int_{0}^{t} \dual{\acal {u}_{*}(s)}{v}{} \, ds
+ \int_{0}^{t} \dual{B({u}_{*}(s),{u}_{*}(s))}{v}{} \, ds \\
& &-\int_{0}^{t} \dual{f(s)}{v}{}\, ds 
 - \int_{0}^{t} \int_{Y} \ilsk{F(s,{u}_{*}(s);y)}{v}{H} \,  {\tilde{\eta }}_{*} (ds,dy) \\
& & - \Dual{ \int_{0}^{t} G(s,{u}_{*}(s)) \, d {W}_{*}(s) }{v}{} =0.
\end{eqnarray*}
Since ${u}_{*}$ is $\zcal $-valued random variable, in particular ${u}_{*}\in \dmath ([0,T];{H}_{w})$, i.e. ${u}_{*}$ is weakly c\`{a}dl\`{a}g. Hence the function on the left-hand side of the above equality is c\`{a}dl\`{a}g with respect to $t$. Since two c\`{a}dl\`{a}g functions equal for $l$-almost all $t \in [0,T]$ must be equal for all $t\in [0,T]$, we infer that for all $t\in [0,T]$ and all $v \in U$
\begin{eqnarray*}
& \ilsk{{u}_{*}(t)}{v}{H} - & \ilsk{{u}_{*}(0)}{v}{H}  +  \int_{0}^{t} \dual{\acal {u}_{*}(s)}{v}{} \, ds
+ \int_{0}^{t} \dual{B({u}_{*}(s),{u}_{*}(s))}{v}{} \, ds \\
& & -\int_{0}^{t} \dual{f(s)}{v}{}\, ds 
 - \int_{0}^{t} \int_{Y} \ilsk{F(s,{u}_{*}(s);y)}{v}{H} \,  {\tilde{\eta }}_{*} (ds,dy) \\
& & - \Dual{ \int_{0}^{t} G(s,{u}_{*}(s)) \, d {W}_{*}(s) }{v}{} =0  \qquad 
\bar{\p } \mbox{-a.s.}
\end{eqnarray*}
Since $U$ is dense in V, we infer that the above equality holds for all $v \in V$.
Putting $\bar{u}:= {u}_{*}$, $\bar{\eta }:= {\eta }_{*}$ and $\bar{W}:= {W}_{*}$, we infer that the system 
$(\bar{\Omega }, \bar{\fcal },\bar{\p }, \bar{\fmath }, \bar{u}, \bar{\eta },\bar{W})$ is a martingale solution of the equation (\ref{E:NS}). The proof of Theorem \ref{T:existence} is thus complete.
\qed

\section{Appendix A}

\subsection{The Aldous condition } \label{S:Aldous_cadlag}
 \noindent
Here $(\smath ,\varrho )$ is a separable and complete metric space.
Let $(\Omega , \fcal ,\p )$ be a probability space with filtration $\mathbb{F}:=({\fcal }_{t}{)}_{t \in [0,T]}$ satisfying the usual hypotheses, see \cite{Metivier_82},
and let $(\Xn {)}_{n \in \nat }$ be a sequence of c\`{a}dl\`{a}g, $\mathbb{F}$-adapted and $\smath $-valued processes.

\begin{definition}  \rm (see \cite{Joffe_Metivier_86})
We say that the sequence $(\Xn )$ of $\smath $-valued random variables
 satifies  condition [$\mathbf{\tilde{T}}$] iff
\begin{itemize}
\item[\mbox{[$\mathbf{\tilde{T}}$]}] $ \, \, \, \,  \forall \, \eps >0 \quad \forall \, \eta >0 \quad \exists \, \delta  >0 $: 
$$ \label{E:cond_modulus_cadlag}
    \sup_{n \in \nat } \, \p \bigl\{  {w}_{[0,T]} (\Xn ,\delta )  > \eta \bigr\}  \le \eps .  
$$ 
\end{itemize} 
\noindent
Let us recall that ${w}_{[0,T]}$ stands for the modulus defined by (\ref{E:modulus_cadlag}).
\end{definition}

\noindent
\bf Remark. \rm Let ${\p }_{n}$ denote the law of $\Xn $ on $\dmath ([0,T], \smath )$.
For fixed $\eta > 0 $ and $\delta >0 $ we denote
$$
    {C}_{\eta ,\delta } := \{ u \in  \dmath ([0,T], \smath ): \, \, \, {w}_{[0,T]} (u ,\delta )  \ge \eta \}  .
$$ 
Then condition 
$$
    \p \bigl\{  {w}_{[0,T]} (\Xn ,\delta )  > \eta \bigr\}  \le \eps 
$$
is equivalent to 
$$
     {\p }_{n} ({C}_{\eta ,\delta }) \le \eps .
$$
\begin{lemma}  \label{L:modulus_cadlag_conv}
Assume that $(\Xn )$ satifies  condition $\mathbf{[\tilde{T}]}$.
Let ${\p }_{n}$ be the  law of $\Xn $ on 

\noindent
$\dmath ([0,T], \smath )$, $n \in \nat $.
Then for every $\eps >0 $ there exists a subset ${A}_{\eps } \subset \dmath ([0,T], \smath ) $
such that
$$
   \sup_{n\in \nat } {\p }_{n} ({A}_{\eps }) \ge 1 - \eps 
$$ 
and 
\begin{equation} \label{E:modulus_cadlag_conv}
   \lim_{\delta \to 0 }  \sup_{u \in {A}_{\eps } }  {w}_{[0,T]} (u,\delta ) =0.
\end{equation}
\end{lemma}
\proof
Fix $\eps >0$. By $\mathbf{[\tilde{T}]}$, for each $k \in \nat $ there exists ${\delta }_{k}>0$
such that
$$
   \sup_{n\in \nat } \p \bigl\{  {w}_{[0,T]} (\Xn ,{\delta }_{k})  > \frac{1}{k} \bigr\}  \le \frac{\eps }{{2}^{k+1}}.
$$
Then
$$
   \sup_{n\in \nat } \p \bigl\{  {w}_{[0,T]} (\Xn ,{\delta }_{k})  \le \frac{1}{k} \bigr\}  \ge 1-  \frac{\eps }{{2}^{k+1}}
$$
or equivalently
$$
     \sup_{n\in \nat }  {\p }_{n} \bigl\{ u \in \dmath ([0,T], \smath ) :  
  \, \, \, {w}_{[0,T]} ( u ,{\delta }_{k})  \le \frac{1}{k} \bigr\}  \ge 1-  \frac{\eps }{{2}^{k+1}}  
$$
Let ${B}_{k} :=  \bigl\{ u \in \dmath ([0,T], \smath ) :  
  \, \, \, {w}_{[0,T]} ( u ,{\delta }_{k})  \le \frac{1}{k} \bigr\}   $  and let 
${A}_{\eps } := \bigcap_{k=1}^{\infty } {B}_{k}$.
We assert that for each $n \in \nat $
$$
   {\p }_{n} \bigl(  {A}_{\eps }  \bigr) \ge 1- \eps . 
$$  
Indeed, we have the following estimate
\begin{eqnarray*}
& {\p }_{n} \bigl( \dmath ([0,T], \smath ) \! \setminus \! {A}_{\eps }  \bigr) 
 & \le  {\p }_{n} \bigl( \dmath ([0,T], \smath ) \! \setminus \! \bigcap_{k=1}^{\infty } {B}_{k}  \bigr) 
 = {\p }_{n} \bigl( \bigcup_{k=1}^{\infty } \bigl( \dmath ([0,T], \smath ) \setminus   {B}_{k} \bigr)  \bigr) \\
& & \le \sum_{k=1}^{\infty } {\p }_{n} \bigl( \dmath ([0,T], \smath ) \setminus   {B}_{k} \bigr) 
  \le \sum_{k=1}^{\infty } \frac{\eps }{{2}^{k+1}}  = \eps .
\end{eqnarray*}
Thus $ {\p }_{n} ({A}_{\eps }) \ge 1 - \eps $.

\bigskip  \noindent
To prove (\ref{E:modulus_cadlag_conv}), let us fix $\tilde{\eps } >0$. Directly from the definition of ${A}_{\eps }$,
we infer that $\sup_{u \in {A}_{\eps }} $ $ {w}_{[0,T]} ( u ,{\delta }_{k}) $ $ \le \frac{1}{k} $ for each $k \in \nat $.
Choose ${k}_{0} \in \nat $ such that $\frac{1}{{k}_{0}} \le \tilde{\eps } $ and let ${\delta }_{0} := {\delta }_{{k}_{0}}$.
Then for every $\delta \le  {\delta }_{0}$ we obtain
$$
    {w}_{[0,T]} (u, \delta ) \le {w}_{[0,T]} (u, {\delta }_{{k}_{0}}  )  \le \tilde{\eps } 
$$
which completes the proof of  (\ref{E:modulus_cadlag_conv}) and the proof of Lemma. \qed

\bigskip  \noindent
Now, we recall the Aldous condition which is connected with condition $\mathbf{[\tilde{T}]}$ 
(see \cite{Joffe_Metivier_86}, \cite{Metivier_88} and \cite{Aldous}). This condition allows to investigate the modulus for the sequence of stochastic processes by means of stopped processes.

\begin{definition}  \rm \label{D:Aldous}
A sequence $({X}_{n}{)}_{n\in \nat }$  satisfies  condition \rm \bf [A] \rm
iff
\begin{itemize}
\item[\mbox{[A]}] $ \, \, \, \,  \forall \, \eps >0 \quad \forall \, \eta >0 \quad \exists \, \delta >0 
$ such that for every sequence $({{\tau}_{n} } {)}_{n \in \nat }$ of $\mathbb{F}$-stopping times with
${\tau }_{n}\le T$ one has
$$
    \sup_{n \in \nat} \, \sup_{0 \le \theta \le \delta }  \p \bigl\{ 
   \varrho \bigl( {X}_{n} ({\tau }_{n} +\theta ),{X}_{n} ( {\tau }_{n}  ) \bigr) \ge \eta \bigr\}  \le \eps .  
$$
\end{itemize} 
\end{definition}

\begin{lemma} \label{L:Aldous_cadlag_equiv}
(See \cite{Joffe_Metivier_86}, Th. 2.2.2)
Condition \rm [\bf A\rm ] \it implies condition $\mathbf{[\tilde{T}]}$. 
\end{lemma}

\bigskip  \noindent
In the following Remark we formulate a certain condition which guaranties that the sequence $({X}_{n}{)}_{n\in \nat }$  satisfies  condition \rm \bf [A]\rm .

\bigskip  
\begin{lemma} \label{L:Aldous_criterion}  \it
Let $(E,\norm{\cdot }{E}{})$ be a separable Banach space and let $({X}_{n}{)}_{n \in \nat }$ be a sequence of $E$-valued random variables.
Assume that for every sequence $({{\tau}_{n} } {)}_{n \in \nat }$ of $\mathbb{F}$-stopping times with
${\tau }_{n}\le T$ and for every $n \in \nat $ and $\theta \ge 0 $ the following condition holds
\begin{equation} \label{E:Aldous_est}
 \e \bigl[ \bigl( \norm{ {X}_{n} ({\tau }_{n} +\theta )-{X}_{n} ( {\tau }_{n}  ) }{E}{\alpha } \bigr] \le C {\theta }^{\beta }
\end{equation}
for some $\alpha ,\beta >0$ and some constant $C >0 $. Then the sequence $({X}_{n}{)}_{n\in \nat }$  satisfies  condition \rm \bf [A] \rm in the space $E$.  \rm
\end{lemma}

\bigskip
\proof
Let us fix $\eps > 0$ and $\eta >0$. By the Chebyshev inequality for every $n \in \nat $ and every $\theta >0$ we have 
$$
 \p \bigl\{ 
   \norm{{X}_{n} ({\tau }_{n} +\theta )-{X}_{n} ( {\tau }_{n}  ) }{E}{} \ge \eta \bigr\} 
   \le \frac{1}{{\eta }^{\alpha }} \e \bigl[ \bigl( \norm{ {X}_{n} ({\tau }_{n} +\theta )-{X}_{n} ( {\tau }_{n}  )}{E}{\alpha } \bigr] 
 \le \frac{C {\theta }^{\beta }}{{\eta }^{\alpha }}. 
$$
Let $\delta := \bigl( \frac{{\eta }^{\alpha } \eps }{C} {\bigr) }^{\frac{1}{\beta }} $. Let us fix $n \in \nat $. Then for every $\theta \in [0,\delta ]$ we have the following inequalities
$$
 \p \bigl\{ 
   \norm{ {X}_{n} ({\tau }_{n} +\theta )-{X}_{n} ( {\tau }_{n}  ) }{E}{} \ge \eta \bigr\} 
   \le \frac{C {\theta }^{\beta }}{{\eta }^{\alpha }}  
 \le \frac{C}{{\eta }^{\alpha }} \Bigl[ 
 \Bigl( \frac{{\eta }^{\alpha } \eps }{C} {\Bigr) }^{\frac{1}{\beta }} {\Bigr] }^{\beta } = \eps .
$$
Hence
$$
    \sup_{0 \le \theta \le \delta }  \p \bigl\{ 
   \norm{ {X}_{n} ({\tau }_{n} +\theta )-{X}_{n} ( {\tau }_{n}  ) }{E}{} \ge \eta \bigr\}  \le \eps .  
$$ 
Since the above inequality holds for every $n \in \nat$, one has
$$
    \sup_{n \in \nat} \, \sup_{0 \le \theta \le \delta }  \p \bigl\{ 
   \norm{ {X}_{n} ({\tau }_{n} +\theta )-{X}_{n} ( {\tau }_{n}  )} {E}{} \ge \eta \bigr\}  \le \eps ,  
$$ 
i.e. condition \rm \bf [A] \rm  is satisfied. This completes the proof. \qed

\subsection{Proof of Corollary \ref{C:tigthness_criterion_cadlag_unbound}} 

 \noindent
Let $\eps >0$. By the  Chebyshev inequality and by (a), we infer that for any $r>0$
$$
  \p \Bigl( \sup_{s \in [0,T]} {|\Xn (s) |}_{H}{} > r  \Bigr) 
  \le \frac{ \e \bigl[ \sup_{s \in [0,T]} {|\Xn (s) |}_{H}\bigr] }{r}
  \le \frac{{C}_{1}}{r}.
$$
Let ${R}_{1}$ be such that $\frac{{C}_{1}}{{R}_{1}} \le \frac{\eps }{3}$. Then
$$
   \p \Bigl( \sup_{s \in [0,T]} {|\Xn (s)|}_{H} > {R}_{1}  \Bigr) \le \frac{\eps }{3}
$$
Let ${B}_{1} := \bigl\{ u \in {\zcal }_{q}: \, \, \sup_{s \in [0,T]} {|u(s)|}_{H}  \le {R}_{1} \bigr\} $.

\bigskip  \noindent
By the  Chebyshev inequality and by (b), we infer that for any $r>0$
$$
  \p \bigl( \norm{\Xn }{{L}^{q}(0,T;V)}{} > r  \bigr)   
  \le \frac{\e \bigl[ \norm{\Xn }{{L}^{q}(0,T;V)}{q} \bigr]  }{{r}^{q}}
  \le \frac{{C}_{2}}{{r}^{q}}.
$$
Let ${R}_{2}$ be such that $\frac{{C}_{2}}{{R}_{2}^{q}} \le \frac{\eps }{3}$. Then
$$
   \p \bigl( \norm{\Xn }{{L}^{q}(0,T;V)}{} > {R}_{2}  \bigr) \le \frac{\eps }{3}. 
$$
Let ${B}_{2} := \bigl\{ u \in {\zcal }_{q}: \, \, \norm{u }{{L}^{q}(0,T;V)}{} \le {R}_{2} \bigr\} $.

\bigskip  \noindent
By Lemmas \ref{L:Aldous_cadlag_equiv} and \ref{L:modulus_cadlag_conv} there exists a subset 
${A}_{\frac{\eps }{3}} \subset \dmath ([0,T], {E}_{1}) $ such that 
${\tilde{\p }}_{n} \bigl( {A}_{\frac{\eps }{3}}\bigr) \ge 1 - \frac{\eps }{3}$ and 
$$
    \lim_{\delta \to 0 }  \sup_{u \in {A}_{\frac{\eps }{3} } }  {w}_{[0,T]} (u,\delta ) =0.
$$
It is sufficient to define ${K}_{\eps } $ as the closure  of the set
${B}_{1} \cap {B}_{2} \cap {A}_{\frac{\eps }{3}}$ in ${\zcal }_{q}$.
By Theorem  \ref{T:Dubinsky_cadlag_unbound}, ${K}_{\eps }$ is compact in ${\zcal }_{q}$. The proof is thus complete. \qed

\section{Appendix B: The Skorokhod Embedding Theorems} \label{S:Skorokhod}
 \noindent
Let us recall the following Jakubowski's version of the Skorokhod Theorem \cite{Jakubowski_1998}, 
see also Brze\'{z}niak and Ondrej\'{a}t  \cite{Brzezniak_Ondrejat_2011}.

\begin{theorem} \label{T:2_Jakubowski} \rm (Theorem 2 in \cite{Jakubowski_1998}). \it
Let $(\xcal , \tau )$ be a topological space such that there exists a sequence $( {f}_{m} ) $ of continuous functions ${f}_{m}:\xcal  \to \rzecz $ that separates points of $\xcal $.
Let $({X}_{n})$ be a sequence of  $\xcal $ valued random variables. Suppose that for every $\eps >0$
there exists a compact subset ${K}_{\eps} \subset \xcal $ such that
$$ 
  \sup_{n \in \nat } \p (\{ {X}_{n} \in {K}_{\eps } \} ) > 1-\eps . 
$$  
Then there exists a subsequence $({X}_{{n}_{k}}{)}_{k\in \nat }$, a sequence $({Y}_{k}{)}_{k\in \nat }$ of $\xcal $ valued random variables and an $\xcal $ valued random variable $Y$ defined on some probability space  $(\Omega , \fcal ,\p )$ such that
$$ 
    \lcal ({X}_{{n}_{k}}) = \lcal ({Y}_{k}), \qquad k=1,2,... 
$$ 
and forall $ \omega \in \Omega $:
$$ 
  {Y}_{k}(\omega ) \stackrel{\tau }{\longrightarrow } Y(\omega )  \quad \mbox{ as } k \to \infty .
$$
\end{theorem}

\bigskip  \noindent
We will use the following version of the Skorokhod Theorem due to Brze\'{z}niak and Hausenblas 
\cite{Brzezniak_Hausenblas_2010}. 

\bigskip  \noindent
\begin{theorem} \rm (Theorem E.1 in \cite{Brzezniak_Hausenblas_2010}) \it \label{T:E.1_[B,H'2010]}
Let ${E}_{1}, {E}_{2}$ be two separable Banach spaces and let ${\pi }_{i}:  {E}_{1}\times {E}_{2} \to {E}_{i} $, $i=1,2$, be the projection onto ${E}_{i}$, i.e.
$$
  {E}_{1}\times {E}_{2} \ni \chi = ({\chi }_{1},{\chi }_{2}) \to {\pi }_{i} (\chi ) \in {E}_{i}.
$$ 
Let $(\Omega ,\fcal ,\p )$ be a probability space and let 
${\chi }_{n}:\Omega \to {E}_{1}\times {E}_{2}$, $n\in \nat $, be a family of random variables such that the sequence $\{ \lcal aw({\chi }_{n}), n \in \nat \} $ is weakly convergent on ${E}_{1}\times {E}_{2}$.
Finally let us assume that there exists a random variable $\rho :\Omega \to {E}_{1}$ such that 
$\lcal aw({\pi }_{1}\circ {\chi }_{n}) = \lcal aw(\rho )$, $\forall \, n \in \nat $.

\bigskip  \noindent
Then there exists a probability space $(\bar{\Omega }, \bar{\fcal }, \bar{\p })$, a family of ${E}_{1}\times {E}_{2}$-valued random variables $\{ {\bar{\chi }}_{n}, \, n \in \nat  \} $ on 
$(\bar{\Omega }, \bar{\fcal }, \bar{\p })$ and a random variable ${\chi }_{*}: \bar{\Omega } \to 
{E}_{1}\times {E}_{2}$ such that 
\begin{itemize}
\item[(i) ] $\lcal aw ({\bar{\chi }}_{n}) = \lcal aw ({\chi }_{n})$, $\, \forall \, n \in \nat $;
\item[(ii) ] ${\bar{\chi }}_{n} \to {\chi }_{*}$ in ${E}_{1}\times {E}_{2}$ a.s.;
\item[(iii) ] ${\pi }_{1} \circ {\bar{\chi }}_{n} (\bar{\omega }) = {\pi }_{1} \circ {\chi }_{*}(\bar{\omega })$ for all $\bar{\omega } \in \bar{\Omega }$.
\end{itemize}
\end{theorem}

\bigskip  \noindent
\bf Remark.  \rm Theorem \ref{T:E.1_[B,H'2010]} remains true if we substitute the Banach spaces ${E}_{1}, {E}_{2}$ by the separable complete metric spaces. 

\bigskip \noindent
Using the ideas due to Jakubowski \cite{Jakubowski_1998}, we can proof the following generalization of Theorem \ref{T:E.1_[B,H'2010]} to the case of nonmetric spaces. Let us notice that in comparison to 
Theorem \ref{T:E.1_[B,H'2010]} we will assume that the sequence $\{ \lcal aw({\chi }_{n}), n \in \nat \} $ is tight. The assumption of the weak convergence of $\{ \lcal aw({\chi }_{n}), n \in \nat \} $ is not sufficient in the case of nonmetric spaces, see  \cite{Jakubowski_1998}.

\bigskip 
\begin{corollary} \label{C:Skorokhod_J,B,H} \rm (Corollary 5.3 in \cite{Motyl_NS_Poisson_pre}) \it 
 Let ${\xcal }_{1}$ be a separable complete metric space and let ${\xcal }_{2}$ be a topological space such that
there exists a sequence  $\{ {f}_{\iota }{ \} }_{\iota \in \nat } $ of continuous functions ${f}_{\iota }:{\xcal }_{2} \to \rzecz $  separating points of ${\xcal }_{2}$.
Let $\xcal := {\xcal }_{1}\times {\xcal }_{2}$ with the Tykhonoff topology induced by the projections 
$$
  {\pi }_{i}: {\xcal }_{1}\times {\xcal }_{2} \to {\xcal }_{i} , \qquad  i=1,2.
$$
Let $(\Omega ,\fcal ,\p )$ be a probability space and let 
${\chi }_{n}:\Omega \to {\xcal }_{1}\times {\xcal }_{2}$, $n\in \nat $, be a family of random variables such that the sequence $\{ \lcal aw({\chi }_{n}), n \in \nat \} $ is tight on ${\xcal }_{1}\times {\xcal }_{2}$.
Finally let us assume that there exists a random variable $\rho :\Omega \to {\xcal }_{1}$ such that 
$\lcal aw({\pi }_{1}\circ {\chi }_{n}) = \lcal aw(\rho )$ for  all $ n \in \nat $.

\bigskip \noindent
Then there exists a subsequence $\bigl( {\chi }_{{n}_{k}} {\bigr) }_{k \in \nat } $,
a probability space $(\bar{\Omega }, \bar{\fcal }, \bar{\p })$, a family of ${\xcal }_{1}\times {\xcal }_{2}$-valued random variables $\{ {\bar{\chi }}_{k}, \, k \in \nat  \} $ on 
$(\bar{\Omega }, \bar{\fcal }, \bar{\p })$ and a random variable ${\chi }_{*}: \bar{\Omega } \to 
{\xcal }_{1}\times {\xcal }_{2}$ such that 
\begin{itemize}
\item[(i) ] $\lcal aw ({\bar{\chi }}_{k}) = \lcal aw ({\chi }_{{n}_{k}})$ for all $ k \in \nat $;
\item[(ii) ] ${\bar{\chi }}_{k} \to {\chi }_{*}$ in ${\xcal }_{1}\times {\xcal }_{2}$ a.s. as $k \to \infty $;
\item[(iii) ] ${\pi }_{1} \circ {\bar{\chi }}_{k} (\bar{\omega }) = {\pi }_{1} \circ {\chi }_{*}(\bar{\omega })$ for all $\bar{\omega } \in \bar{\Omega }$.
\end{itemize}
\end{corollary}

\noindent
For the convenience of the reader we recall the proof.

\proof
Using the ideas due to Jakubowski \cite{{Jakubowski_1998}}, the proof can be reduced to Theorem \ref{T:E.1_[B,H'2010]}.
Let us denote 
$$
    {\chi }_{n} = \bigl( {\chi }^{1}_{n} ,{\chi }^{2}_{n}  \bigr) , 
$$   
where ${\chi }^{i}_{n} := {\pi }_{i} \circ {\chi }_{n}$, $i=1,2$.
Since the sequence $\{ \lcal aw({\chi }_{n}), n \in \nat \} $ is tight on ${\xcal }_{1}\times {\xcal }_{2}$, we infer that the sequence $\{ \lcal aw({\chi }^{2}_{n}), n \in \nat \} $ is tight on ${\xcal }_{2}$.
Let ${K}_{m} \subset {\xcal }_{2}$ be compact subsets such that  ${K}_{m} \subset {K}_{m+1}$, $m=1,2,...$ and 
\begin{equation} \label{E:16_Jakubowski} 
 \sup_{n \in \nat } \p (\{ {\chi }^{2}_{n} \in {K}_{m} \} ) > 1 - \frac{1}{m}.
\end{equation}
\noindent
Let us consider the mapping $\tilde{f} : {\xcal }_{2} \to {\rzecz }^{\nat }$ defined by
$$
  \tilde{f}(z):=\bigl( {f}_{1}(z),{f}_{2}(z),...\bigr) =\bigl( {f}_{\iota }(z){\bigr) }_{\iota \in \nat },
 \qquad z \in {\xcal }_{2}.  
$$
$$
  {\tilde{\mu }}_{n} := \lcal (\tilde{f}({\chi }^{2}_{n}))  \qquad \mbox{and} \qquad 
  {\tilde{K}}_{m}:= \tilde{f} ({K}_{m}).
$$ 
On the set  ${\rzecz }^{\nat }$ let us define the function  
\begin{eqnarray*} 
  &\Phi (y) &:= 
     \min \{ m: \, \, \, y \in {\tilde{K}}_{m} \}  
 \quad \mbox{ if } y \in \bigcup_{m=1}^{\infty }{\tilde{K}}_{m} \nonumber \\
 &\Phi (y) & =  + \infty  \quad  \mbox{ otherwise.} 
\end{eqnarray*} 
Function $\Phi :{\rzecz }^{\nat } \to \nat $ is lower semicontinuous, i.e.: if ${y}_{n}\to {y}_{0}$ in ${\rzecz }^{\nat }$, then 
$$  
  \liminf_{n\to \infty } \Phi ({y}_{n}) \ge \Phi ({y}_{0})
$$
From (\ref{E:16_Jakubowski}) it follows that
\begin{itemize}
\item[$\bullet $] $\Phi < \infty $ (${\tilde{\mu }}_{n}$-p.p.) for all $  n \in \nat $
\item[$\bullet $] and $({\tilde{\mu }}_{n} \circ {\Phi }^{-1})$ is a tight sequence of laws on $\nat $. 
\end{itemize}
Furthermore, the sequence of laws  $\bigl\{ \lcal aw(\tilde{f}\circ {\chi }^{2}_{n}, \Phi \circ \tilde{f}\circ {\chi }^{2}_{n}) , \, n \in \nat \bigr\} $ is tight on ${\rzecz }^{\nat } \times \nat $.

\bigskip  \noindent
Let us consider the product space ${\xcal }_{1} \times ({\rzecz }^{\nat } \times \nat )$ and let 
${P}_{1}:= {\xcal }_{1} \times ({\rzecz }^{\nat } \times \nat ) \to {\xcal }_{1} $ be the projection onto ${\xcal }_{1}$ and ${P}_{2}:= {\xcal }_{1} \times ({\rzecz }^{\nat } \times \nat ) \to  {\rzecz }^{\nat } \times \nat $ be the projection onto ${\rzecz }^{\nat } \times \nat $. Moreover let ${\xi }_{n}$, $n \in \nat $, be ${\xcal }_{1} \times ({\rzecz }^{\nat } \times \nat )$-valued random variables defined by
$$
   {\xi }_{n} := \bigl( {\xi }^{1}_{n}, {\xi }^{2}_{n}\bigr) : \Omega \to {\xcal }_{1} \times ({\rzecz }^{\nat } \times \nat ) , 
$$ 
where
$$
   {\xi }^{1}_{n} := {\chi }^{1}_{n} \qquad \mbox{and} \qquad
   {\xi }^{2}_{n} := (\tilde{f}\circ {\chi }^{2}_{n}, \Phi \circ \tilde{f}\circ {\chi }^{2}_{n}) ,
   \qquad n \in \nat .
$$
Remark that the sequence of laws $\bigl\{ \lcal aw ({\xi }_{n}), n \in \nat \bigr \} $ is tight on 
${\xcal }_{1} \times ({\rzecz }^{\nat } \times \nat )$. By the Prokhorov Theorem we can choose a subsequence $({n}_{k}{)}_{k \in \nat }$ such that  $\bigl\{ \lcal aw ({\xi }_{{n}_{k}})$, $k \in \nat \bigr \} $ is weakly convergent on ${\xcal }_{1} \times ({\rzecz }^{\nat } \times \nat )$.
Thus the subsequence ${\bigl( {\xi }_{{n}_{k}}\bigr) }_{k \in \nat }$ satisfies the assumption of Theorem
\ref{T:E.1_[B,H'2010]}. Hence there exists 
a probability space $(\bar{\Omega }, \bar{\fcal }, \bar{\p })$, a family of ${\xcal }_{1} \times ({\rzecz }^{\nat } \times \nat )$-valued random variables $\{ {\bar{\xi }}_{k}, \, k \in \nat  \} $ on 
$(\bar{\Omega }, \bar{\fcal }, \bar{\p })$ and a random variable ${\xi  }_{*}: \bar{\Omega } \to 
{\xcal }_{1} \times ({\rzecz }^{\nat } \times \nat )$ such that 
\begin{itemize}
\item[(i) ] $\lcal aw ({\bar{\xi }}_{k}) = \lcal aw ({\xi }_{{n}_{k}})$ for all $ k \in \nat $;
\item[(ii) ] ${\bar{\xi }}_{k} \to {\xi }_{*}$ in ${\xcal }_{1} \times ({\rzecz }^{\nat } \times \nat )$ a.s. as $k \to \infty $;
\item[(iii) ] ${P}_{1} \circ {\bar{\xi }}_{k} (\bar{\omega }) = {P }_{1} \circ {\xi }_{*}(\bar{\omega })$ for all $\bar{\omega } \in \bar{\Omega }$.
\end{itemize}
Let us put 
$$ \label{E:bar_chi1_k}
{\bar{\chi }}^{1}_{k} := {P}_{1} \circ {\bar{\xi }}^{}_{k},  \qquad k \in \nat .
$$
Notice that $\bigl( {P}_{2}\circ {\bar{\xi }}_{k} {\bigr) }_{k \in \nat }$ is the Skorokhod representation for the sequence  $\bigl( \tilde{f}\circ {\chi }^{2}_{{n}_{k}}, \Phi \circ \tilde{f}\circ {\chi }^{2}_{{n}_{k}} {\bigr) }_{k \in \nat }$.
Let ${P}_{2}\circ {\bar{\xi }}_{k} = ({\eta }^{1}_{k}, {\eta }^{2}_{k})$, where 
${\eta }^{1}_{k}: \bar{\Omega } \to {\rzecz }^{\nat }$ and ${\eta }^{2}_{k}: \bar{\Omega } \to \nat $, 
$k \in \nat $ and let ${P}_{2} \circ {\xi }_{*} = ({\eta }^{1}_{*},{\eta }^{2}_{*})$, where
${\eta }^{1}_{*} :\bar{\Omega } \to {\rzecz }^{\nat } $ and  ${\eta }^{2}_{*}: \bar{\Omega } \to \nat $.
In the same way as in the proof of Lemma 1 in \cite{Jakubowski_1998}, we can prove that
${\eta }^{2}_{k} = \Phi ({\eta }^{1}_{k})$, $\bar{\p }$-a.s., $k \in \nat $.
Since ${\eta }^{2}_{*} < \infty  $ $\bar{\p }$-a.s., we have
$$
   \sup_{k\in \nat } \Phi ({\eta }^{1}_{k}) < \infty  \qquad  \mbox{$\bar{\p }$-a.s.}
$$
Thus for $\bar{\p } $-almost all $\omega \in \bar{\Omega }$ the values ${\eta }^{1}_{k}(\omega )$ belong to the $\sigma $-compact subspace $\bigcup_{m=1}^{\infty } {\tilde{K}}_{m}$ $ = \tilde{f} \bigl( \bigcup_{m=1}^{\infty } {K}_{m} \bigr)$. Since $\tilde{f} $ restricted to $\sigma $-compact subspace is a measurable homeomorphism, we can define 
$$ \label{E:bar_chi2_k}
{\bar{\chi }}^{2}_{k} := {\tilde{f}}^{-1} ({\eta }^{1}_{k}),  \qquad k \in \nat .
$$
Finally ${\bar{\chi }}_{k}$ is defined by 
$$
    {\bar{\chi }}_{k}:= \bigl( {\bar{\chi }}^{1}_{k}, {\bar{\chi }}^{2}_{k} \bigr) ,  \qquad k \in \nat .
$$
This completes the proof. \qed

\bigskip \noindent
In Section \ref{S:Existence} we  use Corollary \ref{C:Skorokhod_J,B,H} for the space 
$$
 {\xcal }_{2}:={\zcal }_{}:=  {L}_{w}^{2}(0,T;V) \cap {L}^{2}(0,T;{H}_{loc}) \cap \dmath ([0,T];U')
  \cap \dmath ([0,T];{H}_{w}).
$$
So, in the following Remark we will discuss the problem of existence of the  countable family of real valued continuous mappings defined on $\zcal $ and separating points of this space.

\bigskip  
\begin{remark} \ \label{R:separating_maps}
\begin{itemize}
\item[(1)] Since  ${L}^{2}(0,T;{H}_{loc})$ and $\dmath ([0,T];U')$ are separable and completely metrizable spaces, we infer that on each of these spaces there exists a  countable family of continuous real valued mappings separating points, see \cite{Badrikian_70}, expos\'{e} 8.
\item[(2)] For the space ${L}^{2}_{w}(0,T;V)$ it is sufficient to put
$$
    {f}_{m}(u):= \int_{0}^{T} \dirilsk{u(t)}{{v}_{n}(t)}{} \, dt \in \rzecz ,
 \qquad u \in {L}^{2}_{}(0,T;V),\quad m \in \nat ,
$$
where $\{ {v}_{m}, m \in \nat  \} $ is a dense subset of ${L}^{2}(0,T;V)$.
Then $({f}_{m}{)}_{m \in \nat }$ is a sequence of continuous real valued mappings separating points of the space ${L}^{2}_{w}(0,T;V)$.
\item[(3)] Let ${H}_{0} \subset H$ be a countable and dense subset of $H$. Then by 
(\ref{E:D([0,T];H_w)_cadlag}) for each $h \in {H}_{0}$ the mapping
 $$
   \dmath ([0,T];{H}_{w}) \ni u \mapsto \ilsk{u(\cdot )}{h}{H} \in \dmath ([0,T];\rzecz ) 
 $$
is continuous.  Since $\dmath ([0,T];\rzecz )$ is a separable complete metric space, there exists a sequence $({g}_{l}{)}_{l \in \nat }$ of real valued continuous functions defined on $\dmath ([0,T];\rzecz )$ separating points of this space. Then the mappings ${f}_{h,l}$, $h \in {H}_{0}$, $l \in \nat $ defined by 
$$
   {f}_{h,l}(u):= {g}_{l} \bigl( \ilsk{u(\cdot )}{h}{H} \bigr) , \qquad 
    u \in \dmath ([0,T];{H}_{w}),
$$ 
form a countable family of continuous mappings on $\dmath ([0,T];{H}_{w})$ separating points of this space.
\end{itemize}
\end{remark}

\section{Appendix C: Some auxilliary results from functional analysis}

\noindent
The following result can be found in Holly and Wiciak, \cite{Holly_Wiciak_1995}.
We recall it together with the proof.

\begin{lemma} \label{L:2_5_Holly_Wiciak} \rm (see Lemma 2.5, p.99 in \cite{Holly_Wiciak_1995}) \it
Consider a separable Banach space $\Phi $ having the following property
\begin{equation} \label{E:2_6_Holly_Wiciak}
  \mbox{there exists a Hilbert space $\hmath $ such that $\Phi \subset \hmath $ continuously. }
\end{equation}
Then there exists a Hilbert space $\bigl( \hcal , (\cdot |\cdot {)}_{\hcal } \bigr) $ such that
$\hcal \subset \Phi$, $\hcal $ is dense in $\Phi $ and the embedding $\hcal \hookrightarrow \Phi $ is compact.
\end{lemma}

\proof
Without loss of generality we can assume that $\dim \Phi = \infty $ and $\Phi $ is dense in $\hmath $.
Since $\Phi $ is separable, there exists a sequence $({\varphi }_{n}{)}_{n \in \nat } \subset \Phi $
linearly dense in $\Phi $.
Since $\Phi $ is dense in $\hmath $ and the embedding $\Phi \hookrightarrow \hmath $ is continuous,
the subspace
 $\lin \{ {\varphi }_{1}, {\varphi }_{2},... \}  $  is dense in $\hmath $.
After the orthonormalisation of $({\varphi }_{n})$ in the Hilbert space
$\bigl( \hmath , (\cdot |\cdot {)}_{\hmath } \bigr) $ we obtain an orthonormal basis $({h }_{n})$ of this space.
 Furthermore, the sequence $({h }_{n})$ is linearly dense in $\Phi $.
 Since the natural embedding $\iota : \Phi \hookrightarrow \hmath $ is continuous, we infer that
$$
  1 = |{h}_{n} {|}_{\hmath } = |\iota ({h}_{n}){|}_{\hmath } \le |\iota | \cdot |{h}_{n} {|}_{\Phi }
$$
and $$ \frac{1}{|{h}_{n} {|}_{\Phi }} \le |\iota | \qquad \mbox{for all } n \in \nat .$$
Let us take ${\eta }_{0} \in (0,1)$ and define inductively a sequence $({\eta }_{n}{)}_{n \in \nat }$ by
$$
   {\eta }_{n} : = \frac{{\eta }_{n-1} +1}{2} , \qquad n=1,2,...
$$
The sequence $({\eta }_{n} )$ is strongly increasing and $\nlim {\eta }_{n} =1$.
Let us define a sequence $({r}_{n}{)}_{n \in \nat }$ by
$$
   {r}_{n} := \frac{1- {\eta }_{n}}{2 |{h}_{n}{|}_{\Phi }} > 0 ,  \qquad n=1,2,...
$$
Obviously $\nlim {r}_{n} =0$.
Let us consider the set
$$
 \hcal := \Bigl\{ x \in \hmath  : \quad \sum_{n=1}^{\infty }
  \frac{1}{{r}_{n}^{2}} \cdot |(x|{h}_{n}{)}_{\hmath }{|}^{2} < \infty  \Bigr\}
$$
and the Hilbert space ${L}_{\mu }^{2} ({\nat }^{\ast } ,\rzecz )$, where
$\mu : {2}^{{\nat }^{\ast }} \to [0,\infty ]$ is the measure given by the formula
$$
    \mu (M) := \sum_{n \in M } \frac{1}{{r}_{n}^{2}}, \qquad M \subset {\nat }^{\ast } .
$$
The linear operator
$$
   l : {L}_{\mu }^{2} ({\nat }^{\ast }, \rzecz  )  \ni \xi \mapsto
    \sum_{n=1}^{\infty } {\xi }_{n} {h}_{n} \in \hmath
$$
is well defined. Moreover, $l$ is an injection and hence we may introduce the following inner product
$$
 (\cdot |\cdot {)}_{\hcal } := (\cdot |\cdot {)}_{{L}^{2}} \, \underline{\circ } \,  {l}^{-1}
  : \hcal \times \hcal \ni (x,y) \mapsto ({l}^{-1}x|{l}^{-1}y{)}_{{L}^{2}} \in \rzecz  .
$$
Now, $l$ is an isometry onto the pre-Hilbert  space $(\hcal , (\cdot |\cdot {)}_{\hcal }) $
and consequently $\hcal $ is $(\cdot |\cdot {)}_{\hcal }$-complete.
Let us notice that for all $x,y\in \hcal $
$$
   (x|y {)}_{\hcal } = \sum_{n=1}^{\infty }  \frac{1}{{r}_{n}^{2}} \cdot (x|{h}_{n}{)}_{\hmath } (y|{h}_{n}{)}_{\hmath },
   \qquad |x{|}_{\hcal }^{2} = \sum_{n=1}^{\infty } \frac{1}{{r}_{n}^{2}} \cdot |(x|{h}_{n}{)}_{\hmath }{|}^{2}
$$
We will show that $\hcal \subset \Phi $ continuously.
Indeed, let $x\in \hcal $, $|x{|}_{\hcal } \le 1 $. Then for each $i \in \nat $
$$
 |(x|{h}_{i}){h}_{i}{|}_{\Phi } = |(x|{h}_{i})| \cdot |{h}_{i}{|}_{\Phi }
 \le {r}_{i} |{h}_{i}{|}_{\Phi } = \frac{1-{\eta }_{i-1}}{2|{h}_{i}{|}_{\Phi } } |{h}_{i}{|}_{\Phi }
 = \frac{1-{\eta }_{i-1}}{2} \! = {\eta }_{i} - {\eta }_{i-1} .
$$
Thus, for any $k,n \in \nat $, $k< n$, we have the following estimate
$$
 \Bigl| \sum_{i=k+1}^{n}(x|{h}_{i}){h}_{i}{\Bigr| }_{\Phi }
 \le \sum_{i=k+1}^{n} ({\eta }_{i} - {\eta }_{i-1}) = {\eta }_{n} - {\eta }_{k}.
$$
Since in particular, the sequence $\bigl( {s}_{n}:= \sum_{i=1}^{n}(x|{h}_{i}){h}_{i} \bigr) $
is Cauchy in the Banach space $(\Phi , |\cdot {|}_{\Phi })$,
there exists $\varphi \in \Phi $ such that $\nlim | {s}_{n}-\varphi  {|}_{\Phi }=0$.
On the other hand, 

\noindent
${s}_{n} = \sum_{i=1}^{n}(x|{h}_{i}){h}_{i} \to x $ in $\hmath $.
Thus by the uniqueness of the limit $\varphi =x \in \Phi $ and
$$
  \sum_{i=1}^{n}(x|{h}_{i}){h}_{i} \to x \quad \mbox{ in $\Phi $ } .
$$
Moreover,
$$
 |x{|}_{\Phi } \stackrel{\infty \leftarrow n}{\longleftarrow } |{s}_{n}{|}_{\Phi }
 \le {\eta }_{n} - {\eta }_{0}  \stackrel{n \to \infty }{\longrightarrow } 1 - {\eta }_{0} .
$$
Thus $\hcal \subset \Phi $ continuously (with the norm of the embedding not  exceeding $1 - {\eta }_{0}$).
We will show that  the embedding $j: \hcal \hookrightarrow \Phi $ is compact.
It is sufficient to prove that the ball $Z:= \{ x \in \hcal :|x{|}_{\hcal } \le 1 \} $
is relatively compact in $(\Phi , |\cdot {|}_{\Phi }) $.
According to the Hausdorff Theorem  it is sufficient to find (for every  fixed $\eps $)
an $\eps $-net of the set $j(Z)$.

\bigskip  \noindent
Since $\lim_{n \to \infty }{\eta }_{n}=1$, there exists $n \in \nat $ such that $1 - {\eta }_{n} \le \frac{\eps }{2}$.
The linear operator
$$
  {S}_{n} : \hcal \ni x \mapsto \sum_{i=1}^{n}(x|{h}_{i}){h}_{i}  \in \Phi
$$
being finite-dimensional is compact. Therefore ${S}_{n}(Z)$ is relatively compact in $(\Phi , |\cdot {|}_{\Phi }) $
and consequently there is a finite subset $F \subset \Phi $ such that
${S}_{n}(Z) \subset \bigcup_{\varphi \in Z } {\ball }_{\Phi }(\varphi , \frac{\eps }{2})$.

\bigskip  \noindent
We will show that the set $F$ is the $\eps $-net for $j(Z)$. Indeed, let $x \in Z$.
Then ${S}_{N}(x) \to x $ in $(\Phi , |\cdot {|}_{\Phi }) $ and
$$
 |x - {S}_{n}(x){|}_{\Phi } \stackrel{\infty \leftarrow N}{\longleftarrow } |{S}_{N}(x)- {S}_{n}(x){|}_{\Phi }
 \le {\eta }_{N} - {\eta }_{n}  \stackrel{N \to \infty }{\longrightarrow } 1 - {\eta }_{n} \le \frac{\eps }{2} .
$$
On the other hand, ${S}_{n}(x) \in {S}_{n}(Z)$, so, there is $\varphi \in F $ such that
${S}_{n}(x) \in {\ball }_{\Phi }(\varphi , \frac{\eps }{2})$.
Finally,
$$
 |x - \varphi {|}_{\Phi } \le |x - {S}_{n}(x){|}_{\Phi } + |{S}_{n}(x) - \varphi {|}_{\Phi }
 \le \frac{\eps }{2} + \frac{\eps }{2} = \eps ,
$$
i.e. $x \in {\ball }_{\Phi }(\varphi ,\eps )$. Thus
$$
 Z \subset \bigcup_{\varphi \in \Phi } {\ball }_{\Phi }(\varphi ,\eps ) .
$$
The proof is thus complete. \qed

\section{Appendix D: Proof of Lemma \ref{L:B_conv_aux}}  \label{A:Appendix_D}

\proof
Assume first that $\psi \in \vcal $. Then there exists $R>0$ such that $\supp \psi $ is a compact subset of ${\ocal }_{R}$.
Then, using the integration by parts formula, we infer that for every $v ,w \in H $
\begin{eqnarray} 
& & | \dual{B(v,w)}{\psi }{} | = \Bigl| \int_{{\ocal }_{R}} ( v  \cdot \nabla \psi ) w \, dx \Bigr| 
  \nonumber \\
& & \le \norm{u}{{L}^{2}({\ocal }_{R})}{}  \norm{w}{{L}^{2}({\ocal }_{R})}{}
  \norm{\nabla \psi }{{L}^{\infty }({\ocal }_{R})}{}  
  \le C |u{|}_{{H}_{{\ocal }_{R}}}  |w{|}_{{H}_{{\ocal }_{R}}} \norm{\psi }{{V}_{m }}{} .
 \label{E:estimate_B(O_R)_ext}   
\end{eqnarray}
We have
$ B(\un , \un ) - B(u,u) =  B(\un -u , \un ) + B(u,\un -u) $.
Thus, using the estimate (\ref{E:estimate_B(O_R)_ext}) and the H\H{o}lder inequality, we obtain
\begin{eqnarray*}
& & \Bigl| \int_{0}^{t} \dual{ B \bigl( \un (s) ,\un (s) \bigr) }{\psi }{} \, ds 
- \int_{0}^{t} \! \dual{ B \bigl( u(s), u(s)\bigr)  }{\psi }{} \, ds \Bigr|  \\
& &\le \Bigl| \! \int_{0}^{t} \dual{ B \bigl( \un (s) - u(s) ,\un (s) \bigr)  }{\psi }{} \, ds  \Bigr|
  + \Bigl| \int_{0}^{t} \! \dual{ B \bigl( u (s) ,\un (s) - u(s) \bigr)  }{\psi }{} \, ds \Bigr| \\
& &\le C \cdot
\norm{\un - u }{{L}^{2}(0,T;{H}_{{\ocal }_{R}})}{} \bigl( \norm{\un  }{{L}^{2}(0,T;{H}_{{\ocal }_{R}})}{}
+\norm{  u }{{L}^{2}(0,T;{H}_{{\ocal }_{R}})}{} \bigr)
\norm{\psi }{{V}_{m }}{} ,
\end{eqnarray*}
where $C$ stands for some positive constant.
Since $\un \to u $ in ${L}^{2}(0,T;{H}_{loc})$, we infer that for all $ \psi  \in \vcal $
\begin{equation}
    \lim_{n \to \infty } \int_{0}^{t} \dual{ B \bigl( \un (s)\bigr)  }{\psi }{} \, ds
 = \int_{0}^{t} \dual{ B \bigl( u(s)\bigr)  }{\psi }{} \, ds  . \label{E:Appendix_D_B_conv}
\end{equation}
If $\psi \in {V}_{m }$ then for every  $\eps > 0 $ there exists ${\psi }_{\eps } \in \vcal $
such that $\norm{\psi - {\psi }_{\eps }}{{V}_{m }}{} \le \eps $.
Then
\begin{eqnarray*}
& &\bigl| \dual{ B( \un (s )) - B( u(s ) ) }{\psi }{} \bigr| 
  =  \bigl| \dual{ B( \un (s ) ) - B( u(s )) }{
\psi -{\psi }_{\eps }}{} \bigr| \\
& &\quad  + \bigl| \dual{ B(\un (s )) - B( u(s ))  }{{\psi }_{\eps } }{} \bigr|  \\
&  & \le  \bigl( \bigl| B( \un (s ) ){\bigr| }_{{V}_{m }'}
  + \bigl| B( u (s )) {\bigr| }_{{V}_{m }'} \bigr)
  \cdot \norm{\psi -{\psi }_{\eps }}{{V}_{m }}{}  \\
 & &\quad  + \bigl| \dual{ B( \un (s )) - B( u(s ) )  }{{\psi }_{\eps } }{} \bigr| \\
 &  &\le \eps  \bigl( |\un (s ){|}_{H}^{2}
   +  |u (s ){|}_{H}^{2}\bigr)
 + \bigl| \dual{ B( \un (s ) )
   - B( u(s))  }{{\psi }_{\eps } }{} \bigr| .
\end{eqnarray*}
Hence
\begin{eqnarray*}
& &\Bigl| \int_{0}^{t} \dual{ B( \un (s)) - B( u(s )) }{\psi }{}  \, ds \Bigr| \\
&  &\le  \eps \! \int_{0}^{t} \! \bigl( \! |\un (s ){|}_{H}^{2} +  |u (s ){|}_{H}^{2} \! \bigr)  d s 
 + \Bigl| \! \int_{0}^{t} \! \dual{ B( \un (s )) 
 - B( u(s ))  }{{\psi }_{\eps } }{}  ds  \Bigr| \\
& &\le \eps \cdot \bigl( \sup_{n\ge 1}\norm{\un }{{L}^{2}(0,T;H)}{2} +\norm{u }{{L}^{2}(0,T;H)}{2} \bigr) 
+ \bigl| \int_{0}^{t}\dual{ B( \un (s )) - B( u(s ))  }{{\psi }_{\eps } }{}  ds \bigr|
.
\end{eqnarray*}
Passing to the upper limit as $n \to \infty $, we obtain
$$
 \limsup_{n \to \infty }
\bigl| \int_{0}^{t} \dual{ B( \un (s)) - B( u(s )) }{\psi }{}  \, ds \bigr|
 \le   M \eps ,
$$
where $M:=\sup_{n\ge 1}\norm{\un }{{L}^{2}(0,T;H)}{2} +\norm{u }{{L}^{2}(0,T;H)}{2}<\infty $.
Since $\eps >0$ is arbitrary, we infer that (\ref{E:Appendix_D_B_conv}) holds for all $\psi \in {V}_{m}$.
The proof of the Lemma \ref{L:B_conv_aux} is thus complete. \qed

\bigskip 

\section{Appendix E: Proof of Lemma \ref{L:D(0,T,{H}_{w})_conv}}
 \noindent
To prove that $\un \in \dmath ([0,T];{H}_{w})$, it is sufficient to show that for every $h\in H$ the real-valued functions $\ilsk{\un (\cdot )}{h}{}$ are c\`{a}dl\`{a}g on $[0,T]$, i.e. are right continuous and have left limits at every $t \in [0,T]$. Let us fix $n \in \nat $ and ${t}_{0}\in [0,T]$ and let us assume that $h\in U$.
Since $\un \in \dmath ([0,T];U')$, there exists ${a}_{}\in U'$ such that
\begin{equation} \label{E:conv_t0_U'}
   \lim_{t \to {t}_{0}^{-}} \norm{\un (t) - {a}_{}}{U'}{} =0 .
\end{equation}
In fact, $a \in H$. Indeed,  by assumption (i) $\un ([0,T]) \subset H$ and 
$$
    \sup_{s \in [0,T]} |\un (s) {|}_{H} \le r .
$$
Let $({t}_{k}{)}_{k \in \nat } \subset [0,T]$ be a sequence convergent to ${t}_{0}^{-}$.
Since $|\un ({t}_{k}){|}_{H} \le r $, by the Banach-Alaoglu Theorem there exists a subsequence  convergent weakly in $H$ to some $b \in H$, i.e. there exists $({t}_{{k}_{l}}{)}_{l \in \nat }$ such that 
$ \un ({t}_{{k}_{l}}) \to b $ weakly in $ H $ as  $ l \to \infty $. 
Since the embedding $H \hookrightarrow U'$ is continuous, we infer that 
$$
   \un ({t}_{{k}_{l}}) \to b \qquad \mbox{weakly in }  U'  \mbox{ as }  l \to \infty . 
$$
On the other hand, by (\ref{E:conv_t0_U'})
$$
   \un ({t}_{{k}_{l}}) \to a \qquad \mbox{ in } U' \mbox{ as }  l \to \infty . 
$$
Hence $a=b \in H$.

\bigskip  \noindent
We have
\begin{equation} \label{E:estimate_t0_U'}
  \bigl| \ilsk{\un (t)-a}{h}{H} \bigr|  = \bigl| \dual{\un (t)-a}{h}{} \bigr| 
 \le  {\|\un (t) - a \|}_{U'} \cdot \norm{h}{U}{} .
\end{equation}
By (\ref{E:conv_t0_U'}) and (\ref{E:estimate_t0_U'}) we infer that
$\lim_{t \to {t}_{0}}  \ilsk{\un (t)-a}{h}{H}  =0$.
Now, let $h\in H$ and let $\eps >0$. Since $U$ is dense in $H$, there exists $h_{\eps} \in U $ such that ${|h-{h}_{\eps } |}_{H} \le \eps  $.
We have the following inequalities
\begin{eqnarray*}
& \bigl| \ilsk{{u}_{n} (t)-a}{h}{H} \bigr| & 
\le \bigl| \ilsk{{u}_{n} (t)-a}{h-h_{\eps}}{H} \bigr|   
 + \bigl| \ilsk{{u}_{n} (t)-a}{h_{\eps }}{H} \bigr|  \\
& & \le {|\un (t)-a|}_{H} \, {|h-h_{\eps}|}_{H} 
 + \bigl| \ilsk{{u}_{n} (t)-a}{h_{\eps }}{H} \bigr| \\
& & \le 2 \eps \norm{{u}_{n} }{{L}^{\infty }(0,T;H)}{} 
  + \bigl| \ilsk{{u}_{n} (t)-a}{h_{\eps }}{H} \bigr| \\
& & \le 2 \eps r  + \bigl| \ilsk{{u}_{n} (t)-a}{h_{\eps }}{H} \bigr| .
\end{eqnarray*}
Passing to the upper limit as $t \to {t}_{0}^{-}$, we obtain
$$
  \limsup_{t \to {t}_{0}^{-}} \bigl| \ilsk{{u}_{n} (t)-a}{h}{H} \bigr| \le 2 \eps r
$$
Since $\eps $ was chosen in an arbitrary way, we infer that
$$
  \lim_{t \to {t}_{0}^{-}}  \ilsk{{u}_{n} (t)-a}{h}{H} =0.
$$  
The proof of right continuity of $\un $ is analogous.

\bigskip  \noindent
We claim that 
$$
  {u}_{n} \to u \quad \mbox{in} \quad  \dmath ([0,T];{\ball }_{w}) 
\quad \mbox{as } \quad \ninf ,
$$  
i.e. that for all $h\in H $
$$
   \ilsk{\un }{h}{H}  \to \ilsk{u}{h}{H} 
  \quad \mbox{in} \quad  \dmath ([0,T];\rzecz ) . 
$$
By (ii) and Remark \ref{R:cadlag_conv} there exists a sequence $({\lambda }_{n})\subset{\Lambda }_{T}$ converging to identity  uniformly on $[0,T]$ and such that 
$$
   \un \circ {\lambda }_{n} \to u \quad \mbox{in} \quad U' 
$$
uniformly on $[0,T]$.
We will prove that for all $h\in H $
\begin{equation} \label{E:un_lambda_n_uniform_conv}
   \ilsk{\un \circ {\lambda }_{n}}{h}{H}  \to \ilsk{u}{h}{H} 
  \quad \mbox{in} \quad \rzecz   
\end{equation}
uniformly on $[0,T]$.

\bigskip  \noindent
Indeed, let us first  fix $h \in U$. Then for all  $ s \in [0,T]$ we have 
$$
\bigl| \! \ilsk{\un \! \circ \! {\lambda }_{n}(s)-u(s)}{h}{H} \! \bigr|
 = \bigl| \! \dual{\un \! \circ \! {\lambda }_{n}(s)-u(s)}{h}{} \! \bigr|  
 \le  {\|\un \! \circ \! {\lambda }_{n}(s) - u(s) \|}_{U'}  \norm{h}{U}{} .
$$
By Remark \ref{R:cadlag_conv}
$$
  \sup_{s\in [0,T]}  \bigl| \dual{\un \circ {\lambda }_{n}(s)-u(s)}{h}{} \bigr|  
 \le \sup_{s\in [0,T]} \norm{\un \circ {\lambda }_{n}(s)-u(s)}{U'}{} \cdot \norm{h}{U}{} \to 0  
$$ 
as $\ninf $. 
Moreover, since $U$ is dense in $H$,
the desired convergence holds for all $h \in H$. Indeed, let us fix $h \in H$ and $\eps >0$. There exists $h_{\eps} \in U $ such that ${|{h-h_{\eps } }|}_{H} \le \eps  $. Using (i), we infer that for all $s \in [0,T]$ the following estimates hold
\begin{eqnarray*}
& & \bigl|  \ilsk{\un  \circ  {\lambda }_{n}(s)-u(s)}{h}{H}  \bigr|  \\
& &\le | {\un  \circ   {\lambda }_{n}(s)-u(s)|}_{H}  {|h-h_{\eps}|}_{H} 
+ \bigl|  \ilsk{\un \! \circ \! {\lambda }_{n}(s)-u(s)}{h_{\eps}}{H}  \bigr| \\
& &\le \eps \cdot {\| \un \circ {\lambda }_{n}-u\| }_{{L}^{\infty }(0,T;H)} 
+ \bigl| \ilsk{\un \circ {\lambda }_{n}(s)-u(s)}{h_{\eps}}{H}  \bigr|  \\
& &\le 2\eps \cdot \sup_{n \in \nat } {\| \un \| }_{{L}^{\infty }(0,T;H)} + 
\bigl| \ilsk{\un \circ {\lambda }_{n}(s)-u(s)}{h_{\eps}}{H}  \bigr| \\
& &\le 2\eps r  + \sup_{s\in [0,T]} \bigl| \ilsk{\un \circ {\lambda }_{n}(s)-u(s)}{h_{\eps}}{H} \bigr| .
\end{eqnarray*}
Thus
$$
  \sup_{s\in [0,T]}  \bigl| \ilsk{\un \circ {\lambda }_{n}(s)-u(s)}{h}{H}  \bigr| \le 2\eps r
   +\sup_{s\in [0,T]} \bigl| \ilsk{\un \circ {\lambda }_{n}(s)-u(s)}{h_{\eps}}{H} \bigr| .
$$  
Passing to the upper limit as $\ninf $ we obtain
$$
  \limsup_{\ninf }\sup_{s\in [0,T]}  \bigl| \ilsk{\un \circ {\lambda }_{n}(s)-u(s)}{h}{H}  \bigr|  \le 2r\eps .
$$  
Since $\eps $ was chosen in an arbitrary way, we get
$$
  \lim_{\ninf } \sup_{s\in [0,T]} \bigl| \ilsk{\un \circ {\lambda }_{n}(s)-u(s)}{h}{}  \bigr| =0 . 
$$ 
Since $\dmath ([0,T];{\ball }_{w})$ is a complete metric space, we infer that   
$u \in \dmath ([0,T];{\ball }_{w})$ as well. 
By Remark \ref{R:cadlag_conv} this completes the proof of (\ref{E:un_lambda_n_uniform_conv}) and of  Lemma \ref{L:D(0,T,{H}_{w})_conv}. 
\qed

\bigskip  \noindent
\bf Acknowledgements \rm
The author would like to thank Zdzis\l aw Brze\'{z}niak for very helpful comments.

\end{document}